\documentclass[a4paper,11pt]{article}
\usepackage[a4paper,text={16cm,22.5cm},centering]{geometry}
\usepackage{amsmath,amssymb,mathtools}
\usepackage{mathrsfs,bbm}
\usepackage[utf8]{inputenc}
\usepackage[TS1,T1]{fontenc}
\usepackage{lmodern}
\usepackage{enumitem}
\usepackage{microtype}
\usepackage{empheq}
\usepackage{tikz}
\usepackage{authblk}

\usetikzlibrary{arrows} 
\usetikzlibrary{patterns} 
\usetikzlibrary{decorations.markings} 
\usetikzlibrary{decorations.pathmorphing}

\definecolor{vert}{rgb}{0.1,0.7,0.15}

\tikzset{hatch distance/.store in=\hatchdistance, hatch distance=10pt, hatch thickness/.store in=\hatchthickness, hatch thickness=2pt}

\makeatletter 
\pgfdeclarepatternformonly[\hatchdistance,\hatchthickness]{flexible hatch} {\pgfqpoint{-1pt}{-1pt}} {\pgfqpoint{\hatchdistance}{\hatchdistance}} {\pgfpoint{\hatchdistance-1pt}{\hatchdistance-1pt}} {\pgfsetcolor{\tikz@pattern@color} 
\pgfsetlinewidth{\hatchthickness} \pgfpathmoveto{\pgfqpoint{-1pt}{-1pt}} \pgfpathlineto{\pgfqpoint{\hatchdistance}{\hatchdistance}} \pgfusepath{stroke}} 
\makeatother

\makeatletter
\newcommand{\hathat}[1]{
\begingroup
\let\macc@kerna\z@
\let\macc@kernb\z@
\let\macc@nucleus\@empty
\hat{\mathchoice
{\raisebox{.35ex}{\vphantom{\ensuremath{\displaystyle #1}}}}
{\raisebox{.35ex}{\vphantom{\ensuremath{\textstyle #1}}}}
{\raisebox{.28ex}{\vphantom{\ensuremath{\scriptstyle #1}}}}
{\raisebox{.24ex}{\vphantom{\ensuremath{\scriptscriptstyle #1}}}}
\smash{\hat{#1}}}
\endgroup
}
\makeatother

\rmfamily
\DeclareFontShape{T1}{lmr}{b}{sc}{<->ssub*cmr/bx/sc}{}
\DeclareFontShape{T1}{lmr}{bx}{sc}{<->ssub*cmr/bx/sc}{}

\DeclareMathOperator*{\esup}{ess~sup}
\DeclareMathOperator*{\argmin}{arg~min}

\newtheorem{defi}{Definition}[section]
\newtheorem{prop}[defi]{Proposition} 
\newtheorem{lemm}[defi]{Lemma} 
\newtheorem{theo}[defi]{Theorem} 
\newtheorem{coro}[defi]{Corollary} 
\newtheorem{rem}[defi]{Remark} 
\numberwithin{equation}{section}
\newcommand{\qed}{\hglue 0pt\hfill$\square$\par}

\title{Existence Results for the Time-incremental Elastic Contact Problem with Coulomb Friction in 2D}
\author[$\dagger$]{Patrick \textsc{Ballard}}
\author[$\ddag$]{Flaviana \textsc{Iurlano}}
\affil[$\dagger$]{Sorbonne Université, CNRS, Université de Paris, Institut Jean le Rond d'Alembert}
\affil[$\ddag$]{Sorbonne Université, CNRS, Université de Paris, Laboratoire Jacques-Louis Lions}
\begin{document}
\maketitle

\begin{abstract}
	In this article, the structure of the incremental quasistatic contact problem with Coulomb friction in linear elasticity (Signorini-Coulomb problem) is unraveled and sharp existence results are proved for the most general two-dimensional problem with arbitrary geometry and elasticity modulus tensor.  The problem is reduced to a variational inequality involving a nonlinear operator which handles both elasticity and friction.  This operator is proved to fall into the class of the so-called Leray-Lions operators, so that a result of Brézis can be invoked to solve the variational inequality.  It turns out that one property in the definition of Leray-Lions operators is difficult to check and requires proving a new fine property of the linear elastic Neumann-to-Dirichlet operator.  This fine property is only established in the case of the two-dimensional problem, limiting currently our existence result to that case.  In the case of isotropic elasticity, either homogeneous or heterogeneous, the existence of solutions to the Signorini-Coulomb problem is proved for  arbitrarily large friction coefficient.  In the case of anisotropic elasticity, an example of nonexistence of a solution for large friction coefficient is exhibited and the existence of solutions is proved under an optimal condition for the friction coefficient.
\end{abstract}

\section{Introduction}

Among the basic Physics laws that every first-year university student in Science learns, there are the Hooke law (constitutive equation of a linear spring) and the Coulomb law of dry friction.  The study of the coupling between these two basic laws, that is, the coupling between linear elasticity and dry friction in the case of a continuum, amounts to studying linear elasticity with unilateral contact at the boundary, complemented with the Coulomb law of dry friction.  This coupling should be naturally considered within the elastodynamics framework, as it is common experience that the coupling between elasticity and dry friction can excite vibratory responses, such as brake squeal.  This is actually completely out of reach of current knowledge, as no existence of solution has been proved yet for elastodynamics with \emph{frictionless} unilateral contact condition on the boundary (the so-called dynamic Signorini problem).  
As there are many examples of situations in which no vibrations are observed, it makes sense to perform the analysis of the coupling between linear elasticity and dry friction in the \emph{quasistatic} framework, that is, the case where the acceleration terms can be neglected.  As the Coulomb law of dry friction involves velocities, the quasistatic problem is an evolution problem, called the Signorini-Coulomb problem.  It falls back within the class of so-called \emph{rate-independent} processes \cite{FirstMielke,SecondMielke,Mielke} such as perfect elastoplasticity~\cite{DalMasoPlast} or brittle fracture~\cite{FrancfortMarigo}.  It turns out that the Signorini-Coulomb problem raises huge mathematical difficulties which have not been fully overcome yet and that this problem is less understood nowadays than perfect elastoplasticity and brittle fracture are.  

In this paper, we consider the incremental problem, that is, the problem arising on each time step of a time discretization.  The mathematical structure of that problem is unraveled for the most general \emph{two-dimensional} case, yielding optimal existence results.  The extension of the analysis to higher space dimension is not straightforward and is currently under investigation (see Remark \ref{r:3d}).  The complete understanding of the incremental problem is generally considered as a prerequisite for the proof of existence of a time continuous quasistatic evolution.  In the context of elastic contact problems with Coulomb friction, the passage from discrete time to continuous time does not raise any additional difficulty in the case of small friction coefficients \cite{Andersson}. Quasistatic evolution in the case of arbitrary large friction coefficients will be addressed in a future work.

The precise formulation of the problem, the state of the art and the description of our approach are given respectively in the next three sections.

\subsection{Formulation of the Signorini-Coulomb problem}

The Coulomb law of dry friction is an empirical law, proposed by Coulomb during the eighteenth century on the basis of the analysis of the rigid motion of a plate under the action of a prescribed overall force \(\mathbf{t}\).  It relates the slipping velocity \(\dot{\mathbf{u}}_{t}\in \mathbb{R}^{N}\) to the tangential force \(\mathbf{t}_{t}\in\mathbb{R}^{N}\) and the normal force \(t_{n}\leq 0\).  It reads as:
\begin{align*}
	\dot{\mathbf{u}}_{t} = \mathbf{0} & \quad\implies\quad |\mathbf{t}_{t}| \leq -f\,t_{n},\\
	\dot{\mathbf{u}}_{t} \neq \mathbf{0} & \quad\implies\quad\mathbf{t}_{t} = f\,t_{n} \,\frac{\dot{\mathbf{u}}_{t}}{|\dot{\mathbf{u}}_{t}|},
\end{align*}
where \(|\cdot|\) stands for the Euclidean norm and \(f\geq 0\) is a given dimensionless, material-dependent friction coefficient.  The Coulomb law of dry friction is equivalent to the following more concise formulation:
\[
	\forall\hat{\mathbf{v}}\in\mathbb{R}^N,\qquad \mathbf{t}_{t}\cdot\bigl(\hat{\mathbf{v}}-\dot{\mathbf{u}}_{t}\bigr) - f\,t_{n}\bigl(|\hat{\mathbf{v}}| - |\dot{\mathbf{u}}_{t}|\bigr) \geq 0,
\]
where `\(\cdot\)' stands for the usual Euclidean scalar product.  In the particular case \(f=0\) (no friction), the Coulomb law reduces to \(\mathbf{t}_{t}=\mathbf{0}\).

We are now ready to state the formal Signorini-Coulomb problem, which couples elasticity and dry friction.  Consider a smooth bounded open subset \(\Omega\) of \(\mathbb{R}^N\) (\(N=2\) or \(N=3\)), whose boundary is the union of three disjoint parts \(\partial\Omega = \overline{\Gamma}_{U}\cup\overline{\Gamma}_{T}\cup\overline{\Gamma}_{C}\).  We will prescribe respectively Dirichlet conditions on \(\Gamma_U\), Neumann conditions on \(\Gamma_T\) and contact conditions on \(\Gamma_C\).  The domain \(\Omega\) is the reference configuration of a linearly elastic body.  The displacement is denoted by \(\mathbf{u}:\Omega \rightarrow\mathbb{R}^N\), the (linearized) strain by \(\boldsymbol{\varepsilon}(\mathbf{u}):=(\nabla\mathbf{u}+{}^t\nabla\mathbf{u})/2\), the stress by \(\boldsymbol{\sigma}(\mathbf{u})=\boldsymbol{\Lambda}\,\boldsymbol{\varepsilon}(\mathbf{u})\), where the elastic modulus tensor \(\boldsymbol{\Lambda}\) is assumed to satisfy the usual symmetry condition:
\begin{equation}
	\label{eq:reqLamda1}
	\forall\hat{\boldsymbol{\varepsilon}}_1,\hat{\boldsymbol{\varepsilon}}_2\in \mathbb{M}^{N\times N}_\text{sym},\quad \forall x\in \Omega,\qquad \hat{\boldsymbol{\varepsilon}}_1:\boldsymbol{\Lambda}(x)\,\hat{\boldsymbol{\varepsilon}}_2 = \hat{\boldsymbol{\varepsilon}}_2:\boldsymbol{\Lambda}(x)\,\hat{\boldsymbol{\varepsilon}}_1,
\end{equation}
and the strong ellipticity condition:
\begin{equation}
	\label{eq:reqLambda2}
	\exists \alpha>0,\quad \forall \hat{\boldsymbol{\varepsilon}}\in \mathbb{M}^{N\times N}_\text{sym},\quad \forall x\in \Omega,\qquad \hat{\boldsymbol{\varepsilon}}:\boldsymbol{\Lambda}(x)\,\hat{\boldsymbol{\varepsilon}} \geq \alpha\, \hat{\boldsymbol{\varepsilon}}:\hat{\boldsymbol{\varepsilon}}.  
\end{equation}
Above, `:' stands for the Frobenius scalar product between real matrices $A:B:=\sum_{i,j}A_{ij}B_{ij}$.  The outward unit normal to \(\Omega\) will be denoted by \(\mathbf{n}\) and any vector \(\mathbf{v}:\partial\Omega \rightarrow \mathbb{R}^N\) will be split into normal and tangential parts: \(\mathbf{v} = v_n\mathbf{n} + \mathbf{v}_t\) where the scalar product \(\mathbf{n}\cdot \mathbf{v}_t=0\) vanishes.  The loading consists in a given volume force \(\mathbf{F}:\Omega \rightarrow\mathbb{R}^N\) and a given surface force \(\mathbf{T}:\Gamma_T \rightarrow\mathbb{R}^N\) on \(\Gamma_T\).  The initial gap in the direction \(\mathbf{n}\) between the elastic body \(\Omega\) and a given rigid obstacle is represented by a function \(g:\Gamma_C \rightarrow \mathbb{R}\), while the given friction coefficient between the two objects is denoted by \(f:\Gamma_C \rightarrow \left[0,+\infty\right[\).  The quasistatic Signorini-Coulomb problem consists in finding a displacement \(\mathbf{u}(s):\Omega\to\mathbb{R}^N\) defined for time \(s\in [0,T]\), satisfying a given initial condition and, defining the surface traction $\mathbf{t}:= \boldsymbol\sigma(\mathbf{u})\,\mathbf{n}=t_n\mathbf{n}+\mathbf{t}_t$:
\begin{equation}
	\label{eq:SignoriniCoulombCont}
	\left\{\quad
	\begin{aligned}
		& \text{div}\,\boldsymbol{\sigma}(\mathbf{u}) + \mathbf{F} = \mathbf{0}, \qquad & & \text{ in }\Omega,\\
		& \mathbf{u} = \mathbf{0},\qquad\qquad & & \text{ on }\Gamma_U,\\
		& \boldsymbol\sigma(\mathbf{u)\,\mathbf{n}} = \mathbf{T},\qquad & & \text{ on }\Gamma_T,\\
		& u_n - g \leq 0, \qquad t_n \leq 0, \qquad t_n\,(u_n-g) = 0, \qquad & & \text{ on }\Gamma_C,\\
		& \forall \hat{\mathbf{v}}\in\mathbb{R}^N,\qquad \mathbf{t}_t\cdot\bigl(\hat{\mathbf{v}}-\dot{\mathbf{u}}_t\bigr) - ft_n \bigl(|\hat{\mathbf{v}}|-|\dot{\mathbf{u}}_t|\bigr) \geq 0, \qquad & & \text{ on }\Gamma_C,
	\end{aligned}
	\right.
\end{equation}
where the dot stands for the time derivative and \(|\cdot|\) for the Euclidean norm.

As it is usual in the analysis of such rate-independent quasistatic evolutions, we introduce a time discretization with time step $\Delta s>0$, say \(s_i:=i\,\Delta s\) (\(i=0,1,\ldots\)), and consider the problem raised on one time step by replacing the velocity \(\dot{\mathbf{u}}\) by \((\mathbf{u}_i - \mathbf{u}_{i-1})/\Delta s\).  Denoting by \(\mathbf{w}:=\mathbf{u}_{i-1}\) the displacement at the preceding time step, which is supposed to be given, and dropping the index \(i\), the problem on one time step is now to find a (time independent) displacement \(\mathbf{u}:\Omega \rightarrow\mathbb{R}^N\) satisfying:
\begin{subequations}
	\label{eq:SignoriniCoulombDiscret}
	\begin{empheq}[left={\empheqlbrace\quad}]{align}
		& \text{div}\,\boldsymbol{\sigma}(\mathbf{u}) + \mathbf{F} = \mathbf{0}, \qquad & & \text{ in }\Omega,\\
		& \mathbf{u} = \mathbf{0},\qquad\qquad & & \text{ on }\Gamma_U,\\
		& \boldsymbol\sigma(\mathbf{u})\,\mathbf{n} = \mathbf{T},\qquad & & \text{ on }\Gamma_T,\\
		& u_n - g \leq 0, \qquad t_n \leq 0, \qquad t_n\,(u_n-g) = 0, \qquad & & \text{ on }\Gamma_C,\label{eq:Signorini}\\
		& \forall \hat{\mathbf{v}}\in\mathbb{R}^N,\qquad \mathbf{t}_t\cdot\bigl(\hat{\mathbf{v}}-\mathbf{u}_t\bigr) - ft_n \bigl(|\hat{\mathbf{v}}-\mathbf{w}_t|-|\mathbf{u}_t-\mathbf{w}_t|\bigr) \geq 0, \qquad & & \text{ on }\Gamma_C.
		\label{eq:Frict}
	\end{empheq}
\end{subequations}

\subsection{Historical background}

The mathematical analysis of the time-incremental problem~\eqref{eq:SignoriniCoulombDiscret} was first considered by Duvaut \& Lions in~\cite{DuvautLions}.  They observed that when \(-ft_n\) is replaced with a given \(\tau\geq 0\) in line~\eqref{eq:Frict} of problem~\eqref{eq:SignoriniCoulombDiscret} (Coulomb law), then one gets the following variational problem:
\[
\mathbf{u} = \argmin_{\substack{\mathbf{v},\\
\mathbf{v}=\mathbf{0},\text{ on }\Gamma_U,\\
v_{\rm n}\leq g,\text{ on }\Gamma_C,\\
}} \frac{1}{2}\int_{\Omega} \boldsymbol{\varepsilon(\mathbf{v})}:\boldsymbol{\Lambda}\boldsymbol{\varepsilon}(\mathbf{v}) \,{\rm d}x - \int_{\Omega} \mathbf{F}\cdot\mathbf{v} \,{\rm d}x - \int_{\Gamma_T}\mathbf{T}\cdot\mathbf{v} + \int_{\Gamma_C}\tau\,\bigl|\mathbf{v}_t-\mathbf{w}_t\bigr|,
\]
which can be uniquely solved by the direct method of the Calculus of Variations, under appropriate regularity assumptions on the data.  This remark suggests a fixed point strategy applied to the mapping \(\tau\mapsto -ft_n\).  If the mapping \(\tau \mapsto t_n\) were Lipschitz-continuous, then the mapping \(\tau\mapsto -ft_n\) would be a contraction for small friction coefficient \(f\), yielding a unique fixed point and therefore a unique solution of problem~\eqref{eq:SignoriniCoulombDiscret} in that case.  This hope was dashed as no Lipschitz-continuity turns out to be fulfilled.

The next progress came from Jarušek's PhD thesis~\cite{JarusekPhD1,JarusekPhD2} who developed an original idea of Nečas under his supervision.  He was able to run successfully the fixed-point strategy by applying Tikhonov's theorem in a Hilbert space endowed with the weak topology.  The compactness needed to apply this theorem was obtained by requiring additional regularity on the data and proving additional regularity on the solution by use of local rectification together with the shift technique, involving very technical arguments.  In this way, the first existence result for problem~\eqref{eq:SignoriniCoulombDiscret} was proved, provided that the data have little more regularity than usually required, and more importantly, that the friction coefficient is small enough.

An alternative strategy of proof was later designed by Eck and Jarušek in~\cite{EckJarusek,BookJiri}.  They consider a penalty regularization of problem~\eqref{eq:SignoriniCoulombDiscret}, allowing for a small penetration of the body into  the obstacle.  In practice, it amounts to replacing line~\eqref{eq:Signorini} in problem~\eqref{eq:SignoriniCoulombDiscret} by:
\begin{equation*}
	t_n = - \frac{1}{\varepsilon} \bigl\langle u_n - g \bigr\rangle^+,
\end{equation*}
where \(\langle x \rangle^+ := \max\{0,x\}\) stands for the positive part and \(\varepsilon >0\) is the small regularization parameter.  They show that this regularized problem~\eqref{eq:SignoriniCoulombDiscret} can be easily solved by combining the direct method of the Calculus of Variations with Schauder's fixed point theorem, whatever is \(\varepsilon>0\) and for possibly arbitrarily large friction coefficient.  Then, a limit \(\varepsilon \rightarrow 0+\) has to be taken, which involves the product of two weakly convergent sequences.  Some compactness is sought and the local rectification together with the shift technique is therefore used again.  With this approach also, the existence of a solution for problem~\eqref{eq:SignoriniCoulombDiscret} is proved, provided that the data have little more regularity than usually required, and more importantly, that the friction coefficient is small enough.  Here, `small enough' means smaller than a constant which depends only on the elastic modulus tensor.  Explicit optimal (with respect to their strategy of proof) values for that constant in terms of the Poisson ratio in the case of elastic isotropy are provided in~\cite{EckJarusek} for the cases \(N=2\) and \(N=3\).  This result is still the best known up to now about the existence of solution for problem~\eqref{eq:SignoriniCoulombDiscret}.   We underline that it requires that the friction coefficient is smaller than a finite value, in both cases of isotropic and anisotropic elasticity.

This result was later exploited by Andersson in~\cite{Andersson} who proved the existence of solution for the Signorini-Coulomb evolution problem~\eqref{eq:SignoriniCoulombCont} with continuous time.  His proof relies on the estimates of Eck and Jarušek obtained by use of the local rectification together with the shift technique, and  therefore suffers the same limitations: the data have slightly more regularity than usually required, and more importantly, the friction coefficient is small enough.  Here we emphasize that the only results up to now that escape a smallness condition for the friction coefficient are:
\begin{itemize}
	\item the penalty regularization of problem~\eqref{eq:SignoriniCoulombDiscret} considered by Eck \& Jarušek which can be solved for arbitrarily large friction coefficient, thanks to the direct method of the Calculus of Variations combined with Schauder's fixed point theorem (the same is true of its continuous time counterpart in~\cite{Andersson}),
	\item the finite-dimensional counterpart of problem~\eqref{eq:SignoriniCoulombDiscret} (Galerkin approximation) which can be solved for arbitrarily large friction coefficient, thanks to the direct method of the Calculus of Variations combined with Brouwer's fixed point theorem.
\end{itemize}

In this paper, we prove the new result that the two-dimensional time-incremental problem~\eqref{eq:SignoriniCoulombDiscret} can be solved for arbitrarily large friction coefficient in the case of isotropic elasticity.  We also prove that a critical friction coefficient, only depending on the elastic modulus tensor, appears in the case of anisotropic elasticity.  For friction coefficients smaller than this critical value, the existence of solutions is proved. 
We also exhibit an example of nonexistence of solutions for friction coefficients larger than this critical value, in the case of the simpler steady sliding problem.  Hence, our new strategy of proof seems to yield optimal existence results for the two-dimensional time-incremental problem~\eqref{eq:SignoriniCoulombDiscret}.

Incidentally, our proof of nonexistence of solution for the steady problem in the case of anisotropic elasticity with large friction reveals the possible \emph{non-solvability} of (steady) frictional contact problems in the case of anisotropic elasticity, which seems to have remained unsuspected in the engineering literature.  It is worth mentioning though that the engineering literature contains many examples of non-solvability (see \cite{Klarbring,AnderssonKlarbring}, for example) but they seem to all concern the rate-problem (the problem of finding the velocity in a given configuration of the system) in a finite-dimensional setting. In the finite-dimensional setting, the time-incremental problem is always solvable, thanks to Brouwer's fixed point theorem.

\subsection{Strategy of proof and organization of the article}

We consider the time-incremental problem~\eqref{eq:SignoriniCoulombDiscret}.  We choose to take as the principal unknown, the normal component \(t_n\) of the restriction of the surface force \(\mathbf{t}\) to the contact part \(\Gamma_C\) of the boundary.  Note that if \(t_n=\tau\leq 0\) is known, then we can uniquely solve the following minimum problem:
\[
\mathbf{u} = \argmin_{\substack{\mathbf{v},\\
\mathbf{v}=\mathbf{0},\text{ on }\Gamma_U\\
}} \frac{1}{2}\int_{\Omega} \boldsymbol{\varepsilon(\mathbf{v})}:\boldsymbol{\Lambda}\boldsymbol{\varepsilon}(\mathbf{v}) \,{\rm d}x - \int_{\Omega} \mathbf{F}\cdot\mathbf{v} \,{\rm d}x - \int_{\Gamma_T}\mathbf{T}\cdot\mathbf{v} - \int_{\Gamma_C} \tau\,v_n - \int_{\Gamma_C}f\tau\,\bigl|\mathbf{v}_t-\mathbf{w}_t\bigr|,
\]
and find a displacement field \(\mathbf{u}\) in \(\Omega\) satisfying \eqref{eq:SignoriniCoulombDiscret}, except for the fourth line (contact conditions).

Considering the normal component \(u_n\) of the trace of that minimizer \(\mathbf{u}\) on \(\Gamma_C\), we have a nonlinear operator \(A: \tau \mapsto u_n\), which is well-defined for all nonpositive \(\tau\leq 0\).  Then the time-incremental problem~\eqref{eq:SignoriniCoulombDiscret} can be equivalently formulated as a formal variational inequality involving the nonlinear operator \(A\): find \(t\leq 0\) on \(\Gamma_C\) such that:
\begin{equation}
	\label{eq:formalVI}
	\forall \hat{t}\leq 0,\qquad \int_{\Gamma_C} \bigl(At-g) (\hat{t}-t) \geq 0.
\end{equation}
The nonlinear operator \(A\) describes both elasticity and dry friction.  However, it turns out that this operator is monotone only in the case of no friction \(f=0\).  We are therefore led to seek more general conditions on the nonlinear operator \(A\) that ensure the solvability of the variational inequality~\eqref{eq:formalVI}.  

This was precisely the aim of the article \cite{BrezisIeq} of Brézis, who sought weak conditions on an abstract nonlinear operator \(A\) ensuring that a sequence of Galerkin approximations for the solution to~\eqref{eq:formalVI} can be built, based on Brouwer's theorem, and that a limit to a solution of~\eqref{eq:formalVI} can be taken.  Brézis was able to identify a new class of operators that he called \emph{pseudomonotone}, for which the variational inequality~\eqref{eq:formalVI} can be solved.  The precise definition of pseudomonotone operators is recalled in Appendix~A.  The class of pseudomonotone operators contains a smaller class previously introduced by Leray and Lions \cite{LerayLions} to solve a family of nonlinear elliptic boundary value problems.  In a nutshell, the Leray-Lions operators \(A\) are basically nonlinear operators which have the form:
\begin{equation}
	\label{eq:genFormLL}
	At = \mathscr{A}(t,t),
\end{equation}
where the operator \(\mathscr{A}\) has a monotonicity property with respect to the second variable when the first variable is frozen, satisfies suitably weak continuity properties separately in the two variables, and a condition concerning the passage to the limit in the product of two weakly convergent sequences.  The precise definition of Leray-Lions operators is recalled in Appendix~A and involves four properties \textit{(i)}, \textit{(ii)}, \textit{(iii)} and \textit{(iv)}.  Brézis results \cite{BrezisIeq} made it possible to solve variational inequalities based on Leray-Lions operators.

Coming back to the formal variational inequality associated with the time-incremental problem~\eqref{eq:SignoriniCoulombDiscret}, we see that the operator \(A\) has naturally the general form~\eqref{eq:genFormLL}: it suffices to consider the energy functional:
\[
E_{t,\tau}(v) := \frac{1}{2}\int_{\Omega} \boldsymbol{\varepsilon(\mathbf{v})}:\boldsymbol{\Lambda}\boldsymbol{\varepsilon}(\mathbf{v}) \,{\rm d}x - \int_{\Omega} \mathbf{F}\cdot\mathbf{v} \,{\rm d}x - \int_{\Gamma_T}\mathbf{T}\cdot\mathbf{v} - \int_{\Gamma_C} \tau\,v_n - \int_{\Gamma_C}ft\,\bigl|\mathbf{v}_t-\mathbf{w}_t\bigr|,
\]
and to define \(\mathscr{A}(t,\tau)\) as the normal component of the trace on \(\Gamma_C\) of the unique minimizer of \(E_{t,\tau}\) on the space of those \(\mathbf{v}\) such that \(\mathbf{v}_{|\Gamma_U}=\mathbf{0}\).  Hence, \(At:= \mathscr{A}(t,t)\) and by the results of Brézis surveyed in Appendix~A, we are brought back to prove that the operator \(A\) is Leray-Lions and coercive to prove the existence of a solution for the time-incremental problem~\eqref{eq:SignoriniCoulombDiscret}.  It turns out that it is very easy to prove that \(A\) is unconditionally coercive in any dimension \(N\) and that properties \textit{(i)}, \textit{(ii)}, \textit{(iii)} in the definition of Leray-Lions operator (Definition~\ref{thm:defLerayLions} in Appendix~A) are unconditionally fulfilled in any dimension.

The challenging part is to prove that property \textit{(iv)} in the definition of Leray-Lions operator is also fulfilled.  This requires to pass to the limit in a product of two weakly converging sequences, and therefore to prove somehow a compactness result.  We prove such a compactness result, first in the restrictive case of two-dimensional (2D), that is $N=2$, isotropic elasticity, based on a new fine property of the elastic Neumann-to-Dirichlet operator.  As a consequence we obtain the solvability of the 2D time-incremental problem~\eqref{eq:SignoriniCoulombDiscret} in the case of isotropic elasticity without any smallness condition on the friction coefficient.

The handling of anisotropy raises new challenges as it turns out that the operator \(A\), as defined in the isotropic case, needs no longer to be Leray-Lions in the case of anisotropic elasticity.  The definition of the operator \(A\) has to be slightly adapted, so that it can be proved to be Leray-Lions even in the anisotropic case.  But this adaptation requires a critical condition on the friction coefficient (not needed in the isotropic case).  We verify that this condition is optimal by proving a non-existence result for larger friction coefficients.

The first use of pseudomonotone operators in the sense of Brézis in the context of frictional contact problems seems to be \cite{BallardIurlano}.  There, the steady motion of a 2d elastic half-space is analyzed, under frictional contact with a given obstacle.   We underline that in the present paper, our invocation of Leray-Lions operators is innovative and is made in a completely different context from the one (a class of nonlinear elliptic boundary value problems) which had motivated the original definition by Leray and Lions \cite{LerayLions}.

The article will be split into two parts: the case of isotropic elasticity will be extensively discussed in Section~\ref{sec:isotropy} and the generalization to anisotropic elasticity will be discussed in Section~\ref{sec:anisotropy}.  Section~\ref{sec:isotropy} will consist in three subsections.  In Subsection~\ref{sec:regData}, we state the regularity assumptions on the data and give easy preliminary lemmas.  Subsection~\ref{sec:existProofIsotropy} gathers the existence proof for the 2D isotropic elastic case.  The proof of the required fine property of the elastic Neumann-to-Dirichlet operator, which is a more technical matter, is postponed to Subsection~\ref{sec:fineProperty}.  The discussion of anisotropic elasticity in Section~\ref{sec:anisotropy} is split into three subsections.  In Subsection~\ref{sec:nonexistence}, an example of frictional contact problem with no solution for large friction coefficient is presented.  Subsection~\ref{sec:existAnisotrop} contains the discussion about the generalization of the solvability result previously obtained in isotropic elasticity to the case of anisotropy.  The detailed proof of the structure of the Neumann-to-Dirichlet operator for the homogeneous anisotropic elastic half-space, which plays a central role in the analysis, is postponed to Subsection~\ref{sec:N2Danisotrop}.   Appendix~A recalls general results on pseudomonotone and Leray-Lions operators.

\section{The case of isotropic elasticity}

\label{sec:isotropy}

The framework of this section is that of isotropic elasticity.  We rigorously formulate the time-incremental contact problem with Coulomb friction condition and prove the existence of at least one solution in the 2D case.

\subsection{Data regularity and trace spaces}

\label{sec:regData}

In all the sequel, the bounded open set \(\Omega\subset \mathbb{R}^N\) will be supposed connected and of class \(C^{1,1}\), which entails in particular that \(\mathbf{n}\in W^{1,\infty}(\partial\Omega;\mathbb{R}^N)\).  We consider three nonintersecting open subsets \(\Gamma_U\), \(\Gamma_T\) and \(\Gamma_C\) of the boundary \(\partial \Omega\) of \(\Omega\), such that \(\partial\Omega=\overline{\Gamma}_U\cup\overline{\Gamma}_T\cup\overline{\Gamma}_C\) (where \(\overline{\Gamma}\) stands for the closure of \(\Gamma\) in \(\partial\Omega\)).  For the sake of simplicity, we make the following additional hypotheses.
\begin{itemize}
	\item The subset \(\Gamma_U\) has positive surface measure \(\mathcal{H}^{N-1}\).  As a homogeneous Dirichlet condition will be prescribed on \(\Gamma_U\), this hypothesis is made for convenience, to obtain some coercivity.  This hypothesis is not essential and could be dropped at the price of additional complexity of the presentation.
	\item The subset \(\Gamma_C\subset \partial\Omega\) is of class \(C^{0,1}\) (that is, Lipschitz) and \(\text{dist}(\Gamma_U,\Gamma_C)>0\).  Again, this hypothesis is made for convenience and is probably not essential.  It ensures that the space of traces on \(\Gamma_C\) of functions in \(H^1(\Omega)\) that vanish on \(\Gamma_U\):
	\[
	H := \Bigl\{ u_{|\Gamma_C} \bigm| u\in H^1(\Omega) \; \text{ and }\;u=0\;\text{ on }\Gamma_U\Bigr\}
	\]
	is exactly \(H^{1/2}(\Gamma_C)\).  Indeed, the space of traces of functions in \(H^1(\Omega)\) is \(H^{1/2}(\partial\Omega)\) and restrictions to \(\Gamma_C\) of functions in \(H^{1/2}(\partial\Omega)\) are in \(H^{1/2}(\Gamma_C)\).  Vice versa, when \(\text{dist}(\Gamma_U,\Gamma_C)>0\), it is always possible to find an extension to \(\partial\Omega\) of a given function in \(H^{1/2}(\Gamma_C)\), so that this extension vanishes identically in \(\Gamma_U\).  When this hypothesis is not made, the above trace space can be strictly smaller than \(H^{1/2}(\Gamma_C)\) (see \cite[Section~11]{LionsMagenes}).
\end{itemize}

Considering an arbitrary open \(C^1\) manifold \(\omega\) with dimension \(d\), the spaces \(L^2(\omega)\), \(L^\infty(\omega)\), \(H^1(\omega)\), \(W^{1,\infty}(\omega)\) and \(H^{1/2}(\omega)\) are defined on the basis of their norms:
\begin{align}
	\|v\|_{L^2(\omega)} & := \left(\int_{\omega} |v|^2\,{\rm d}x\right)^{1/2},\nonumber\\
	\|v\|_{L^\infty(\omega)} & := \esup_{x\in\omega} |v(x)|,\nonumber\\
	\|v\|_{H^1(\omega)} & := \left(\|v\|_{L^2(\omega)}^2 + \|\nabla v\|_{L^2(\omega)}^2\right)^{1/2},\nonumber\\
	\|v\|_{W^{1,\infty}(\omega)} & := \max\left\{\|v\|_{L^\infty(\omega)},\|\nabla v\|_{L^\infty(\omega)}\right\},\nonumber\\
	\|v\|_{H^{1/2}(\omega)} & := \left(\int_{\omega} |v|^2\,{\rm d}x+\int_{\omega\times\omega}\frac{|v(x)-v(y)|^2}{|x-y|^{d+1}}\,{\rm d}x\,{\rm d}y\right)^{1/2}.\label{eq:defH12}
\end{align}
From the definition of these norms, we have the inequalities:
\begin{equation}
	\label{eq:multLischitz}
	\Bigl\| |v|\Bigr\|_{H^{1/2}(\omega)} \leq \|v\|_{H^{1/2}(\omega)},\qquad \|fv\|_{H^{1/2}(\omega)} \leq C\|f\|_{W^{1,\infty}(\omega)}\|v\|_{H^{1/2}(\omega)},
\end{equation}
for some constant \(C\) depending only on \(\omega\).  We define the space \(H^{-1/2}(\omega)\) as the dual space of \(H^{1/2}(\omega)\) and it will be systematically endowed with the dual norm. Recalling that the \(C^\infty\) functions with compact support are dense in \(H^{1/2}(\omega)\), the space \(H^{-1/2}(\omega)\) is a space of distributions on \(\omega\).   Here and henceforth $\langle\cdot,\cdot\rangle$ will denote the duality pairing.

Given \(\boldsymbol\Lambda\in L^\infty(\Omega)\) satisfying requirements~\eqref{eq:reqLamda1} and \eqref{eq:reqLambda2}, we have the following standard result on the continuity of the Dirichlet operator, based on Korn's inequality.  
\begin{prop}
	\label{thm:vectorTrace}
	Let $J_\Lambda:H^{1/2}(\Gamma_C;\mathbb{R}^N)\to \{\mathbf v\in H^1(\Omega,\mathbb{R}^N):\mathbf{v}_{|\Gamma_U}=\mathbf{0}\}$ be defined as	
	\[
	J_{\Lambda}(\mathbf{u}) := \argmin_{\substack{\mathbf{v}\in H^1(\Omega;\mathbb{R}^N),\\
	\mathbf{v}_{|\Gamma_U}=\mathbf{0},\\
	\mathbf{v}_{|\Gamma_C}=\mathbf{u}\\
	}} \int_{\Omega} \boldsymbol\varepsilon(\mathbf{v}):\boldsymbol\Lambda\boldsymbol\varepsilon(\mathbf{v})\,{\rm d}x, \qquad\text{for } \mathbf{u}\in H^{1/2}(\Gamma_C;\mathbb{R}^N).
	\]
	Then, \(\| J_\Lambda(\mathbf{u})\|_{H^1(\Omega,\mathbb{R}^N)}\) is a norm on \(H^{1/2}(\Gamma_C;\mathbb{R}^N)\) which is equivalent to that of \(H^{1/2}(\Gamma_C;\mathbb{R}^N)\).
	
	In an analogous way, let $J'_\Lambda:H^{-1/2}(\Gamma_C;\mathbb{R}^N)\to \{\mathbf v\in H^1(\Omega,\mathbb{R}^N):\mathbf{v}_{|\Gamma_U}=\mathbf{0}\}$ be defined as
	\[
	J_\Lambda'(\mathbf{t}) := \argmin_{\substack{\mathbf{v}\in H^1(\Omega;\mathbb{R}^N)\\
	\mathbf{v}_{|\Gamma_U}=0}}\biggl\{\frac{1}{2}\int_{\Omega}\boldsymbol\varepsilon(\mathbf{v}):\boldsymbol\Lambda\boldsymbol\varepsilon( \mathbf{v})\,{\rm d}x - \bigl\langle \mathbf{t},\mathbf{v}_{|\Gamma_C}\bigr\rangle\biggr\}, \qquad\text{for }\mathbf{t}\in H^{-1/2}(\Gamma_C;\mathbb{R}^N).
	\]
	Then,
	\[
	\forall \mathbf{v}\in H^1(\Omega;\mathbb{R}^N), \quad \text{with }\mathbf{v}_{|\Gamma_U}=0,\qquad  \int_\Omega \boldsymbol\varepsilon(\mathbf{J'_{\Lambda}(\mathbf t)}):\boldsymbol\Lambda\boldsymbol\varepsilon( \mathbf{v})\,{\rm d}x = \bigl\langle \mathbf{t},\mathbf{v}_{|\Gamma_C} \bigr\rangle
	\]
	and \(\mathbf{t}\mapsto\|J'_\Lambda(\mathbf{t})\|_{H^1(\Omega,\mathbb{R}^N)}\) is a norm on $H^{-1/2}(\Gamma_C;\mathbb{R}^N)$ which is equivalent to the norm of $H^{-1/2}(\Gamma_C;\mathbb{R}^N)$.
\end{prop}

\noindent\textbf{Proof.} As we have made the assumption that $d(\Gamma_C,\Gamma_U)>0$, it follows from the definition~\eqref{eq:defH12} of the $H^{1/2}$-norm that the extension by $0$ on $\Gamma_U$ of any function in $H^{1/2}(\Gamma_C)$ is in $H^{1/2}(\Gamma_C\cup \Gamma_U)$.  We recall that outside the hypothesis $d(\Gamma_C,\Gamma_U)>0$ (which has been made for conveniency), things are a little bit more intricated (see \cite[Section~11]{LionsMagenes}).  The first minimum problem in the statement of Proposition~\ref{thm:vectorTrace} is a mixed boundary value problem for the elasticity operator.  Standard results in \cite{LionsMagenes} for such mixed problems state that the solution mapping $J_\Lambda$ is continuous from \(H^{1/2}(\Gamma_C;\mathbb{R}^N)\) into \(H^1(\Omega,\mathbb{R}^N)\).  Combining this result with the continuity of the trace operators, we obtain the first claimed quivalence of norms.  The second one follows.\qed

\bigskip
Recalling that \(\Omega\) has been assumed connected and of class \(C^{1,1}\), we can combine Proposition~\ref{thm:vectorTrace} with estimate~\eqref{eq:multLischitz} to yield:

\begin{prop}
	\label{thm:dualTrace}
	Let \(\mathbf{u}\in H^{1/2}(\Gamma_C;\mathbb{R}^N)\).  Then, \(u_n:=\mathbf{u}\cdot\mathbf{n}\in H^{1/2}(\Gamma_C)\) and the mappings:
	\[
	\left\{
	\begin{array}{rcl}
		H^{1/2}(\Gamma_C;\mathbb{R}^N) & \rightarrow & H^{1/2}(\Gamma_C)\\
		\mathbf{u} & \mapsto & \mathbf{u}\cdot\mathbf{n}
	\end{array}
	\right.
	\hspace*{2cm}
	\left\{
	\begin{array}{rcl}
		H^{1/2}(\Gamma_C) & \rightarrow & H^{1/2}(\Gamma_C;\mathbb{R}^N)\\
		u & \mapsto & u\,\mathbf{n}
	\end{array}
	\right.
	\]
	are continuous.  On the dual side, given \(\mathbf{t}\in H^{-1/2}(\Gamma_C;\mathbb{R}^N)\), we define \(t_n\in H^{-1/2}(\Gamma_C)\) by the formula:
	\[
	\forall v\in H^{1/2}(\Gamma_C),\qquad \langle t_n , v\rangle := \langle \mathbf{t},v\mathbf{n}\rangle,
	\]
	and the mapping \(\mathbf{t}\mapsto t_n\) is continuous from \(H^{-1/2}(\Gamma_C;\mathbb{R}^N)\) into \(H^{-1/2}(\Gamma_C)\).  In particular, there exists \(C>0\) such that:
	\[
	\forall \mathbf{t}\in H^{-1/2}(\Gamma_C;\mathbb{R}^N), \qquad \|t_n\|_{H^{-1/2}(\Gamma_C)} \leq C \bigl\| J'_\Lambda(\mathbf{t}) \bigr\|_{H^1(\Omega,\mathbb{R}^N)}.
	\]
\end{prop}

\subsection{Existence in the case of isotropic elasticity}

\label{sec:existProofIsotropy}

Our strategy for proving the existence of a solution to the Signorini-Coulomb problem on one time step is to express this problem in terms of a variational inequality applying to \(t_n\in H^{-1/2}(\Gamma_C)\), the unknown normal part of the surface traction on \(\Gamma_C\).  This variational inequality is based on a nonlinear operator \(A\) whose role is to handle both elasticity in the body and dry friction on the boundary.

More precisely, we will assume the gap with the obstacle to be \(g\in H^{1/2}(\Gamma_C)\), and we will set:
\begin{equation}
	\label{eq:defK}
	K:= \Bigl\{\hat{t}\in H^{-1/2}(\Gamma_C)\;\bigm|\; \hat{t} \leq 0\Bigr\}.
\end{equation}
Above and henceforth, a distribution is said \emph{nonpositive} (notation \(\hat{t}\leq 0\)) whenever it takes nonpositive values on all nonnegative test functions.  A nonpositive distribution is classically a nonpositive measure.  We will then put the Signorini-Coulomb problem on one time step under the form of finding \(t_n\in K\) such that:
\[
\forall \hat{t}\in K,\qquad \bigl\langle At_n -g \,,\, \hat{t}-t_n \bigr\rangle \geq 0,
\]
for an appropriate nonlinear operator \(A:K \rightarrow H^{1/2}(\Gamma_C)\).  Then, it will be sufficient to prove that \(A\) is pseudomonotone in the sense of Brézis and coercive (see Appendix~A) to yield the existence of a solution for this variational inequality.  Actually, the operator will be proved to belong to a subclass of pseudomonotone operators, the so-called subclass of Leray-Lions operators (see Appendix~A).

\begin{defi}
	\label{defu}
	\label{thm:defAscr}
	We assume that \(\Omega\subset \mathbb{R}^N\) is open, bounded, connected and of class \(C^{1,1}\), \(\boldsymbol\Lambda\in L^\infty(\Omega)\), \(\mathbf{F}\in L^2(\Omega;\mathbb{R}^N)\), \(\mathbf{T}\in L^2(\Gamma_C;\mathbb{R}^N)\), \(f\in W^{1,\infty}(\Gamma_C;\left[0,+\infty\right[)\) and \(\mathbf{w}_t\in H^{1/2}(\Gamma_C;\mathbb{R}^N)\) such that \(\mathbf{w}_t\cdot\mathbf{n}=0\) a.e. in $\Gamma_C$.  Given \(\tau,t\in K\), the convex functional defined by:
	\[
	E_{t,\tau}(\mathbf{v}) := \frac{1}{2}\int_\Omega \boldsymbol\varepsilon(\mathbf{v}):\boldsymbol\Lambda\boldsymbol\varepsilon(\mathbf{v})\,{\rm d}x - \int_\Omega \mathbf{F}\cdot\mathbf{v} \,{\rm d}x - \int_{\Gamma_T}\mathbf{T}\cdot\mathbf{v}_{|\Gamma_T} - \int_{\Gamma_C} \tau\, v_n - \int_{\Gamma_C} ft\, \bigl| \mathbf{v}_t - \mathbf{w}_t \bigr|,
	\]
	where $v_n$ and $\mathbf{v}_t$ are defined by the relation $\mathbf{v}_{|\Gamma_C}=v_n\mathbf{n}+\mathbf{v}_t$,
	has a unique minimizer on the following closed subspace of \(H^1(\Omega;\mathbb{R}^N)\):
	\begin{equation}
		\label{eq:defV}
		V := \Bigl\{ \mathbf{v}\in H^1(\Omega;\mathbb{R}^N)\;\bigm| \mathbf{v}=\mathbf{0},\; \text{ on }\Gamma_U\Bigr\}.
	\end{equation}
	The unique minimizer \(\mathbf{u}\) of \(E_{t,\tau}\) on \(V\) is characterized by the variational inequality:
	\begin{multline}
		\label{eq:charact_ineq}
		\forall \mathbf{v}\in V,\qquad \int_\Omega \boldsymbol\varepsilon(\mathbf{u}):\boldsymbol\Lambda\boldsymbol\varepsilon(\mathbf{v}-\mathbf{u})\,{\rm d}x \geq \int_\Omega \mathbf{F}\cdot(\mathbf{v}-\mathbf{u}) \,{\rm d}x + \int_{\Gamma_T}\mathbf{T}\cdot(\mathbf{v}_{|\Gamma_T}-\mathbf{u}_{|\Gamma_T}) + \mbox{} \\
		\mbox{} + \int_{\Gamma_C} \tau\,( v_n-u_n) + \int_{\Gamma_C}ft\Bigl(\bigl| \mathbf{v}_t - \mathbf{w}_t \bigr|-\bigl| \mathbf{u}_t - \mathbf{w}_t \bigr|\Bigr).
	\end{multline}
	We define \(\mathscr{A}(t,\tau):=\mathbf{u}_{|\Gamma_C}\cdot\mathbf{n}\) the normal part of the trace of \(\mathbf{u}\) on \(\Gamma_C\) and \(At:=\mathscr{A}(t,t)\), so that:
	\[
	\mathscr{A}: K\times K \rightarrow H^{1/2}(\Gamma_C),\qquad A: K \rightarrow H^{1/2}(\Gamma_C).
	\]
\end{defi}

\begin{prop}
	\label{thm:boundedProp}
	Denoting by \(\mathbf{u}\) the minimizer of \(E_{\tau, t}\) on \(V\) introduced in Definition \ref{defu}, we have:
	\[
	\bigl\| \mathbf{u} \bigr\|_{H^1(\Omega,\mathbb{R}^N)} \leq C_1 + C_2\bigl(\|\tau\|_{H^{-1/2}}+\|t\|_{H^{-1/2}}\bigr),
	\]
	for two positive constants \(C_1,C_2>0\) that are independent of \(t,\tau\in K\).  In particular, the nonlinear operators \(\mathscr{A}\) and \(A\) are bounded (in the sense of Definition~\ref{thm:genDefi}).
\end{prop}

\noindent\textbf{Proof.} Taking \(\mathbf{v}=\mathbf{0}\) in \eqref{eq:charact_ineq}, we have:
\[
	\int_\Omega \boldsymbol\varepsilon(\mathbf{u}):\boldsymbol\Lambda\boldsymbol\varepsilon(\mathbf{u})\,{\rm d}x \leq \int_\Omega \mathbf{F}\cdot\mathbf{u} \,{\rm d}x + \int_{\Gamma_T}\mathbf{T}\cdot\mathbf{u}_{|\Gamma_T} + \int_{\Gamma_C} \tau\,u_n - \int_{\Gamma_C}ft\bigl| \mathbf{u}_t \bigr|.
\]
Hence, by estimates~\eqref{eq:multLischitz} and Proposition~\ref{thm:dualTrace},
\[
	\int_\Omega \boldsymbol\varepsilon(\mathbf{u}):\boldsymbol\Lambda\boldsymbol\varepsilon(\mathbf{u})\,{\rm d}x \leq \Bigl(\bigl\|\mathbf{F}\bigr\|_{L^2}+\bigl\|\mathbf{T}\bigr\|_{L^2}\Bigr)\bigl\|\mathbf{u}\bigr\|_{H^1} + C\Bigl(\bigl\|\tau\bigr\|_{H^{-1/2}}+\bigl\|t\bigr\|_{H^{-1/2}}\Bigr)\|\mathbf{u}_{|\Gamma_C}\|_{H^{1/2}}.
\]
where we have used that $\|\mathbf{u}_{|\Gamma_T}\|_{H^{1/2}}\leq C\|\mathbf u\|_{H^1}$ and that $f\in W^{1,\infty}(\Gamma_C)$.  Since $\|\mathbf{u}_{|\Gamma_C}\|_{H^{1/2}}\leq C\|\mathbf u\|_{H^1}$, the claim is a direct consequence of Korn's inequality. \qed

\bigskip
We now handle the task of proving that the operator \(A\) of Definition~\ref{thm:defAscr} is a Leray-Lions operator, that is, that the corresponding operator \(\mathscr{A}\) satisfies properties \textit{(i)}, \textit{(ii)}, \textit{(iii)} and \textit{(iv)} of Definition~\ref{thm:defLerayLions} in Appendix~A.  We will see that properties \textit{(i)}, \textit{(ii)} and \textit{(iii)} are easily checked and that all the difficulty concentrates on the proof of property~\textit{(iv)}.

\begin{prop}
	\label{thm:strongMonotProp}
	Let $\mathscr A$ be as in Definition \ref{defu} and let \(t\in K\) be arbitrary.  The mapping \(\tau\mapsto \mathscr{A}(t,\tau)\) is Lipschitz-continuous from \(K\subset H^{-1/2}\) to \(H^{1/2}\).  It also satisfies the strong monotonicity property:
	\[
		\forall t,\tau\in K, \qquad \bigl\langle\mathscr{A}(t,t)-\mathscr{A}(t,\tau)\,,\,t-\tau\bigr\rangle \geq C \bigl\| t-\tau \bigr\|_{H^{-1/2}}^2, 
	\]
	for some positive constant \(C>0\) independent of \(\,t\) and \(\tau\).  In particular, property~(i) of Definition~\ref{thm:defLerayLions} is fulfilled.
\end{prop}

\noindent\textbf{Proof.}
Let \(t,\tau^1,\tau^2\in K\).  We denote by \(\mathbf{u}^i\) the minimizer of \(E_{t,\tau^i}\) on \(V\) introduced in Definition \ref{defu}.  Taking \(\mathbf{u}^2\) as a test function in the variational inequality \eqref{eq:charact_ineq} characterizing \(\mathbf{u}^1\), \(\mathbf{u}^1\) as a test function in the variational inequality characterizing \(\mathbf{u}^2\) and taking the sum of the two corresponding inequalities, we get:
\begin{equation}
	\label{eq:monotIneq}
	\int_\Omega \boldsymbol\varepsilon(\mathbf{u}^1-\mathbf{u}^2):\boldsymbol\Lambda\boldsymbol\varepsilon(\mathbf{u}^1-\mathbf{u}^2)\,{\rm d}x \leq \int_{\Gamma_C} \bigl(\tau^1-\tau^2\bigr)\bigl(u_n^1-u_n^2\bigr).
\end{equation}
Since \(\mathbf{u}^1-\mathbf{u}^2=J_{{\Lambda}}(\mathbf{u}_{|\Gamma_C}^1-\mathbf{u}_{|\Gamma_C}^2)\), we get by Proposition \ref{thm:vectorTrace},
\[
\bigl\| u_n^1-u_n^2 \bigr\|_{H^{1/2}}^2 \leq C \int_\Omega \boldsymbol\varepsilon(\mathbf{u}^1-\mathbf{u}^2):\boldsymbol\Lambda\boldsymbol\varepsilon(\mathbf{u}^1-\mathbf{u}^2)\,{\rm d}x,
\]
for some positive constant \(C\) independent of \(u^1,u^2\).  It follows \(\|u_n^1-u_n^2\|_{H^{1/2}} \leq C \| \tau^1 - \tau^2\|_{H^{-1/2}}\), which is exactly the Lipschitz-continuity of the mapping \(\tau\mapsto \mathscr{A}(t,\tau)\).

Let \(\boldsymbol\tau\in H^{-1/2}(\Gamma_C;\mathbb{R}^N)\) be defined by:
\[
\forall \mathbf{v}\in V,\qquad \bigl\langle \boldsymbol\tau\,,\, \mathbf{v_{|\Gamma_C}} \bigr\rangle := \int_\Omega \boldsymbol\varepsilon(\mathbf{u}^1-\mathbf{u}^2):\boldsymbol\Lambda\boldsymbol\varepsilon(\mathbf{v})\,{\rm d}x.
\]
Integrating by parts, the right-hand side of the previous expression is seen to depend only on $\mathbf v_{|\Gamma_C}$, so that $\boldsymbol{\tau}$ is well-defined.
In addition, \(J_\Lambda'(\boldsymbol\tau)=\mathbf{u}^1-\mathbf{u}^2\) and \(\tau_n=\tau^1-\tau^2\), with $J'_\Lambda$ and $\tau_n$ introduced in Propositions \ref{thm:vectorTrace} and \ref{thm:dualTrace} respectively.  Combining Proposition~\ref{thm:dualTrace} with inequality~\eqref{eq:monotIneq}, we obtain:
\[
\bigl\| \tau^1-\tau^2 \bigr\|_{H^{-1/2}}^2\leq C	\Bigl\langle\mathscr{A}(t,\tau^1)-\mathscr{A}(t,\tau^2)\,,\,\tau^1-\tau^2\Bigr\rangle, 
\]
for some positive constant \(C>0\), independent of \(t,\tau^1,\tau^2\in K\).

Since the mapping \(\tau\mapsto \mathscr{A}(t,\tau)\) is bounded by Proposition~\ref{thm:boundedProp} and {since} Lipschitz-continuity obviously implies hemicontinuity, property~\textit{(i)} of Definition~\ref{thm:defLerayLions} is proved. \qed

\begin{prop}
	\label{thm:contProp}
	Let $\mathscr A$ be as in Definition \ref{defu} and let \(\tau\in K\) be arbitrary.  The mapping \(t\mapsto \mathscr{A}(t,\tau)\) is bounded and continuous from \(K\subset H^{-1/2}\) to \(H^{1/2}\).  In particular, property~(ii) of Definition~\ref{thm:defLerayLions} is fulfilled.
\end{prop}

\noindent\textbf{Proof.} 
The boundedness has already been proved in Proposition~\ref{thm:boundedProp}.  To prove the continuity, we first pick an arbitrary \(R>0\) and \(t^1,t^2\in K\), such that \(\|t^i\|_{H^{-1/2}}\leq R\), for \(i=1,2\).  Denoting by \(\mathbf{u}^i\) the minimizer of \(E_{t^i,\tau}\) on \(V\) introduced in Definition \ref{defu}, Proposition~\ref{thm:boundedProp} shows that \(\|\mathbf{u}^i\|_{H^1}\) (\(i=1,2\)) are bounded by a constant depending on \(\mathbf{F},\mathbf{T},\tau\) and \(R\), but independent of \(t^1\) and \(t^2\).  The variational inequalities characterizing \(\mathbf{u}^1\) and \(\mathbf{u}^2\) (see \eqref{eq:charact_ineq}) yield:
\begin{align*}
	\int_\Omega \boldsymbol\varepsilon(\mathbf{u}^1-\mathbf{u}^2):\boldsymbol\Lambda\boldsymbol\varepsilon(\mathbf{u}^1-\mathbf{u}^2)\,{\rm d}x 
	& \leq \int_{\Gamma_C} f\bigl(t^1-t^2\bigr)\Bigl(\bigl| \mathbf{u}_t^1 - \mathbf{w}_t \bigr| - \bigl| \mathbf{u}_t^2 - \mathbf{w}_t \bigr|\Bigr),\\
	& \leq \bigl\| f(t^1-t^2) \bigr\|_{H^{-1/2}} \Bigl(\bigl\| |\mathbf{u}_t^1 - \mathbf{w}_t| \bigr\|_{H^{1/2}} + \bigl\| |\mathbf{u}_t^2 - \mathbf{w}_t| \bigr\|_{H^{1/2}} \Bigr), \\
	& \leq C_1\bigl\| t^1-t^2 \bigr\|_{H^{-1/2}} \Bigl(\bigl\| \mathbf{u}_t^1 - \mathbf{w}_t \bigr\|_{H^{1/2}} + \bigl\| \mathbf{u}_t^2 - \mathbf{w}_t \bigr\|_{H^{1/2}} \Bigr), \\
	& \leq C_2\bigl\| t^1-t^2 \bigr\|_{H^{-1/2}} ,
\end{align*}
which shows that the mapping \(t\mapsto \mathscr{A}(t,\tau)\) is locally Hölder-continuous of exponent \(1/2\) from \(K\subset H^{-1/2}\) to \(H^{1/2}\). \qed

\begin{prop}
	\label{thm:propiii}
	Let $\mathscr A$ be as in Definition \ref{defu}.
	Let \(t^k\rightharpoonup t\) be a weakly converging sequence in \(K\) such that \({\lim_{k \rightarrow +\infty} \langle \mathscr{A}(t^k,t^k)-\mathscr{A}(t^k,t),t^k-t\rangle =0}\).  Then, for all \(\tau\in K\), the sequence \(\mathscr{A}(t^k,\tau)\) converges strongly in \(H^{1/2}\) towards \(\mathscr{A}(t,{\tau})\).  In particular, property~(iii) of Definition~\ref{thm:defLerayLions} is fulfilled.
\end{prop}

\noindent\textbf{Proof.}  The strong monotonicity property of Proposition~\ref{thm:strongMonotProp} entails that the sequence \(t^k\) converges strongly in {\(H^{-1/2}\) towards \(t\)}.  The continuity property of Proposition~\ref{thm:contProp} yields the strong convergence of \(\mathscr{A}(t^k,\tau)\) towards \(\mathscr{A}(t,\tau)\).\qed

\bigskip
At this point, there {remains} only to prove that \(\mathscr{A}\) also fulfills property~\textit{(iv)} of Definition~\ref{thm:defLerayLions} to reach the conclusion that the operator \(A:K \rightarrow H^{1/2}\) is Leray-Lions and therefore pseudomonotone.  If it is coercive in the sense of Definition~\ref{thm:defPseudoMonot}, in addition, {then} the solvability of the Signorini-Coulomb problem will follow from Theorem~\ref{thm:Brezis}.  As it turns out that the coercivity is easy to obtain and that property~\textit{(iv)} is very difficult, we first establish the coercivity result.

\begin{prop}
	\label{thm:propCoer}
	The operator \(A\) in Definition~\ref{defu} is coercive in the sense:
	\[
	\lim_{\substack{\|t\|_{H^{-1/2}}}\rightarrow +\infty,\\
	t\in K} \frac{\langle At,t\rangle}{\|t\|_{H^{-1/2}}} = +\infty.
	\]
\end{prop}

\noindent\textbf{Proof.}  Let \(t\in K\) and \(\mathbf{u}\) be the minimizer of \(E_{t,t}\) on \(V\), introduced in Definition~\ref{defu}.  Taking \(\mathbf{v}=\mathbf{0}\) as a test function in the variational inequality~\eqref{eq:charact_ineq} characterizing \(\mathbf{u}\), we obtain:
\begin{align*}
	\langle At,t\rangle & \geq  \int_\Omega \boldsymbol\varepsilon(\mathbf{u}):\boldsymbol\Lambda\boldsymbol\varepsilon(\mathbf{u})\,{\rm d}x - \int_\Omega \mathbf{F}\cdot\mathbf{u} \,{\rm d}x - \int_{\Gamma_T}\mathbf{T}\cdot\mathbf{u} \mbox{}+ \int_{\Gamma_C}ft\bigl| \mathbf{w}_t \bigr|,\\
	& \geq \int_\Omega \boldsymbol\varepsilon(\mathbf{u}):\boldsymbol\Lambda\boldsymbol\varepsilon(\mathbf{u})\,{\rm d}x - \int_\Omega \mathbf{F}\cdot\mathbf{u} \,{\rm d}x - \int_{\Gamma_T}\mathbf{T}\cdot\mathbf{u} \mbox{}- C\,\bigl\| t \bigr\|_{H^{-1/2}},
\end{align*}
for some positive constant \(C>0\) independent of \(t\) and \(\mathbf{u}\).  Let \(\mathbf{U}\) be the minimizer on \(V\) of:
\[
\mathbf{v}\mapsto \frac{1}{2}\int_\Omega \boldsymbol\varepsilon(\mathbf{v}):\boldsymbol\Lambda\boldsymbol\varepsilon(\mathbf{v})\,{\rm d}x - \int_\Omega \mathbf{F}\cdot\mathbf{v} \,{\rm d}x - \int_{\Gamma_T}\mathbf{T}\cdot\mathbf{v}_{|\Gamma_T} ,
\]
and \(\boldsymbol\tau\in H^{-1/2}(\Gamma_C;\mathbb{R}^N)\) be defined by:
\[
\forall \mathbf{v}\in V,\qquad \bigl\langle \boldsymbol\tau , \mathbf{v}_{|\Gamma_C}\bigr\rangle := \int_\Omega \boldsymbol\varepsilon(\mathbf{u}-\mathbf{U}):\boldsymbol\Lambda\boldsymbol\varepsilon(\mathbf{v})\,{\rm d}x,
\]
so that the definition is well-posed, \(J'_\Lambda(\boldsymbol\tau)=\mathbf{u}-\mathbf{U}\) and \(\tau_n=t\), with $J'_\Lambda$ and $\tau_n$ introduced in Propositions~\ref{thm:vectorTrace} and~\ref{thm:dualTrace} respectively.  Proposition~\ref{thm:dualTrace} yields:
\[
\bigl\| t \bigr\|_{H^{-1/2}}^2 \leq C \int_\Omega \boldsymbol\varepsilon(\mathbf{u}-\mathbf{U}):\boldsymbol\Lambda\boldsymbol\varepsilon(\mathbf{u}-\mathbf{U})\,{\rm d}x,
\]
for some positive constant \(C>0\) independent of \(t\) and \(\mathbf{u}\).  Bringing all together, we obtain:
\begin{align*}
	\langle At,t\rangle \mbox{}+ C\,\bigl\| t \bigr\|_{H^{-1/2}} & \geq  \int_\Omega \boldsymbol\varepsilon(\mathbf{u}):\boldsymbol\Lambda\boldsymbol\varepsilon(\mathbf{u})\,{\rm d}x - \int_\Omega \boldsymbol\varepsilon(\mathbf{U}):\boldsymbol\Lambda\boldsymbol\varepsilon(\mathbf{u})\,{\rm d}x ,\\
	& \geq \int_\Omega \boldsymbol\varepsilon(\mathbf{u}-\mathbf{U}):\boldsymbol\Lambda\boldsymbol\varepsilon(\mathbf{u}-\mathbf{U})\,{\rm d}x + \int_\Omega \boldsymbol\varepsilon(\mathbf{U}):\boldsymbol\Lambda\boldsymbol\varepsilon(\mathbf{u}-\mathbf{U})\,{\rm d}x,\\
	& \geq \frac{\|t\|_{H^{-1/2}}^2}{C} + \bigl\langle t,U_n + f |\mathbf{U}_t|\bigr\rangle,\\
	& \geq \frac{\|t\|_{H^{-1/2}}^2}{C} - \bigl\| t \bigr\|_{H^{-1/2}} \Bigl\|U_n + f |\mathbf{U}_t|\Bigr\|_{H^{1/2}},
\end{align*}
where in the third inequality we have used~\eqref{eq:charact_ineq} applied with the choice \(\mathbf{v}=\mathbf{u}+\mathbf{U}\).  The conclusion follows. \qed

\bigskip
Let now deal with property~\textit{(iv)} of Definition~\ref{thm:defLerayLions}.  It involves the handling of a product of two weakly converging sequences which turns out to raise huge difficulties.  The next theorem shows how to overcome such difficulties in the 2D case \(N=2\), under the assumptions that the elastic modulus tensor \(\boldsymbol\Lambda\) is that of isotropic elasticity at each point and that it is Lipschitz-continuous.  Its proof is postponed to Section~\ref{sec:fineProperty}.

Since $N=2$, the tangential direction is (pointwisely) unique, so that $\mathbf t_t$ and $\mathbf v_t$ are parallel vectors.  Hence, $\mathbf t_t$ and $\mathbf v_t$ will be replaced below by their scalar components $t_t$ and $v_t$.  We will denote by $\mathscr{M}(\overline \Gamma_C)$ the set of signed Radon measures with support in $\overline \Gamma_C$, and, for $\mu\in \mathscr{M}(\overline \Gamma_C)$, we will denote by $|\mu|$ the total variation of $\mu$.

\begin{theo}
	\label{thm:cornerstone2D}
	Let \(N=2\) and let \(\boldsymbol\Lambda\in W^{1,\infty}(\Omega)\) and isotropic.  For \(\tilde{t}\in H^{-1/2}(\Gamma_C)\cap \mathscr{M}(\overline\Gamma_C)\), we denote by \(\mathbf{u}\) the unique minimizer on \(V\) of:
	\[
	\mathbf{v} \mapsto \frac{1}{2}\int_\Omega \boldsymbol\varepsilon(\mathbf{v}):\boldsymbol\Lambda\boldsymbol\varepsilon(\mathbf{v})\,{\rm d}x - \int_{\Gamma_C} \tilde{t} \,v_t,
	\]
	and by \(u_n\in H^{1/2}(\Gamma_C)\) the normal part of the trace of \(\mathbf{u}\) on \(\Gamma_C\).  The linear mapping 
	\(\tilde{t}\mapsto u_n=L\tilde{t}\) is continuous for the strong topologies of \(H^{-1/2}(\Gamma_C)\) and \(H^{1/2}(\Gamma_C)\) (and also for the weak topologies, as it is linear). 
	
	Then, for any sequence \(t^k\in K\) converging weakly in \(H^{-1/2}(\Gamma_C)\) towards a limit \(t\in K\), and any sequence \(\tilde{t}^k\rightharpoonup \tilde{t}\) in \(H^{-1/2}(\Gamma_C)\cap \mathscr{M}(\overline\Gamma_C)\) such that \(|\tilde{t}^k| \leq -t^k\), we have:
	\[
	\lim_{k \rightarrow +\infty} \bigl\langle L\tilde{t}^k , t^k\bigr\rangle = \bigl\langle L\tilde{t} , t\bigr\rangle.
	\]
\end{theo}

\begin{rem}
	\label{thm:remCEMurat}
	To prove Theorem~\ref{thm:cornerstone2D}, it would be natural to try to prove that the injection of \(K\) into \(H^{-1/2}(\Gamma_C)\) is compact, in the spirit of the result of Murat~\cite{Murat}.  Actually, we can adapt the counterexample given at the second page of~\cite{Murat} to prove that the injection of \(K\) into \(H^{-1/2}(\Gamma_C)\) is \emph{not} compact.  In particular, the nonnegative function \(\epsilon|\log\epsilon|^{-1/2}/(x^2+\epsilon^2)\) converges weakly in \(H^{-1/2}(-1,1)\) towards \(0\), as \(\epsilon \rightarrow 0+\), but its \(H^{-1/2}(-1,1)\) norm converges towards a positive value (see Proposition~\ref{thm:BallardJarusek}).  Theorem~\ref{thm:cornerstone2D} is actually a fine property of the elasticity operator, as shown in Section~\ref{sec:fineProperty}.
\end{rem}

\begin{coro}
	\label{thm:coroPropiv}
	Let \(N=2\), \(\boldsymbol\Lambda\in W^{1,\infty}(\Omega)\) and isotropic, and let \(\mathscr{A}\) be as in Definition \ref{defu}.  Let \(\tau\in K\), \(t^k\rightharpoonup t\) be a weakly convergent sequence in \(K\) such that \(\mathscr{A}(t^k,\tau)\rightharpoonup F\) converges weakly in \(H^{1/2}(\Gamma_C)\) towards some limit \(F\).  Then, \(\lim_{k \rightarrow+\infty} \langle \mathscr{A}(t^k,\tau),t^k\rangle = \langle F,t\rangle\).  In other words, under the previous hypotheses, property (iv) of Definition~\ref{thm:defLerayLions} holds true and the operator \(A\) is a Leray-Lions operator.
\end{coro}

\noindent\textbf{Proof.} Let \(\mathbf{u}^k\in V\) be the minimizer on \(V\) of \(E_{t^k,\tau}\) (see Definition~\ref{defu}).  The idea is to split $\mathbf{u}^k$, and then $\mathscr{A}(t^k,\tau)={u}^k_n$, into two parts, the first not depending on $k$ and having $\tau$ as normal traction on $\Gamma_C$; the second having the same tangent traction  as \(\mathbf{u}^k\) on $\Gamma_C$.  The normal component of this second part will be written in terms of the linear operator $L$ introduced in Theorem \ref{thm:cornerstone2D}.

Let \(\mathbf{U}\in V\) be the minimizer on \(V\) of:
\[
\mathbf{v} \rightarrow \frac{1}{2} \int_\Omega \boldsymbol\varepsilon(\mathbf{v}):\boldsymbol\Lambda\boldsymbol\varepsilon(\mathbf{v})\,{\rm d}x - \int_\Omega \mathbf{F}\cdot\mathbf{v}\,{\rm d}x - \int_{\Gamma_T}\mathbf{T}\cdot\mathbf{v} - \int_{\Gamma_C}\tau v_n.
\]
Let \(\mathbf{\tilde{t}}^k\) be the element of \(H^{-1/2}(\Gamma_C;\mathbb{R}^N)\) defined by:
\begin{equation}
	\label{eq:deftildeT}
	\forall \mathbf{v}\in V,\qquad \bigl\langle \mathbf{\tilde{t}}^k,\mathbf{v}_{|\Gamma_C}\bigr\rangle := \int_\Omega \boldsymbol\varepsilon(\mathbf{u}^k-\mathbf{U}):\boldsymbol\Lambda\boldsymbol\varepsilon(\mathbf{v})\,{\rm d}x.
\end{equation}
Then, as both \(\mathbf{u}^k\) and \(\mathbf{U}\) have $\tau$ as normal traction on $\Gamma_C$, the normal part vanishes: \(\tilde{t}_n^k=\tau-\tau=0\).  We denote by \(\tilde{t}^k\) the (scalar-valued) component of \(\mathbf{\tilde{t}}^k\) along the tangent vector (as \(N=2\)).  Applying inequality~\eqref{eq:charact_ineq} to \(\mathbf{u^k}\), with \(\mathbf{u^k}+\mathbf{v}\) as test function, and taking~\eqref{eq:deftildeT} into account, we get:
\[
\forall v\in H^{1/2}(\Gamma_C),\qquad  \bigl\langle \tilde{t}^k,v\bigr\rangle - \bigl\langle ft^k,|v|\bigr\rangle \geq 0,
\]
which, as \(ft^k\) is a nonpositive measure, entails that \(|\langle\tilde{t}^k,v\rangle| \leq -\langle ft^k,|v|\rangle\) and therefore that \(\tilde{t}^k\) is a measure with total variation satisfying \(|\tilde{t}^k| \leq -ft^k\).  As the sequence \(t^k\) is assumed to weakly converge in \(H^{-1/2}(\Gamma_C)\), it is bounded in the same space.  By the inequality \(|\langle\tilde{t}^k,v\rangle| \leq -\langle ft^k,|v|\rangle \leq  \|f\|_{W^{1,\infty}}\|t^k\|_{H^{-1/2}}\|v\|_{H^{1/2}}\), the sequence \(\tilde{t}^k\) is also bounded in \(H^{-1/2}(\Gamma_C)\).  Also, by linearity, we have:
\[
\mathscr{A}(t^k,\tau) = U_n+(u^k_n-U_n)=\mbox{}U_n + L\tilde{t}^k,
\]
where \(L\) denotes the linear operator defined in the statement of Theorem~\ref{thm:cornerstone2D}.  The sequence 
\(\langle \mathscr{A}(t^k,\tau),t^k\rangle = \langle  U_n + L\tilde{t}^k,t^k \rangle\) is bounded.  Extracting a subsequence if necessary, it converges.  Possibly extracting another subsequence, we may assume that \(\tilde{t}^k\) converges weakly in \(H^{-1/2}(\Gamma_C)\) towards some limit \(\tilde{t}\).  As \(\mathscr{A}(t^k,\tau)\rightharpoonup F\) and \(L\) is continuous for the weak topologies, we have the identity \(F=U_n+L\tilde{t}\).  Now, Theorem~\ref{thm:cornerstone2D} yields:
\[
\lim_{k \rightarrow +\infty} \bigl\langle \mathscr{A}(t^k,\tau),t^k\bigr\rangle = \langle F,t \rangle.
\]
As the limit is the same for all extracted subsequences, the convergence of the whole limit holds true.\qed

\bigskip
Bringing Propositions \ref{thm:strongMonotProp}, \ref{thm:contProp}, \ref{thm:propiii}, \ref{thm:propCoer} and Corollary \ref{thm:coroPropiv} together, and invoking Proposition~\ref{thm:LerayLionsImpliesPseudomonotone} in Appendix~A, we have proved:

\begin{theo}
	Let $N=2$ and let \(\boldsymbol\Lambda\in W^{1,\infty}(\Omega)\) and isotropic.  Then, the operator \(A\) in Definition~\ref{thm:defAscr} is a Leray-Lions operator and is therefore pseudomonotone.  In addition, it is coercive.
\end{theo}

We can now rely on Brézis's theorem (Theorem~\ref{thm:Brezis} in Appendix~A) to solve variational inequalities based on the operator \(A\).

\begin{theo}
	\label{theoexist}
	Let \(N=2\) and let \(\boldsymbol\Lambda\in W^{1,\infty}(\Omega)\) and isotropic, \(f\in W^{1,\infty}(\Gamma_C)\), \(\mathbf{F}\in L^2(\Omega,\mathbb{R}^N)\), \(\mathbf{T}\in L^2(\Gamma_T,\mathbb{R}^N)\), \(g\in H^{1/2}(\Gamma_C)\) and \(\mathbf{w}_t\in H^{1/2}(\Gamma_C;\mathbb{R}^N)\).  Let also $A$ be as in Definition~\ref{defu}.  Then, there exists \(t\in K\) satisfying the variational inequality:
	\[
	\forall \hat{t}\in K, \qquad \bigl\langle At-g,\hat{t}-t\bigr\rangle \geq 0.
	\]
\end{theo}

It is readily checked that the unique minimizer \(\mathbf{u}\in V\) of \(E_{t,t}\) on \(V\) (see Definition~\ref{defu}) with $t$ given by Theorem \ref{theoexist}, solves the formal time-incremental Signorini-Coulomb problem~\eqref{eq:SignoriniCoulombDiscret}.  This means that, whenever the minimizer is smooth enough, the conditions are fulfilled pointwisely.  In particular, condition~\eqref{eq:Frict} of problem~\eqref{eq:SignoriniCoulombDiscret} gives on \(\Gamma_C\):
\begin{eqnarray*}
	u_t-w_t=0 &\Rightarrow& |t_t|\leq -ft_n,\\
	u_t-w_t\neq 0 &\Rightarrow& t_t= -ft_n\frac{u_t-w_t}{|u_t-w_t|}.
\end{eqnarray*}

\subsection{Proof of Theorem \ref{thm:cornerstone2D}}

\label{sec:fineProperty}

In this section, we give a detailed proof of Theorem~\ref{thm:cornerstone2D} which is the cornerstone of our proof of the existence of solutions for the Signorini-Coulomb problem.

In order to both assess the difficulty of the matter involved and to display the core idea of our strategy of proof, it is useful to examine first the particular case of the isotropic homogeneous 2d half-space.  Indeed, in this case, the explicit knowledge of the fundamental solution makes it possible to express the solution of the standard elastic Neumann problem for the half-space in terms of a convolution product.  In other words, the operator \(L\) of Theorem~\ref{thm:cornerstone2D} reduces in this case to a convolution operator with a function which is known explicitly.

\begin{lemm}
	[Fundamental solution of the 2d isotropic half-space]
	\label{thm:fundSol2d}
	Consider the half-space \(\{(x,y)\in \mathbb{R}^2\;|\;y<0\}\) filled with an isotropic homogeneous linearly elastic material of Young modulus \(E>0\) and Poisson ratio \(\nu\in\left]-1,1/2\right[\).  Then, all the tempered distributional displacements \(\mathbf{u}^0\) satisfying the elastic equilibrium equations with zero body forces and surface traction equal to \(\mathbf{T}\delta\), where \(\delta\) is the Dirac measure at \((0,0)\) and \(\mathbf{T}=(T_x,T_y)\) is a given element of \(\mathbb{R}^2\), are of the form:
	\begin{multline*}
		\mathbf{u}^0(x,y) = T_x \mathbf{u}^{0,x}(x,y) + T_y \mathbf{u}^{0,y}(x,y) +\mbox{}\\
		\mbox{}+ \bigl(U_x-\Theta y +(1-\nu^2)\Sigma x\bigr)\mathbf{e}_x + \bigr(U_y +\Theta x -\nu (1+\nu)\Sigma y\bigr)\mathbf{e}_y ,
	\end{multline*}
	where:
	\begin{align*}
		u_x^{0,x}(x,y) & := \frac{1}{\pi E} \biggl[-(1-\nu^2)\log(x^2+y^2) - (1+\nu)\frac{y^2}{x^2+y^2}\biggr], \\
		u_x^{0,y}(x,y) & := \frac{1}{\pi E} \biggl[-(1-2\nu)(1+\nu)\arctan\frac{x}{y} + (1+\nu)\frac{xy}{x^2+y^2}\biggr], \\
		u_y^{0,x}(x,y) & := \frac{1}{\pi E} \biggl[(1-2\nu)(1+\nu)\arctan\frac{x}{y} + (1+\nu)\frac{xy}{x^2+y^2}\biggr],  \\
		u_y^{0,y}(x,y) & := \frac{1}{\pi E} \biggl[-(1-\nu^2)\log(x^2+y^2) + (1+\nu)\frac{y^2}{x^2+y^2}\biggr],
	\end{align*}
	and \(U_x,U_y,\Theta,\Sigma\) denote four arbitrary real constants.  
\end{lemm}

\noindent\textbf{Proof.} This is a classical result, obtained by use of the Fourier transform, and the verification is left to the reader (this is in fact a particular case of the proof of Theorem \ref{thm:anisotropicD2N} in Section \ref{sec:N2Danisotrop}). \qed

\begin{rem}
	The constant \(\Sigma\) is readily seen to be a component of the stress tensor at infinity.  The three remaining constants represent an arbitrary overall rigid displacement.   This arbitrary affine displacement plays no role in the sequel where \(U_x=U_y=\Theta=\Sigma=0\) will be systematically chosen.  
\end{rem}

Given a compactly supported surface traction distribution \(\mathbf{\tilde{t}}\in H^{-1/2}(-1,1;\mathbb{R}^2)\) with normal and tangential parts \(\tilde{t}_n\) and \(\tilde{t}_t\), the normal part of the surface displacement reads as:
\[
u_n = u_y = - \frac{2(1-\nu^2)}{\pi E} \,\log|\cdot| * \tilde{t}_n - \frac{(1-2\nu)(1+\nu)}{2 E} \,\text{sgn}(\cdot) * \tilde{t}_t,
\]
up to an arbitrary overall rigid displacement component which plays no role in the sequel and is therefore omitted.  Above, \(\text{sgn}(\cdot)\) stands for the sign function and \(*\) for the convolution product.  Hence, in the case of the half-space, the normal component of the surface displacement induced by a tangential traction distribution is obtained as a convolution with the `sign' function.  Therefore, the operator \(L:H^{-1/2}(-1,1)\rightarrow H^{1/2}(-1,1)\) in the statement of Theorem~\ref{thm:cornerstone2D}, is given by:
\[
L\tilde{t} = \text{sgn}(\cdot) * \tilde{t},
\]
where the extension by zero of \(\tilde{t}\) defines a distribution in \(H^{-1/2}(\mathbb{R})\) and   \(\text{sgn}(\cdot) * \tilde{t}\) denotes actually the restriction to \(\left]-1,1\right[\) of \(\text{sgn}(\cdot) * \tilde{t}\in H_{\rm loc}^{1/2}(\mathbb{R})\).  The (nonzero) multiplicative constant has been dropped.  This operator is  linear continuous from \(H^{-1/2}(-1,1)\) into \(H^{1/2}(-1,1)\).  The range of this operator contains all the closed subspaces of elements in \(H^{1/2}(-1,1)\) having support contained in a fixed compact of \(\left]-1,1\right[\). Hence, the operator \(L:H^{-1/2}(-1,1)\rightarrow H^{1/2}(-1,1)\) is linear continuous, but not compact.  In the particular case of the isotropic homogeneous 2d half-space, Theorem~\ref{thm:cornerstone2D} reduces to the following result for which we are going to provide a direct proof.

\begin{theo}
	\label{thm:cornerstoneHalfSpace}
	Let \(t^k\geq 0\) be a sequence in the nonnegative cone of \(H^{-1/2}(-1,1)\) that converges weakly in \(H^{-1/2}(-1,1)\) towards a limit \(t\).  Let also \(\tilde{t}^k\) be a sequence in \(H^{-1/2}(-1,1)\cap \mathscr{M}([-1,1])\) that converges weakly in \(H^{-1/2}(-1,1)\) towards a limit \(\tilde{t}\) and such that \(|\tilde{t}^k|\leq t^k\), for all \(k\).  Then:
	\[
	\lim_{k\rightarrow+\infty} \bigl\langle \text{\rm sgn}(\cdot)*\tilde{t}^k , t^k\bigr\rangle = \bigl\langle \text{\rm sgn}(\cdot)*\tilde{t} , t\bigr\rangle.
	\]
\end{theo}

\noindent\textbf{Proof.}  Since \(\|t^k\|_{\mathscr{M}([-1,1])} = \langle t^k,1 \rangle\) and the weakly converging sequence \(t^k\) is bounded in \(H^{-1/2}(-1,1)\), it is also bounded in \(\mathscr{M}([-1,1])\).  It therefore converges weakly-* in \(\mathscr{M}([-1,1])\) towards \(t\). 
As \(|\tilde{t}^k|\leq t^k\), the sequence \(\tilde{t}^k\) is also bounded in \(\mathscr{M}([-1,1])\) and therefore converges weakly-* in \(\mathscr{M}([-1,1])\) towards \(\tilde{t}\).  

The convolution with the sign function is an antisymmetric operation, that is 
\[
\bigl\langle \text{\rm sgn}(\cdot)*\tilde{t}^k , t^k\bigr\rangle=
-\bigl\langle \tilde{t}^k , \text{\rm sgn}(\cdot)*t^k\bigr\rangle.
\]
For \(t^k\) in \(H^{-1/2}(-1,1)\cap \mathscr{M}([-1,1])\), the primitives \(u^k:=\text{sgn}(\cdot)*t^k\) are functions with bounded variation that are also in \(H_{\rm loc}^{1/2}(\mathbb{R})\). Such functions have a left and right limit at each point, as they have bounded variation.  But the definition~\eqref{eq:defH12} of the \(H^{1/2}(-1,1)\)-norm implies that a function with a jump is not in \(H^{1/2}(-1,1)\).  Hence, any function in \(H_{\rm loc}^{1/2}(\mathbb{R})\) with bounded variation has equal left and right limits at each point and is therefore continuous:
\[
(H^{1/2}\cap BV)(-1,1)\subset C^0([-1,1])
\]
and all the \(u^k:=\text{sgn}(\cdot)*t^k\) are continuous functions.  They are also nondecreasing, as \(t^k\geq 0\).  We can therefore apply the two following arguments which yield the strong convergence in \(C^0([-1,1])\) of the sequence \(u^k\) towards \(u\).
\begin{enumerate}
	\item The weak-* convergence in \(\mathscr{M}([-1,1])\) of \(t^k\) entails the pointwise convergence of the sequence \(\int_{-1}^xt^k\) and therefore of the sequence \(u^k:=\text{sgn}(\cdot)*t^k\) towards \(u\), for all \(x\in [-1,1]\).  To see it, consider \(x\in\left]-1,1\right[\) and \(\epsilon>0\).  As the measure \(t\) has no atoms, we can build a continuous function \(\varphi_\epsilon:[-1,1]\rightarrow[0,1]\), supported in \([-1,x]\) and a continuous function \(\psi_\epsilon:[-1,1]\rightarrow[0,1]\) that takes the value \(1\) all over \([-1,x]\), such that:
	\[
	(1-\epsilon)\int_{-1}^{x}t \leq \int_{-1}^{1}\varphi_\epsilon t \leq \int_{-1}^{x}t \leq \int_{-1}^{1}\psi_\epsilon t \leq (1+\epsilon) \int_{-1}^{x}t.
	\]
	As:
	\[
	\forall k\in \mathbb{N},\qquad \int_{-1}^{1}\varphi_\epsilon t^k \leq \int_{-1}^x t^k \leq \int_{-1}^{1}\psi_\epsilon t^k,
	\]
	we obtain that the sequence of functions \(x\mapsto \int_{-1}^x t^k\) converges pointwisely towards \(x\mapsto \int_{-1}^x t\). 
	\item Each function $u^k$ is monotone and pointwisely converging to a continuous function in the compact set $[-1,1]$.  This in fact forces $u^k$ to converge uniformly.  Indeed, the pointwise limit \(u\), being continuous on the compact \([-1,1]\), is uniformly continuous.  Given \(\epsilon>0\), there exist finitely many \(-1=x_0<x_1\cdots<x_n=1\) such that, for all \(i\) and all \(x\in [x_i,x_{i+1}]\), one has \(|u(x)-u(x_i)|<\epsilon\).  Also, for sufficiently large \(k\), and all \(i\), one has \(|u(x_i)-u^k(x_i)|<\epsilon\).  Using the fact that all the \(u^k\) are nondecreasing, we get \(|u^k(x)-u(x)|< 5\epsilon\), for sufficiently large \(k\) and all \(x\in [-1,1]\).  That is, the sequence \(u^k\) converges strongly in \(C^0([-1,1])\) towards \(u\).
\end{enumerate}
This entails \(\lim_{k\to+\infty} \langle \tilde{t}^k,u^k\rangle =\langle \tilde{t},u\rangle\), and therefore the claim.\qed

\begin{rem}
	In Theorem~\ref{thm:cornerstoneHalfSpace}, the function \(\text{sgn}(\cdot)\) could be replaced by any monotone function.  As any function with bounded variation is the difference of two monotone functions, it is readily checked that Theorem~\ref{thm:cornerstoneHalfSpace} holds true, with the function \(\text{sgn}(\cdot)\) replaced by an arbitrary function \(h\in \textsl{BV}([-r,r])\), with $r>2$. 
\end{rem}

\begin{rem}
	The similar convolution operator \(\tilde{t} \mapsto \log|\cdot|\,*\,\tilde{t}\) is an isomorphism from \(H^{-1/2}(-1,1)\) onto \(H^{1/2}(-1,1)\).  But we cannot replace the function \(\text{sgn}(\cdot)\) in the above theorem by the function \(\log|\cdot|\), otherwise the conclusion would break.  To see it, recall that the mapping \(t\mapsto \sqrt{-\langle \log|\cdot|*t,t\rangle}\) is a norm on \(H^{-1/2}(-1,1)\) which is equivalent to the norm of \(H^{-1/2}\) (see, for example, \cite[Theorem 3]{bibiJiri} or Theorem~\ref{thm:BallardJarusek} in this paper).  In the particular case where \(\tilde{t}^k=t^k\), the conclusion of the theorem, with \(\text{sgn}(\cdot)\) replaced by \(\log|\cdot|\), would actually be the convergence of the norms, which together with weak convergence, classically yields strong convergence.  In other terms, the above theorem, with \(\text{sgn}(\cdot)\) replaced by \(\log|\cdot|\), would imply that the injection of the nonnegative cone of \(H^{-1/2}(-1,1)\) into \(H^{-1/2}(-1,1)\) is compact.  But the injection of the nonnegative cone of \(H^{-1/2}(-1,1)\) into \(H^{-1/2}(-1,1)\) is not compact (a counterexample is given in Remark~\ref{thm:remCEMurat}).  Therefore, one cannot replace the function \(\text{sgn}(\cdot)\) in the above theorem by the function \(\log|\cdot|\).  In other words,  the fact that the operator \(L\) maps the \emph{tangential} traction to the \emph{normal} displacement plays an essential role.  In that sense, the above theorem is a fine property of the elasticity operator rather than a general topological property.
\end{rem}

\begin{rem}
	\label{r:3d}
	The proof of Theorem \ref{thm:cornerstoneHalfSpace} can be rephrased by saying that the inclusion
	\[
	(H^{1/2}\cap BV)(-1,1)\subset C^0([-1,1])
	\]
	is compact in the following sense: for $u^k,u\in(H^{1/2}\cap BV)(-1,1)$, we have
	\[
	u^k\rightharpoonup u \quad \text{ weakly-}H^{1/2} \text{ and strictly-}BV \quad \implies \quad u^k\to u \quad \text{strongly-}L^\infty,\]
	where 'strictly-$BV$' means that $u^k\to u$ strongly-$L^1$ and $|u^k{}'|([-1,1])\to|u'|([-1,1])$.  This compactness property fails in higher dimension.
	
	The only point of the article where we need $N=2$ is precisely Theorem \ref{thm:cornerstoneHalfSpace}.  In dimension $N=3$, in the simplest case, one is led to pass to the limit in the following expression
	\[
	\lim_{k\rightarrow+\infty} \bigl\langle \frac{x}{x^2+y^2}*\tilde{t}^k , t^k\bigr\rangle
	\]
	where $\tilde t^k$ and $t^k$ are supported in $\left]-1,1\right[^2$ and satisfy hypotheses analogous to those of Theorem \ref{thm:cornerstoneHalfSpace}.  Due to the previous discussion, the 2D argument does not generalize directly to this case and further investigations are needed.			
\end{rem}

\begin{rem}
	The above argument is reminiscent to that employed in~\cite{Michelle}, where the frictionless dynamic Signorini problem was formulated as a variational inequality.  Then, in the case of the half-space, the authors employed partial Fourier techniques to prove some monotonicity properties of the Steklov-Poincar\'e operator (the analogue of the Neumann-to-Dirichlet operator in the case of the wave equation).  However, the results that they obtained for the half-space did not generalize to the case of bounded bodies.  We are now going to prove that our analysis does not meet this difficulty.
\end{rem}

To extend Theorem~\ref{thm:cornerstoneHalfSpace} to more general geometries, that is, to the case of bounded bodies, we will rely on the classical idea that the fundamental solution of the bounded body is locally that of the half-space with the addition of a smooth correction.  This general idea amounts to describing the Neumann-to-Dirichlet operator as a pseudo-differential operator. In that context, the computation of the fundamental solution of the half-space really means to calculate the principal symbol of the operator near a point of the boundary while the remaining part is compact. The following lemma precisely implements this idea, while avoiding to rely on the full machinery of pseudo-differential operators.

\begin{lemm}
	\label{thm:represGreen}
	Let \(N=2\) and \(\boldsymbol\Lambda\in W^{1,\infty}\) be isotropic at each \(x\in\Omega\).  Let also \(t\in H^{-1/2}(\Gamma_C)\cap \mathscr{M}(\overline{\Gamma}_C)\).   We denote by \(\mathbf{u}\) the unique minimizer on \(V\) of:
	\[
	\mathbf{v} \mapsto \frac{1}{2}\int_\Omega \boldsymbol\varepsilon(\mathbf{v}):\boldsymbol\Lambda\boldsymbol\varepsilon(\mathbf{v})\,{\rm d}x - \int_{\Gamma_C} t \,v_t.
	\]
	Then, the normal part \(u_n\) of the trace of \(\mathbf{u}\) on \(\Gamma_C\) can be put under the form:
	\[
	\text{for a.a. }x\in\Gamma_C,\qquad u_n(x) = \int_{x'\in\Gamma_C} t\,g(x,x') 
	\]
	for some function \(g:\Gamma_C\times\Gamma_C \rightarrow\mathbb{R}\), independent of \(t\), of the form:
	\[
	g(x,x') = \frac{(1-2\nu(x'))(1+\nu(x'))}{2 E(x')}\,h(x,x') + \tilde{g}(x,x'),
	\]
	where the function \(h\) takes only the two values \(-1,+1\) and is continuous on \(\Gamma_C\times\Gamma_C\setminus \Delta\), with \(\Delta:=\{(x,x)|\:x\in\Gamma_C\}\) being the diagonal, and where \(\tilde{g}: x' \rightarrow g(\cdot,x')\in C^0(\overline{\Gamma}_C;H^{1/2}(\Gamma_C))\).
\end{lemm}

\noindent\textbf{Proof.}

\noindent\textbf{Step 1.} \textit{Construction of a fundamental solution for the bounded body.  Fix \(x'\in \Gamma\), with \(\Gamma := \partial\Omega\setminus\overline{\Gamma}_U\), and \(\mathbf{\tilde{T}}=(\tilde{T}_x,\tilde{T}_y)\in \mathbb{R}^2\).  Then, there exist unique \((\tilde{\mathbf{u}},\tilde{\mathbf{u}}_{|\Gamma})\in L^1(\Omega)\times L^1(\Gamma)\)  such that:
\begin{equation}
	\label{eq:veryWeak}
	\forall\mathbf{v}\in W^{2,\infty}(\Omega;\mathbb{R}^2)\; \text{with }\mathbf{v}_{|\Gamma_U}=\mathbf{0},\quad \int_\Omega \tilde{\mathbf{u}}\cdot \text{\rm div}\,\boldsymbol\sigma(\mathbf{v})\,{\rm d}x = -\mathbf{\tilde{T}}\cdot \mathbf{v}(x')+\int_\Gamma \tilde{\mathbf{u}}_{|\Gamma}\cdot\boldsymbol{\sigma}(\mathbf{v})\,\mathbf{n}.
\end{equation}
The function \(\tilde{\mathbf{u}}\) has the additional regularity \(\tilde{\mathbf{u}}\in H^1(\Omega\setminus\mathscr{V}_{x'})\) where \(\mathscr{V}_{x'}\) is an arbitrary neighborhood of \(x'\), and \(\tilde{\mathbf{u}}_{|\Gamma}\) coincides with the trace of \(\tilde{\mathbf{u}}\) on \(\Gamma\).  In the sequel, the notation \(\mathbf{g}(x,x'):=\tilde{\mathbf{u}}(x)\) will be used.
}

\medskip
This is a classical matter and the proof is sketched only for the sake of completeness.  Consider the mapping:
\[
\mathscr{L}
\left\{
\begin{array}{rcl}
	L^2(\Omega;\mathbb{R}^2)\times  L^2(\Gamma;\mathbb{R}^2) & \rightarrow & L^2(\Omega;\mathbb{R}^2)\times L^2(\Gamma;\mathbb{R}^2) \\
	(\mathbf{F},\mathbf{T}) & \mapsto & (\mathbf{u},\mathbf{u}_{|\Gamma})
\end{array}
\right.
\]
where \(\mathbf{u}\) denotes the unique minimizer on \(V=\{\mathbf{v}\in H^1(\Omega;\mathbb{R}^2)|\mathbf{v}_{|\Gamma_U}=\mathbf{0}\}\) of 
\[
\mathbf{v} \mapsto \frac{1}{2}\int_{\Omega}\boldsymbol\varepsilon(\mathbf{v}):\boldsymbol\Lambda\boldsymbol\varepsilon(\mathbf{v})\,{\rm d}x - \int_\Omega \mathbf{F}\cdot\mathbf{v}\,{\rm d}x - \int_\Gamma \mathbf{T}\cdot\mathbf{v}.
\]
The linear mapping \(\mathscr{L}=\mathscr{L}^*\) is self-adjoint.  By standard regularity results for elliptic boundary-value problems, if \((\mathbf{F},\mathbf{T})\in L^p(\Omega)\times L^p(\Gamma)\) for some \(p>2\), then \(\mathbf{u}\in W^{1,p}(\Omega)\subset C^0(\overline{\Omega})\) (as the dimension of space is \(2\)).  Hence, \(\mathscr{L}:L^p(\Omega)\times L^p(\Gamma) \rightarrow C^0(\overline{\Omega})\times C^0(\overline{\Gamma})\).  By duality, \(\mathscr{L}^*:\mathscr{M}(\overline{\Omega})\times\mathscr{M}(\overline{\Gamma})\rightarrow L^{p'}(\Omega)\times L^{p'}(\Gamma) \) can be seen as an extension of \(\mathscr{L}\), as \(\mathscr{L}\) is self-adjoint.  Then, it is readily checked that \((\tilde{\mathbf{u}},\tilde{\mathbf{u}}_{|\Gamma}):=\mathscr{L}^*(\mathbf{0},\mathbf{\tilde{T}}\delta_{x'})\) (where \(\delta_{x'}\) is the Dirac measure at \(x'\)) satisfies~\eqref{eq:veryWeak}.  Uniqueness and the fact that \(\tilde{\mathbf{u}}\in H^1(\Omega\setminus\mathscr{V}_{x'})\) is a consequence of the standard local regularity results for elliptic problems which apply to very weak solutions such as the one satisfying~\eqref{eq:veryWeak}. 

\medskip
\noindent\textbf{Step 2.} \textit{Regularity with respect to \(x\) of the trace of \(\mathbf{g}(x,x')\) on \(\Gamma:=\partial\Omega\setminus\overline{\Gamma}_U\).  Given \(x'\in \Gamma\) and \(\mathbf{\tilde{T}=(\tilde{T}_x,\tilde{T}_y)}\in \mathbb{R}^2\setminus \{\mathbf{0}\}\) a unit tangent vector at \(x'\) (that is, \(\mathbf{\tilde{T}}\cdot\mathbf{n}(x')=0\) and \(\tilde{T}_x^2+\tilde{T}_y^2=1\)), there exists \(x\mapsto \tilde{g}(x,x')\in H^{1/2}(\Gamma)\), such that:
\[
\text{for a.a. }x\in \Gamma,\qquad g_n(x,x') = \frac{(1-2\nu(x'))(1+\nu(x'))}{2 E(x')}\,h(x,x') + \tilde{g}(x,x'),
\]
where the function \(h\) takes only the two values \(-1,+1\) and is continuous on \(\Gamma\times\Gamma\setminus \Delta\), with \(\Delta:=\{(x,x)|\:x\in\Gamma\}\) being the diagonal.}

Pick \(x'\in \Gamma\).  There exist an open neighborhood \(\mathscr{V}_{x'}\) of \(x'\) in \(\mathbb{R}^2\) with \(\mathscr{V}_{x'}\cap\Gamma_U=\varnothing\), an open ball \(B\) centered at \(\mathbf 0\) and a \(W^{2,\infty}\)-diffeomorphism \(\phi : \mathscr{V}_{x'} \rightarrow B\), such that $\phi(x')=\mathbf{0}$, \(\mathscr{R}:=\nabla\phi(x')\) is a rotation and such that \(\phi\) induces a diffeomorphism of \(\mathscr{V}_{x'}\cap\Omega\) onto \(D:=B\cap\{y<0\}\).  We take \(\mathbf{T}=\mathscr{R}\mathbf{\tilde{T}}\) (so that \(T_y=0\)), \(E=E(x')\), \(\nu=\nu(x')\) and \(U_x=U_y=\Omega=\Sigma=0\) in the definition of \(\mathbf{u}^0\) in Lemma~\ref{thm:fundSol2d}.  Careful examination of the expression of \(\mathbf{u}^0\) shows that \(\mathbf{u}^0\) is in \(W^{1,p}(D;\mathbb{R}^2)\), for all \(p\in\left[1,2\right[\), and is \(C^\infty\) in \(\overline{D}\setminus\{0\}\).  Therefore the displacement \(\mathbf{u}^0\circ\phi\) is in \(W^{1,p}(\mathscr{V}_{x'}\cap\Omega)\).  Considering an appropriate cut-off function \(\varphi\) in \(C_c^\infty(\mathbb{R}^2)\) which is identically equal to \(1\) on a neighborhood of \(x'\), the displacement \(\varphi\mathscr{R}^{-1}(\mathbf{u}^0\circ\phi)\) has compact support in \(\mathscr{V}_{x'}\).  We denote by \(\mathbf{\tilde{u}}^0\) the extension to \(\overline{\Omega}\) by zero.  In the particular case where \(\mathbf{\tilde{T}}=(\tilde{T}_x,\tilde{T}_y)\in \mathbb{R}^2\setminus \{\mathbf{0}\}\) is a unit tangent vector at \(x'\), the following property of \(\mathbf{\tilde{u}}^0\) is readily checked from the explicit expression of \(\mathbf{u}^0\) in Lemma~\ref{thm:fundSol2d}.
\begin{quote}
	Since $\mathbf{\tilde{T}}$ is tangent, the normal trace \(\tilde{u}_n^0\) on \(\partial\Omega\) is in \(W^{1,\infty}(\partial\Omega\setminus\{x'\})\) (with support in \(\mathscr{V}_{x'}\)).  It therefore has limits as \(x \rightarrow x'-\) and \(x \rightarrow x'+\), but these two limits are different (as \(u_2^0(0-,0-)\neq u_2^0(0+,0-)\)).  Actually, we have:
	\[
	\Bigl| \tilde{u}_n^0(x'-) - \tilde{u}_n^0(x'+) \Bigr| = \frac{(1-2\nu(x'))(1+\nu(x'))}{E(x')}.
	\]
\end{quote}
To get the claim of Step 2, it is now sufficient to prove that \(\mathbf{g}(\cdot,x')-\mathbf{\tilde{u}}^0(\cdot)\in H^1(\Omega)\).  By the definition of \(\mathbf{g}(\cdot,x')\), this is going to be a consequence of the following statement.
\begin{quote}
	There exist \(\boldsymbol\tau(\cdot,x')\in L^p(\Gamma;\mathbb{R}^2)\), for all \(p\in \left[1,2\right[\), and \(\boldsymbol\Phi(\cdot,x')\in L^p(\Omega;\mathbb{R}^2)\), for all \(p\in \left[1,2\right[\), such that:
	\begin{multline}
		\forall\mathbf{v}\in W^{2,\infty}(\Omega;\mathbb{R}^2)\quad \text{such that: }\mathbf{v}_{|\Gamma_U}=\mathbf{0},\\
		\int_\Omega \boldsymbol\varepsilon(\mathbf{\tilde{u}}^0):\boldsymbol\Lambda\boldsymbol\varepsilon(\mathbf{v})\,{\rm d}x = \mathbf{\tilde{T}}\cdot \mathbf{v}(x') + \int_\Omega \boldsymbol \Phi(x,x')\cdot \mathbf{v}(x)\,{\rm d}x + \int_{\Gamma} \boldsymbol\tau(x,x')\cdot\mathbf{v}(x). 
		\label{eq:trucToProve}
	\end{multline}
\end{quote}
The key fact will be the observation that the displacement \(\mathbf{u}^0\) in Lemma~\ref{thm:fundSol2d} is in \(W^{1,p}(D;\mathbb{R}^2)\), for all \(p\in \left[1,2\right[\), and is \(C^\infty\) in \(\overline{D}\setminus\{0\}\).  We observe that the function \(\varphi\circ\phi^{-1}\) is a cut-off function on \(D\) which is identically equal to 1 on a neighborhood of \(\mathbf{0}\).  Denoting by \(\boldsymbol\Lambda^0\) the elastic modulus tensor of the isotropic homogeneous body of Young modulus \(E(x')\) and Poisson ratio \(\nu(x')\), we have:
\[
\forall\mathbf{v}\in W^{2,\infty}(\mathbb{R}^2;\mathbb{R}^2)\quad \text{with compact support},\qquad \int_{\{y<0\}} \boldsymbol\varepsilon(\mathbf{u}^0):\boldsymbol\Lambda^0\boldsymbol\varepsilon(\mathbf{v})\,{\rm d}x = \mathbf{T}\cdot \mathbf{v}(0) ,
\]
(where \(\mathbf{T}=\mathscr{R}\mathbf{\tilde{T}}\)).  Therefore, we have:
\begin{multline*}
	\forall\mathbf{v}\in W^{2,\infty}(\mathbb{R}^2;\mathbb{R}^2)\quad \text{with compact support},\qquad \\
	\int_{\{y<0\}} \!\boldsymbol\varepsilon\bigl((\varphi\circ\phi^{-1})\mathbf{u}^0\bigr):\boldsymbol\Lambda^0\boldsymbol\varepsilon(\mathbf{v})\,{\rm d}x = \mathbf{T}\cdot \mathbf{v}(0)+ \int_{\{y<0\}}\boldsymbol\Phi^0(x)\cdot\mathbf{v}(x)\,{\rm d}x + \int_{\{y=0\}}\boldsymbol\tau^0(x)\cdot\mathbf{v}(x),
\end{multline*}
for some functions \(\boldsymbol\Phi^0\in L^\infty\) and \(\boldsymbol\tau^0\in L^\infty\) that depend on \(\mathbf{u}^0\), \(\varphi\) and \(\phi\), where we have used the fact that $\varphi\circ\phi^{-1}\in W^{2,\infty}$ and it is $1$ in a neighborhood of $\mathbf{0}$.  Using \(\phi\) to perform a change of variable in the integrals of the above identity and applying the rotation $\mathscr R^{-1}$, we get by a standard calculation:
\begin{multline*}
	\forall\mathbf{v}\in W^{2,\infty}(\Omega;\mathbb{R}^2),\qquad \int_\Omega\boldsymbol\varepsilon(\mathbf{\tilde{u}}^0):\tilde{\boldsymbol\Lambda}\boldsymbol\varepsilon(\mathbf{v})\,{\rm d}x = \mbox{} \\
	\mathbf{\tilde{T}}\cdot \mathbf{v}(x')+ \int_{\Omega}\tilde{\boldsymbol\Phi}(x)\cdot\mathbf{v}(x)\,{\rm d}x + \int_{\Gamma}\tilde{\boldsymbol\tau}(x)\cdot\mathbf{v}(x),
\end{multline*}
where \(\tilde{\boldsymbol\Phi}\) and \(\tilde{\boldsymbol\tau}\) are \(L^{\infty}\) functions with support in \(\mathscr{V}_{x'}\).  The conveyed elastic modulus tensor \(\tilde{\boldsymbol\Lambda}\) satisfies the usual symmetry and ellipticity properties, but does not need to be isotropic homogeneous anymore.  However, it is a \(W^{1,\infty}\) function.  Note that we have \(\mathbf{\tilde{u}}^0\in W^{1,p}(\Omega;\mathbb{R}^2)\), for all \(p\in \left[1,2\right[\).  Finally, we can write:
\[
\boldsymbol\Lambda\boldsymbol\varepsilon(\mathbf{\tilde{u}}^0) = \tilde{\boldsymbol\Lambda}\boldsymbol\varepsilon(\mathbf{\tilde{u}}^0) +  \Bigl( \boldsymbol\Lambda\tilde{\boldsymbol\Lambda}^{-1} - \textbf{Id} \Bigr)\tilde{\boldsymbol\Lambda}\boldsymbol\varepsilon(\mathbf{\tilde{u}}^0).
\]
As \(\nabla \phi(x')\) is a rotation, \(\tilde{\boldsymbol\Lambda}(x')=\boldsymbol\Lambda(x')=\boldsymbol\Lambda^0\), and the second term of the right-hand side is the product of a \(W^{1,\infty}\) function which vanishes at \(x=x'\) and of an \(L^{p}\) (for all \(p\in \left[1,2\right[\)) stress field \(\tilde{\boldsymbol\sigma}:=\boldsymbol\varepsilon(\mathbf{\tilde{u}}^0)\) such that \(\text{div}\,\tilde{\boldsymbol\sigma}+\tilde{\boldsymbol\Phi}=\mathbf{0}\) in \(\Omega\) and \(\tilde{\boldsymbol\sigma}\mathbf{n}=\tilde{\mathbf{T}}\delta_{x'}+\tilde{\boldsymbol\tau}\) on \(\Gamma\).  This is sufficient to prove~\eqref{eq:trucToProve}. 

\medskip
\noindent\textbf{Step 3.} \textit{Regularity with respect to \(x'\) of the trace of \(\mathbf{g}(\cdot,x')\) on \(\Gamma_C\).  The function \(x'\mapsto \tilde{g}(\cdot,x')\in C^0(\overline{\Gamma}_C ; H^{1/2}(\Gamma_C))\).  In addition, the function \(x'\mapsto \mathbf{g}(\cdot,x') \in C^0(\overline{\Gamma}_C;W^{1,p}(\Omega;\mathbb{R}^2))\), for all \(p\in\left[1,2\right[\)}.

In Step~2, it was noted that the fundamental solution \(\mathbf{u}^0\) of Lemma~\ref{thm:fundSol2d} satisfies \(\mathbf{u}_{|D}^0\in W^{1,p}(D;\mathbb{R}^2)\), for all \(p\in\left[1,2\right[\).  Taking \(\bar{x}\) on the axis \(y=0\) in a neighborhood of \(\mathbf{0}\), this classically entails that \(\bar{x}\mapsto \mathbf{u}^0(\cdot-\bar{x})\in C^0(\left]-\epsilon,\epsilon\right[;W^{1,p}(D;\mathbb{R}^2))\).  Picking \(x'\in \Gamma:=\partial\Omega\setminus\overline{\Gamma}_U\) and the diffeomorphism \(\phi : \mathscr{V}_{x'} \rightarrow B\) as in Step~2, we consider an arbitrary \(x''\in \mathscr{V}_{x'}\cap \Gamma\) and \(\bar{x}:=\phi(x'')\).  Using the fundamental solution \(\mathbf{u}^0(\cdot-\bar{x})\) (with \(E(x'')\) and \(\nu(x'')\)) instead of \(\mathbf{u}^0(\cdot)\) in the reasoning made in Step~2, we can check that the functions \(x'\mapsto \boldsymbol\Phi(\cdot,x')\) and \(x'\mapsto \boldsymbol\tau(\cdot,x')\) are continuous functions of \(x'\) with values in \(L^p\) (\(p\in\left[1,2\right[\)).  Therefore, by~\eqref{eq:veryWeak} and \eqref{eq:trucToProve} and the standard regularity theory for elliptic boundary value problems, the mapping \(x'\mapsto \mathbf{g}(\cdot,x')-\mathbf{\tilde{u}}^0(\cdot,x')\) is continuous with respect to the \(H^1\) norm.  This is sufficient to obtain the conclusion of Step~3.

\medskip
\noindent\textbf{Step 4.} \textit{Conclusion.}

Consider an arbitrary \(t\in H^{-1/2}(\Gamma_C)\cap \mathscr{M}(\overline{\Gamma}_C)\).   Setting
\[
\mathbf{u}(x) := \int_{x'\in\Gamma_C} t\, \mathbf{g}(x,x'),\qquad \text{for }x\in\Omega,
\]
where \(\mathbf{g}(x,x')\) is the function defined in Step~1, we have \(\mathbf{u}\in W^{1,p}(\Omega;\mathbb{R}^2)\), for all \(p\in\left[1,2\right[\), by Step~3.  We have also:
\[
\forall \mathbf{v}\in C^\infty(\mathbb{R}^2;\mathbb{R}^2)\quad \text{with } \mathbf{v}_{|\Gamma_U}=\mathbf{0}, \qquad \int_\Omega \boldsymbol\varepsilon(\mathbf{u}):\boldsymbol\Lambda\boldsymbol\varepsilon(\mathbf{v})\,{\rm d}x = \int_{\Gamma_C} t v_t.
\]
This entails that \(\mathbf{u}\in H^1(\Omega;\mathbb{R}^2)\) and that the linear mapping \(t\mapsto \mathbf{u}\) is continuous with respect to the strong topologies of \(H^{-1/2}(\Gamma_C)\) and \(H^1(\Omega;\mathbb{R}^2)\).  Clearly \(\mathbf{u}\) is the unique minimizer on \(V\) of:
\[
\mathbf{v} \mapsto \frac{1}{2}\int_\Omega \boldsymbol\varepsilon(\mathbf{v}):\boldsymbol\Lambda\boldsymbol\varepsilon(\mathbf{v})\,{\rm d}x - \int_{\Gamma_C} t \,v_t.
\]
This concludes the proof. \qed

\bigskip
\noindent\textbf{Proof of Theorem~\ref{thm:cornerstone2D}.}

Let us consider a sequence \(t^k\in K\) (we recall that \(K\) is the nonpositive cone of \(H^{-1/2}(\Gamma_C)\)) converging weakly in \(H^{-1/2}(\Gamma_C)\) towards a limit \(t\in K\), and a sequence \(\tilde{t}^k\) in \(H^{-1/2}(\Gamma_C)\cap \mathscr{M}(\overline{\Gamma}_C)\) such that \(|\tilde{t}^k| \leq -t^k\) and converging weakly in \(H^{-1/2}(\Gamma_C)\) towards a limit \(\tilde{t}\).  Both \(t^k\) and \(\tilde{t}^k\) are bounded measures on \(\overline{\Gamma}_C\).

The proof is based on the expression of the operator \(L:t\in H^{-1/2}(\Gamma_C)\cap\mathscr{M}(\overline\Gamma_C)\mapsto u_n\) provided by Lemma~\ref{thm:represGreen}:
\[
\bigl\langle  t^k, L(\tilde{t}^k)\bigr\rangle = \int_{x\in \Gamma_C} t^k \int_{x'\in \Gamma_C} \Bigl[ c(x')\, h(x,x') +\tilde{g}(x,x')\Bigr]  \tilde{t}^k,
\]
where \(c\in W^{1,\infty}(\Gamma_C)\), \(h\) takes only the two values \(-1,+1\) and is continuous on \(\Gamma_C\times\Gamma_C\setminus \Delta\), with \(\Delta:=\{(x,x)|\:x\in\Gamma_C\}\) being the diagonal, and where \(x' \rightarrow \tilde{g}(\cdot,x')\in C^0(\overline{\Gamma}_C;H^{1/2}(\Gamma_C))\).  As \(t^k \in H^{-1/2}(\Gamma_C)\), the sequence of functions:
\[
x'\mapsto \int_{x\in \Gamma_C} \tilde{g}(x,x')\,  t^k, 
\]
is a sequence of continuous functions on \(\overline{\Gamma}_C\) that converges pointwisely towards the continuous function \(x'\mapsto \int_{x\in \Gamma_C} \tilde{g}(x,x')\, t\).  We are going to prove that the convergence is uniform on \(\overline{\Gamma}_C\).  As the function \(x' \rightarrow \tilde{g}(\cdot,x')\) is uniformly continuous on \(\overline{\Gamma}_C\), given \(\epsilon >0\), there exists a finite covering of $\overline{\Gamma}_C$ with balls $B_{x_i}$, \(x_i'\in \overline{\Gamma}_C\), such that:
\[
\forall x'\in B_{x_i}\cap \overline{\Gamma}_C,\qquad  \bigl\| \tilde{g}(\cdot,x') - \tilde{g}(\cdot,x_i')\bigr\|_{H^{1/2}(\Gamma_C)} < \epsilon. 
\]
Also, for large enough \(k\), we have:
\[
\forall i, \qquad \biggl|\int_{x\in \Gamma_C} \tilde{g}(x,x_i')\,t^k - \int_{x\in \Gamma_C} \tilde{g}(x,x_i')\, t \biggr| < \epsilon.
\]
Denoting by \(M>0\) an upper bound of the bounded sequence \(\| t^k\|_{H^{-1/2}(\Gamma_C)}\) and gathering all together, we obtain, for \(k\) large enough:
\[
\forall x'\in \overline{\Gamma}_C,\qquad \biggl|\int_{x\in \Gamma_C} \tilde{g}(x,x')\, t^k - \int_{x\in \Gamma_C} \tilde{g}(x,x')\, t \biggr| < (2M+1)\epsilon.
\]
This is sufficient to conclude that:
\[
\lim_{k\rightarrow+\infty}\int_{x'\in\Gamma_C}\tilde{t}^k\biggl(\int_{x\in \Gamma_C} \tilde{g}(x,x')\,t^k \biggr) = \int_{x'\in\Gamma_C}\tilde{t}\biggl(\int_{x\in \Gamma_C} \tilde{g}(x,x')\, t \biggr).
\]
Now, there remains to obtain the same conclusion with \(c(x')h(x,x')\) instead of \(\tilde{g}(x,x')\).  Replacing \(\tilde{t}^k\) by \(c\tilde{t}^k\), we may assume that \(c\equiv 1\) without loss of generality.  Then, the proof runs essentially as the one of Theorem~\ref{thm:cornerstoneHalfSpace}.
\begin{enumerate}
	\item The sequence \(t^k\) is a sequence which converges weakly in \(H^{-1/2}(\Gamma_C)\) and weakly-* in \(\mathscr{M}(\overline{\Gamma}_C)\) towards \(t\).  All these measures have no atoms in \(\overline\Gamma_C\).  Since \(x\mapsto h(x,x')\) has bounded variation on \(\overline{\Gamma}_C\), this entails the pointwise convergence of the following sequence of continuous functions:
	\[
	\forall x'\in \Gamma_C,\qquad \lim_{k \rightarrow+\infty} \int_{x\in \Gamma_C} h(x,x')\, t^k = \int_{x\in \Gamma_C} h(x,x')\, t,
	\]
	as in the proof of Theorem~\ref{thm:cornerstoneHalfSpace}.
	\item The pointwise convergence is actually uniform on \(\overline{\Gamma}_C\).  This fact relies on the two following properties. 
	\begin{itemize}
		\item  The pointwise limit \(x'\mapsto \int_{x\in \Gamma_C} h(x,x')\, t\), being continuous on the compact \(\overline{\Gamma}_C\), is uniformly continuous.
		\item  The same reasoning as in Theorem~\ref{thm:cornerstoneHalfSpace} shows that the sequence \(x'\mapsto \int_{x\in \Gamma_C} h(x,x')\, t^k\) converges strongly in \(C^0(\overline{\Gamma}_C)\) towards \(x'\mapsto \int_{x\in \Gamma_C} h(x,x')\, t\).   More precisely, as \(\Gamma_C\) is an open subset of the smooth curve \(\partial\Omega\), it has countably many connected components \(I_i\).  If there are infinitely many \(I_i\)'s, then, given \(\epsilon >0\), there exists \(\bar{n}>0\) such that:
		\[
		-\sum_{i=\bar{n}}^{+\infty}\int_{x\in I_i}t < \epsilon /3,
		\]
		and therefore, there exists \(\bar k>0\) such that:
		\[
		\forall k\geq \bar k,\qquad -\sum_{i=\bar{n}}^{+\infty}\int_{x\in I_i}t^k < 2\epsilon /3.
		\]
		Hence,
		\begin{multline*}
			\forall k\geq \bar k,\quad \forall x'\in \overline{\Gamma}_C,\qquad \\
			\Biggl|\int_{x\in \bigcup_{i=\bar{n}}^{\infty}I_i} h(x,x')\, t^k - \int_{x\in \bigcup_{i=\bar{n}}^{\infty}I_i} h(x,x')\, t\Biggr| \leq -\sum_{i=\bar{n}}^{+\infty}\int_{I_i}t^k - \sum_{i=\bar{n}}^{+\infty}\int_{I_i}t < \epsilon,
		\end{multline*}
		as \(h\) takes values in \(\{-1,+1\}\) and \(t^k\leq 0\).  Therefore, it is now sufficient to prove that the sequence of functions \(x'\mapsto \int_{x\in I_0}h(x,x')t^k\) converges uniformly on \(\overline{\Gamma}_C\) to obtain the expected conclusion.  As,
		\[
		\forall x'\in \overline{\Gamma}_C\setminus I_0, \qquad \int_{x\in I_0}h(x,x')t^k - \int_{x\in I_0}h(x,x')t = \pm \int_{x\in I_0}t^k -t,
		\]
		we only have to prove that the sequence of functions \(x'\mapsto \int_{x\in I_0}h(x,x')t^k\) converges uniformly on \(\overline{I}_0\).   As \(I_0\) is an open connected subset of the smooth curve \(\partial\Omega\), it can be parametrized by the arc-length denoted by \(s\in \left]0,l\right[\).   With this notation, we have either \(h(x,x')=\text{sgn}(s-s')\) or \(h(x,x')=-\text{sgn}(s-s')\) and it is now sufficient to prove that the sequence of functions \(s'\mapsto \int_{s\in \left]0,l\right[}\text{sgn}(s-s')\,t^k\) converges uniformly on \([0,l]\).  But the proof of this fact was already given in Theorem~\ref{thm:cornerstoneHalfSpace}.
	\end{itemize}
\end{enumerate}
Hence, we can conclude that:
\[
\lim_{k\rightarrow+\infty}\int_{x'\in\Gamma_C}\tilde{t}^k\biggl(\int_{x\in \Gamma_C} h(x,x')\,  t^k \biggr) = \int_{x'\in\Gamma_C}\tilde{t}\biggl(\int_{x\in \Gamma_C} h(x,x')\,  t \biggr),
\]
and that the same result holds true with \(h(x,x')\) replaced by \(c(x')h(x,x')+\tilde{g}(x,x')\).  Recalling the definition of \(c,h,\tilde g\) in the statement of Lemma~\ref{thm:represGreen}, it turns out that the above equality is nothing but the conclusion of Theorem~\ref{thm:cornerstone2D} for the sequences \( t^k\) and \(\tilde{t}^k\). \qed

\section{The case of anisotropic elasticity}

\label{sec:anisotropy}

In Section~\ref{sec:isotropy}, the formal time-incremental Signorini-Coulomb problem~\eqref{eq:SignoriniCoulombDiscret} was reformulated as a variational inequality in terms of a nonlinear operator \(A\).  Then, the path followed to prove the existence of a solution to that variational inequality, and therefore to problem~\eqref{eq:SignoriniCoulombDiscret}, was to prove that this operator is coercive and Leray-Lions, that is, fulfills the four properties defining a Leray-Lions operator (see Definition~\ref{thm:defLerayLions} in Appendix~A).   Coercivity and the first three properties were easily checked to hold in any space dimension and any case of isotropic or anisotropic heterogeneous elasticity.  However, the checking of the fourth property defining a Leray-Lions operator turned out to be very tricky and required to prove a new fine property of the elastic Neumann-to-Dirichlet operator.  The fine property was proved to hold true for the 2D case (\(N=2\)) in the case of isotropic elasticity (either homogeneous or heterogeneous).  In this section, we discuss completely the extension to anisotropic elasticity (still in the 2D case).  Surprisingly, the picture turns out to be different from the isotropic case as a condition on the friction coefficient for the solvability  of the Signorini-Coulomb problem shows up: the friction coefficient must be smaller than a critical value.  We show that our condition is optimal by displaying a proof of nonexistence of a solution for friction coefficients violating the condition.  This surprising feature of the coupling between friction and anisotropic elasticity seems to have remained unnoticed up to now.

This section is organized as follows.  As the structure of the Neumann-to-Dirichlet operator of the elastic half-space was shown in Subsection~\ref{sec:fineProperty} to play a crucial role in the proof of the fine property for the bounded body of arbitrary geometry, we first provide the structure of the Neumann-to-Dirichlet operator of the anisotropic 2D half-space at the beginning of Subsection~\ref{sec:nonexistence} (Theorem~\ref{thm:anisotropicD2N}).  This is then exploited to prove that the steady-sliding frictional contact problem for the anisotropic half-space can have no solution in the case of large friction coefficients, whereas this problem has a unique solution for arbitrarily large friction coefficients in the isotropic case.  The optimal condition on the friction coefficient is supplied.  In Subsection~\ref{sec:existAnisotrop}, the solvability of the general anisotropic 2D Signorini-Coulomb problem is proved provided that this condition is fulfilled.  The proof of Theorem~\ref{thm:anisotropicD2N}, being technical, is postponed to Subsection~\ref{sec:N2Danisotrop}.

\subsection{A frictional contact problem on the half-space with no solution}

\label{sec:nonexistence}

In this section, \(\boldsymbol\Lambda:\mathbb{M}^{2\times 2}_\text{sym}\to \mathbb{M}^{2\times 2}_\text{sym}\) denotes an arbitrary linear operator, satisfying requirements~\eqref{eq:reqLamda1} and \eqref{eq:reqLambda2}, and constant (independent of \(x\)). 

The corresponding Neumann-to-Dirichlet operator of the 2D half-space has the structure given by the following theorem, whose proof is postponed to Section~\ref{sec:N2Danisotrop}.  We will identify the boundary $\mathbb R\times \{0\}$ of the half-space $\mathbb R\times\mathbb{R}^-$ with $\mathbb R$.  We will also use the indices $n$ and $t$ to denote the normal and tangential components of distributions defined on $\mathbb R\times \{0\}$ with values in $\mathbb R^2$.

\begin{theo}
	\label{thm:anisotropicD2N}
	Let \(\mathbf{t}\in H^{-1/2}(-1,1;\mathbb{R}^2)\) be a compactly supported surface traction distribution, prescribed at the boundary of the 2D half-space.  Then, all the tempered distributional displacements $\mathbf{u}:\mathbb R\times\mathbb{R}^-\to\mathbb{R}^2$ that satisfy the anisotropic elastic equilibrium equations with vanishing body forces are in \(H_{\rm loc}^1(\mathbb R\times\mathbb{R}^-;\mathbb R^2)\).  They therefore have a trace on the boundary of the half-space.  The corresponding surface displacement reads as:
	\begin{align*}
		u_n & = -\mathscr{C}_1\,\log|x|*t_n - \mathscr{C}_2\,\log|x|*t_t - \mathscr{C}_3\,\text{\rm sgn}(x)*t_t+a_n,\\
		u_t & = - \mathscr{C}_4\,\log|x|*t_t- \mathscr{C}_2\,\log|x|*t_n + \mathscr{C}_3\,\text{\rm sgn}(x)*t_n+a_t,
	\end{align*}
	for \(x\in\mathbb{R}\), where \(\mathscr{C}_1,\mathscr{C}_2,\mathscr{C}_3,\mathscr{C}_4\) are uniquely determined real constants that are independent of \(\mathbf{t}\), and \(\mathbf{a}\) is an affine function that depends on four arbitrary real constants (three of these constants correspond to an overall rigid motion and the last one to a component of stress at infinity).  The constants \(\mathscr{C}_1\) and \(\mathscr{C}_4\) are positive.  Finally, the mapping \(\mathbf{\Lambda} \mapsto (\mathscr{C}_1,\mathscr{C}_2,\mathscr{C}_3,\mathscr{C}_4)\) is of class \(C^\infty\).
\end{theo}

In the particular case of isotropic elasticity,
\[
\mathscr{C}_1=\mathscr{C}_4 = \frac{2(1-\nu^2)}{\pi E},\qquad \mathscr{C}_2=0,\qquad \mathscr{C}_3=\frac{(1-2\nu)(1+\nu)}{2E},
\]
where \(\nu\in\left]-1,1/2\right[\) is the Poisson ratio and \(E>0\) the Young modulus.  The key difference between isotropic and anisotropic elasticity is the fact that \(\mathscr{C}_2\) can be nonzero in some cases of anisotropic elasticity.

To demonstrate the impact of a possibly nonzero \(\mathscr{C}_2\) on the existence of solutions to frictional contact problems, we consider the frictional contact problem raised by a moving indentor steadily sliding along the boundary of the half-space with given velocity \(w\neq 0\).  This problem was already studied in the particular case of isotropic elasticity in \cite{BallardIurlano}.  The shape of the indentor is represented by a given function \(g:\left]-1,1\right[\to\mathbb{R}\) so that the contact conditions in a frame moving with the indentor read as:
\[
t_n \leq 0,\qquad u_n\leq g,\qquad t_n (u_n-g)=0.
\]
As the indentor slides steadily along the boundary of the half-space, the Coulomb friction law reduces to the following linear condition:
\[
t_t = - \text{\rm sgn}(w)\,f\,t_n.
\]
In the sequel, the friction coefficient \(f\in \mathbb{R}^+\) will be assumed to be constant, that is, independent of \(x\).  The arbitrary affine function \(\mathbf{a}\) in the Neumann-to-Dirichlet operator plays no role and will be taken as zero.   The steady sliding frictional contact problem is now that of finding \(u_n,t_n:\left]-1,1\right[\to\mathbb{R}\) such that:
\begin{equation*}
	\left\{
	\begin{aligned}
		-\bigl(\mathscr{C}_1-\text{\rm sgn}(w)\,f\,\mathscr{C}_2\bigr)\;\log|x|*t_n  + \text{\rm sgn}(w)\,f\,\mathscr{C}_3\; \text{\rm sgn}(x) * t_n = u_n, & \qquad \text{in }\left]-1,1\right[,\\
		\quad  t_n \leq 0,\qquad u_n\leq g,\qquad t_n (u_n-g)=0, & \qquad\text{in }\left]-1,1\right[.
	\end{aligned}
	\right.
\end{equation*}
where the convolution product is understood in terms of the extension by zero of \(t_n\) to the whole line.  The above formal contact problem is going to be formulated as a variational inequality.  The following proposition, which we recall here from \cite{bibiJiri}, provides the appropriate functional framework.

\begin{prop}
	\label{thm:BallardJarusek}
	For \(t\in H^{-1/2}(-1,1)\), we have \(\log|x|*t\in H_{\rm loc}^{1/2}(\mathbb{R})\) and \(\text{\rm sgn}(x)*t\in H_{\rm loc}^{1/2}(\mathbb{R})\), so that the bilinear forms:
	\[
	S(t_1,t_2):=   -\bigl\langle \log|x| * t_1, t_2  \bigr\rangle,\qquad  A(t_1,t_2):=\bigl\langle \text{\rm sgn}(x) * t_1 , t_2  \bigr\rangle,
	\]
	are well-defined on \(H^{-1/2}(-1,1)\).  They are also continuous on \(H^{-1/2}(-1,1)\).  The bilinear form \(A\) is skew-symmetric.  The bilinear form \(S\) is symmetric and positive definite.  It therefore induces a norm on \(H^{-1/2}(-1,1)\) and this norm is equivalent to that of \(H^{-1/2}(-1,1)\).
\end{prop}

\noindent\textbf{Proof.} For \(\varphi\in L^1(\mathbb{R})\), we define the Fourier transform of \(\varphi\) as:
\begin{equation}
	\label{eq:convFourier}
	\mathscr{F}[\varphi](\xi) := \frac{1}{\sqrt{2\pi}}\int_{-\infty}^{\infty} e^{-i\xi x}\,\varphi(x)\,{\rm d}x.
\end{equation}
With this convention, the distributional Fourier transforms of the locally integrable functions \(\log|x|\) and \(\text{\rm sgn}(x)\) read as:
\begin{equation}
	\label{eq:FourierTransform}
	\mathscr{F}\Bigl[\log|x|\Bigr] = -\sqrt{\frac{\pi}{2}}\,\text{\rm fp}\frac{1}{|\xi|} - \sqrt{2\pi}\,\Gamma_{\rm Eul}\,\delta,\qquad \qquad 
	\mathscr{F}\Bigl[\text{\rm sgn}(x)\Bigr] = - i \sqrt{\frac{2}{\pi}}\,\text{\rm pv}\frac{1}{\xi},
\end{equation}
where \(\delta\) is the Dirac measure at \(0\), \(\Gamma_{\rm Eul}:=\lim_{n\rightarrow\infty}(-\log n+\sum_{k=1}^n 1/k)\) is the Euler-Mascheroni constant, and \(\text{\rm fp}1/|x|\) and \(\text{\rm pv}1/x\) stand for the distributional derivatives of \(\text{\rm sgn}(x)\,\log|x|\) and \(\log|x|\), respectively (`fp' stands for `finite part' and `pv' for `principal value').  Identifying \(H^{-1/2}(-1,1)\) with the space of elements in \(H^{-1/2}(\mathbb{R})\) whose support is contained in \([-1,1]\), and recalling the definition of \(H^{-1/2}(\mathbb{R})\) in terms of the Fourier transform:
\[
H^{-1/2}(\mathbb{R}) := \Bigl\{ t\in \mathscr{S}'(\mathbb{R})\;\bigm|\;  |\mathscr{F}[t](\xi)|/(1+|\xi|^2)^{1/4}\in L^2(\mathbb{R})\Bigr\},
\]
where \(\mathscr{S}'\) denotes the space of tempered distributions, we can check that the mapping \(t\mapsto (\log|x|*t)'\) is continuous on \(H^{-1/2}(\mathbb{R})\).  Hence, \(\log|x|*t\in H_{\rm loc}^{1/2}(\mathbb{R})\) and \(\text{\rm sgn}(x)*t\in H_{\rm loc}^{1/2}(\mathbb{R})\) whenever \(t\in H^{-1/2}(-1,1)\), and the bilinear forms \(S\) and \(A\) are well-defined and continuous on \(H^{-1/2}(-1,1)\).   As the function \(\log|x|\) is even and \(\text{\rm sgn}(x)\) is odd, one has that \(S\) is symmetric and that \(A\) is skew-symmetric.  Considering the restriction of S to the codimension 1 subspace:
\[
\mathscr{H} := \Bigl\{t\in H^{-1/2}(-1,1)\;\bigm|\; \langle t,1\rangle = 0\Bigr\},
\]
where \(1\) denotes the function in \(H^{1/2}(-1,1)\) that takes the constant value \(1\) all over \(\left]-1,1\right[\), we have by Plancherel's formula:
\[
\forall t\in \mathscr{H},\qquad S(t,t) = \pi \int_{-\infty}^{\infty}\frac{|\mathscr{F}[t](\xi)|^2}{|\xi|}\,{\rm d}\xi \geq \pi\int_{-\infty}^{\infty}\frac{|\mathscr{F}[t](\xi)|^2}{\sqrt{1+|\xi|^2}}\,{\rm d}\xi,
\]
where the first integral is well-defined as the Fourier transform of any \(t\in \mathscr{H}\) is a \(C^\infty\) function which vanishes at \(\xi=0\) and the second integral can be identified with \(\|t\|_{H^{-1/2}}^2\).  Hence, \(S\) is a scalar product on \(\mathscr{H}\) and, by the open mapping theorem, it induces a norm on \(\mathscr{H}\) that is equivalent to that of \(H^{-1/2}(-1,1)\).  To get the final conclusion, it remains only to prove that \(S\) is positive definite on \(H^{-1/2}(-1,1)\).  Taking an arbitrary \(t\in H^{-1/2}(-1,1)\), we have the decomposition:
\[
t = \biggl(t-\frac{\langle t,1\rangle}{\pi\sqrt{1-x^2}}\biggr) + \frac{\langle t,1\rangle}{\pi\sqrt{1-x^2}},
\]
where the first term belongs to \(\mathscr{H}\).  Using (see, for example, \cite[Step 3 in proof of Theorem 3]{bibiJiri})
\[
\forall x\in \left]-1,1\right[,\qquad \int_{-1}^1 \frac{\log|s-x|}{\sqrt{1-s^2}}{\rm d}s = -\pi\log 2,
\]
we get:
\[
S(t,t) = S\biggl(t-\frac{\langle t,1\rangle}{\pi\sqrt{1-x^2}},t-\frac{\langle t,1\rangle}{\pi\sqrt{1-x^2}}\biggr) + \log2\;\langle t,1\rangle^2,
\]
which shows that \(S\) is positive definite on \(H^{-1/2}(-1,1)\), indeed. \qed

\bigskip
The following consequence is immediate.

\begin{coro}
	\label{thm:propCoerFormBil}
	The following bilinear form:
	\[
	a(t_1,t_2):= \bigl(\mathscr{C}_1-\text{\rm sgn}(w)\,f\,\mathscr{C}_2\bigr)\;S(t_1,t_2) + \text{\rm sgn}(w)\,f\,\mathscr{C}_3\; A(t_1,t_2),
	\]
	is continuous on \(H^{-1/2}(-1,1)\).  If \(\mathscr{C}_1-\text{\rm sgn}(w)\,f\,\mathscr{C}_2>0\), then \(a\) is coercive on \(H^{-1/2}(-1,1)\).  If \(\mathscr{C}_1-\text{\rm sgn}(w)\,f\,\mathscr{C}_2<0\), then \(-a\) is coercive on \(H^{-1/2}(-1,1)\).
\end{coro}

Under the hypothesis \(g\in H^{1/2}(-1,1)\), the steady sliding frictional contact problem can now be precisely formulated as a variational inequality.

\medskip
\noindent\textbf{Problem }$\mathscr{P}$.  Find \(t_n\in H^{-1/2}(-1,1)\), such that \(t_n\leq 0\), and:
\[
\forall \hat{t}\in H^{-1/2}(-1,1),\;\text{such that }\hat{t}\leq 0,\qquad a(t_n,\hat{t}-t_n) \geq \bigl\langle g,\hat{t}-t_n\bigr\rangle.
\]

When condition \(\mathscr{C}_1-\text{\rm sgn}(w)\,f\,\mathscr{C}_2>0\) is fulfilled, the bilinear form \(a\) is continuous and coercive and the Lions-Stampacchia theorem \cite[Theorem 4.4]{Troianiello} ensures that problem~$\mathscr{P}$ has one and only one solution.  This happens in particular for any value of the friction coefficient \(f\) in the case of isotropic elasticity (as \(\mathscr{C}_1>0\) and \(\mathscr{C}_2=0\) in that case).  However, for cases of anisotropic elasticity where \(\mathscr{C}_2\neq 0\), one can have \(\mathscr{C}_1-\text{\rm sgn}(w)\,f\,\mathscr{C}_2<0\) whenever the friction coefficient is large enough.  In that case, the bilinear form \(-a\) is coercive and there are many examples of smooth shapes \(g\) of the indentor for which problem~$\mathscr{P}$ has no solution, as seen in the next proposition.

\begin{prop}
	Let \(H\) be a Hilbert space, \(K\subset H\) a closed convex cone such that \(K\setminus\{\mathbf{0}\}\neq\varnothing\) and \(K\cap (-K)=\{\mathbf{0}\}\).  Let also \(a:H\times H\to\mathbb{R}\) be a continuous bilinear form such that \(-a\) is coercive and let \(\tilde{g}\in -K\setminus\{\mathbf{0}\}\).
	
	Then, there is no \(t\in K\) such that:
	\[
	\forall \hat{t}\in K,\qquad a(t-\tilde{g},\hat{t}-t)\geq 0.
	\]
\end{prop}

\noindent\textbf{Proof.} Assume that there is such \(t\in K\).  As \(-\tilde{g}\in K\) and \(K\) is a cone, \(t-(\tilde{g}-t)\in K\).  Taking \(t-(\tilde{g}-t)\) as a test function in the variational inequality, we get:
\[
a(\tilde{g}-t,\tilde{g}-t) \geq 0.
\]
As \(\tilde{g}\notin K\), we have \(\tilde{g}\neq t\), and we obtain a contradiction since \(-a\) is coercive. \qed

\bigskip
As a consequence, we get an optimal condition on the friction coefficient to ensure the solvability of the steady sliding frictional contact problem.

\begin{coro}
	If \(\mathscr{C}_1-f|\mathscr{C}_2|>0\), then problem~$\mathscr{P}$ has a unique solution for all indentor shapes \(g\in H^{1/2}(-1,1)\).  If \(\mathscr{C}_1-f|\mathscr{C}_2|\leq 0\), then there are smooth indentor shapes \(g\) for which problem~$\mathscr{P}$ has no solution.
\end{coro}

\noindent\textbf{Proof.} The case \(\mathscr{C}_1-f|\mathscr{C}_2| \neq 0\) follows from the previous discussion.  In the case  \(\mathscr{C}_1-f|\mathscr{C}_2| = 0\), the bilinear form is skew-symmetric, which easily gives examples of nonexistence. \qed

\subsection{Existence in the case of anisotropic elasticity}

\label{sec:existAnisotrop}

As seen in Theorem~\ref{thm:anisotropicD2N}, the normal component of the surface displacement induced by a compactly supported surface traction \((t_n,t_t)\) applied at the boundary of a (possibly anisotropic) homogeneous elastic 2D half-space is given by:
\[
u_n = -\mathscr{C}_1 \log|x|*\bigl(t_n + t_t \,\mathscr{C}_2/\mathscr{C}_1\bigr) - \mathscr{C}_3 \;\text{\rm sgn}(x) * t_t,
\]
where \(\mathscr{C}_1>0\), \(\mathscr{C}_2\) and \(\mathscr{C}_3\) are constants that depend only on the elastic modulus operator \(\boldsymbol\Lambda\).  Setting \(\alpha:=\mathscr{C}_2/\mathscr{C}_1\), we have:
\[
u_n = -\mathscr{C}_1 \log|x|*\bigl(t_n + \alpha\,t_t\bigr) - \mathscr{C}_3\; \text{\rm sgn}(x) * t_t,
\]
so that the difference between isotropic and anisotropic elasticity lies only in the possibility of nonzero \(\alpha\).  It was seen in Section~\ref{sec:nonexistence} that nonzero \(\alpha\) has a dramatic influence on the existence of solutions to the frictional contact problem.

To be able to adapt Theorem~\ref{thm:cornerstone2D} to the anisotropic case, we need to have \(t_n\) replaced by \(t_n+\alpha t_t\) in the Coulomb friction law.  This drives us to use a nonorthogonal local basis on \(\Gamma_C\) instead of the `natural' basis \((\mathbf{n},\boldsymbol\tau)\).  This is the motivation of the following definitions.  We consider the following primal and dual basis:
\begin{align*}
	\mathbf{e}_1 & := \mathbf{n}+\alpha\boldsymbol\tau,\qquad & \mathbf{e}_2 & := \boldsymbol\tau,\\
	\mathbf{e}^1 & := \mathbf{n},\qquad & \mathbf{e}^2 & := \boldsymbol\tau-\alpha\mathbf{n},
\end{align*}
which clearly fulfills \(\mathbf{e}_i\cdot\mathbf{e}^j = \delta_i^j\), with \(\delta_i^j\) denoting the Kronecker symbol.  Using this pair of bases, we have:
\[
\mathbf{u} = u^1\mathbf{e}_1+u^2\mathbf{e}_2,\qquad \mathbf{t} = t_1\mathbf{e}^1+t_2\mathbf{e}^2,\qquad \mathbf{t}\cdot\mathbf{u} = t_1 u^1+t_2 u^2,
\]
with:
\begin{align*}
	u^1 & := \mathbf{u}\cdot\mathbf{e}^1 = u_n,\qquad & u^2 & := \mathbf{u}\cdot\mathbf{e}^2 = u_t-\alpha u_n,\\
	t_1 & := \mathbf{t}\cdot\mathbf{e}_1 = t_n +\alpha t_t,\qquad & t_2 & := \mathbf{t}\cdot\mathbf{e}_2 = t_t.
\end{align*}
Here, we emphasize that these identities are valid with \(\alpha\) possibly varying with \(x\in \Gamma_C\).  This will be the case as \(\alpha\) is determined by \(\boldsymbol\Lambda(x)\) and heterogeneous material will be possibly considered.

The aim of the following proposition is to reformulate the unilateral contact condition and the friction law in terms of \(t_1,t_2,u^1,u^2\) instead of \(t_n,t_t,u_n,u_t\).

\begin{prop}
	\label{thm:equivFrictionLaw}
	Let $f>0$, $g,\alpha\in \mathbb{R}$, and \(\mathbf{t},\mathbf{u}\in \mathbb{R}^2\).  We assume that \(0\leq |\alpha|f<1\).  Then, we have $\textit{(i)}\iff\textit{(ii)}$ with:
	\begin{align*}
		\textit{(i)} & \left\{
		\begin{array}{l}
			u_n - g \leq 0, \qquad t_n \leq 0, \qquad t_n\,(u_n-g) = 0, \\[1.5ex]
			\forall \hat{v}\in\mathbb{R},\qquad t_t\bigl(\hat{v}-u_t\bigr) - ft_n \bigl(|\hat{v}|-|u_t|\bigr) \geq 0,
		\end{array}
		\right. \\[1.5ex]
		\textit{(ii)} & \left\{
		\begin{array}{l}
			u^1 - g \leq 0, \qquad t_1 \leq 0, \qquad t_1\,(u^1-g) = 0, \\[1.5ex]
			\forall \hat{v}\in\mathbb{R},\qquad t_2\Bigl[\hat{v}-(u^2+\alpha g)\Bigr] \\[1.2ex]
			\displaystyle
			\hspace{3cm}\mbox{}- \frac{ft_1}{1+\alpha f}\Bigl[\langle\hat{v}\rangle^+ - \langle u^2+\alpha g\rangle^+\Bigr] - \frac{ft_1}{1-\alpha f}\Bigl[\langle\hat{v}\rangle^- - \langle u^2+\alpha g\rangle^-\Bigr]\geq 0,
		\end{array}
		\right.
	\end{align*}
	where \(\langle x\rangle^+:=\max\{0,x\}\) and \(\langle x\rangle^-:=\max\{0,-x\}\) denote respectively the positive and the negative parts.
\end{prop}

\begin{rem}
	The second line of (ii) is equivalent to the following conditions
	\begin{gather*}
		\frac{ft_1}{1+\alpha f} \leq t_2 \leq -\frac{ft_1}{1-\alpha f},\\
		u^2+\alpha g>0  \implies 	t_2=\frac{ft_1}{1+\alpha f},\;\mbox{ }\\
		u^2+\alpha g<0  \implies 	t_2=-\frac{ft_1}{1-\alpha f}.
	\end{gather*}
\end{rem}

\noindent\textbf{Proof.}

\noindent$\textit{(i)}\implies\textit{(ii)}$

Assuming that \(\textit{(i)}\) is fulfilled, we obviously have \(u^1- g\leq 0\) as \(u^1=u_n\).   By \textit{(i)}, we have \(|t_t|\leq -ft_n\) which entails that \(t_1=t_n+\alpha t_t\leq 0\), since \(0\leq |\alpha|f<1\).  If \(u^1=u_n<g\) then \(t_n=t_t=0\) by \textit{(i)}, so that \(t_1=0\).  The first line of \textit{(ii)} is proved.  Furthermore,
\[
ft_n<t_t<-ft_n \iff \frac{f(t_n+\alpha t_t)}{1+\alpha f}<t_t<-\frac{f(t_n+\alpha t_t)}{1-\alpha f}\iff  \frac{ft_1}{1+\alpha f}<t_2<-\frac{ft_1}{1-\alpha f}.
\]
If one, and then, each of the previous conditions is satisfied, we have \(u_n-g=u_t=0\) by \textit{(i)}, so that \(u^2+\alpha g=0\).  If \(t_n<0\) and $t_t=-ft_n$, that is, $t_2=-ft_1/(1-\alpha f)$, then \(u_n=g\) and \(u^2+\alpha g = u_t>0\).  Analogously, if \(t_n<0\) and $t_t=ft_n$, then $t_2=ft_1/(1+\alpha f)$ and \(u^2+\alpha g = u_t<0\).  Finally, the second line of \textit{(ii)} is fulfilled.

\noindent$\textit{(ii)}\implies\textit{(i)}$

The proof follows the same lines as that of 
\noindent$\textit{(i)}\implies\textit{(ii)}$ and is left to the reader.\qed

\bigskip
\begin{defi}
	\label{defAlpha}
	Let $\Omega\subset \mathbb R^N$ be open, bounded, connected and of class $C^{1,1}$.  Consider an arbitrary \(\mathbf{\Lambda}\in W^{1,\infty}(\Omega)\) satisfying requirements~\eqref{eq:reqLamda1} and \eqref{eq:reqLambda2}.  Taking an arbitrary \(x\in\partial\Omega\), one can always consider the Neumann-to-Dirichlet operator for the \emph{homogeneous} elastic 2D half-space with normal \(\mathbf{n}(x)\) and with elastic modulus taken as \(\mathbf{\Lambda}(x)\).  Theorem~\ref{thm:anisotropicD2N} provides associated functions \(\mathscr{C}_1(x),\mathscr{C}_2(x)\in W^{1,\infty}(\partial\Omega)\), with \(\mathscr{C}_1(x)>0\).  Setting \(\alpha(x):=\mathscr{C}_2(x)/\mathscr{C}_1(x)\), we have a well-defined function \(\alpha\in W^{1,\infty}(\partial\Omega)\) which is completely determined by \(\mathbf{\Lambda}\in W^{1,\infty}(\Omega)\).   In the case of an isotropic (possibly heterogeneous) material, the function \(\alpha\) vanishes identically on \(\partial\Omega\).
\end{defi}

We assume \(N=2\) (2D problem).  For \(\mathbf{\Lambda}\in W^{1,\infty}(\Omega)\) and \(f\in W^{1,\infty}(\Gamma_C)\) satisfying the conditions:
\begin{equation}
	\label{eq:condAnisotropic}
	f\geq 0,\qquad \sup_{x\in \Gamma_C} f(x)|\alpha(x)| < 1,
\end{equation}
and for \(t,\tau\in K\) (\(K\) was defined in formula~\eqref{eq:defK}), we define the modified energy:
\begin{multline*}
	E_{t,\tau}^\alpha(\mathbf{v}) := \frac{1}{2}\int_\Omega \boldsymbol\varepsilon(\mathbf{v}):\boldsymbol\Lambda\boldsymbol\varepsilon(\mathbf{v})\,{\rm d}x - \int_\Omega \mathbf{F}\cdot\mathbf{v} \,{\rm d}x - \int_{\Gamma_T}\mathbf{T}\cdot\mathbf{v}_{|\Gamma_T} \,{\rm d}x- \int_{\Gamma_C} \tau\, v_n  \\
	\mbox{} - \int_{\Gamma_C} \frac{ft}{1+\alpha f}\, \bigl\langle v_t - w_t - \alpha v_n +\alpha g\bigr\rangle^+ - \int_{\Gamma_C} \frac{ft}{1-\alpha f}\, \bigl\langle v_t - w_t - \alpha v_n +\alpha g\bigr\rangle^-,
\end{multline*}
The functional \(E_{t,\tau}^\alpha\) is clearly strictly convex and continuous on \(V\) (the space \(V\) was defined in formula~\eqref{eq:defV}).  Its unique minimizer \(\mathbf{u}\) on \(V\) satisfies, in particular, \(t_1=t_n+\alpha t_t=\tau\) on \(\Gamma_C\) and:
\begin{multline*}
	\forall \hat{v}\in H^{1/2}(\Gamma_C),\qquad \int_{\Gamma_C} t_t\bigl(\hat{v}-u_t+w_t+\alpha u_n-\alpha g\bigr)
	\\
	\mbox{} - \int_{\Gamma_C} \frac{ft}{1+\alpha f}\, \Bigl(\bigl\langle\hat{v}\bigr\rangle^+- \bigl\langle u_t - w_t - \alpha u_n +\alpha g\bigr\rangle^+\Bigr) - \int_{\Gamma_C} \frac{ft}{1-\alpha f}\, \Bigl(\bigl\langle\hat{v}\bigr\rangle^-- \bigl\langle u_t - w_t - \alpha u_n +\alpha g\bigr\rangle^-\Bigr)\geq 0,
\end{multline*}
where \(\mathbf{t}:=\boldsymbol\sigma(\mathbf{u})\mathbf{n}\) denotes the associated surface traction on \(\Gamma_C\).  To obtain the above variational inequality, we have decomposed displacements and traction along the bases $\{e_i\}$ and $\{e^j\}$ and used the change of variables $\hat v=v_t-\alpha v_n+\alpha g-w_t$.  We define \(\mathscr{A}^\alpha(t,\tau):=\mathbf{u}_{|\Gamma_C}\cdot\mathbf{n}=u^1\) as the normal part of the trace of \(\mathbf{u}\) on \(\Gamma_C\) and \(A^\alpha t:=\mathscr{A}^\alpha(t,t)\), so that:
\[
\mathscr{A}^\alpha: K\times K \rightarrow H^{1/2}(\Gamma_C),\qquad A^\alpha: K \rightarrow H^{1/2}(\Gamma_C).
\]
Thanks to Proposition~\ref{thm:equivFrictionLaw}, the time-incremental Signorini-Coulomb problem~\eqref{eq:SignoriniCoulombDiscret} reduces, under condition~\eqref{eq:condAnisotropic}, to find \(t\in K\) such that:
\[
\forall \hat{t}\in K,\qquad \bigl\langle A^\alpha t -g \,,\, \hat{t}-t \bigr\rangle \geq 0,
\]
where \(g\in H^{1/2}(\Gamma_C)\).

Hence, the problem is reduced to prove that \(A^\alpha: K \rightarrow H^{1/2}(\Gamma_C)\) is coercive and is a Leray-Lions operator.  This is going to be done along the same lines as in the isotropic case in Section~\ref{sec:existProofIsotropy}.  It is readily checked that the conclusions of Propositions~\ref{thm:boundedProp}, \ref{thm:strongMonotProp}, \ref{thm:contProp}, \ref{thm:propiii} and \ref{thm:propCoer} still holds true with the operator \(\mathscr{A}\) replaced with \(\mathscr{A}^\alpha\).

The key reason for considering operators \(A^\alpha\) and \(\mathscr{A}^\alpha\) instead of \(A\) and \(\mathscr{A}\) in the anisotropic case is that Theorem~\ref{thm:cornerstone2D} does not hold in the anisotropic case, due to the different form of the Neumann-to-Dirichlet operator of the homogeneous elastic 2D half-space (Theorem~\ref{thm:anisotropicD2N}).  Instead, the following modified version holds true.

\begin{theo}
	\label{thm:cornerstone2Dalpha}
	Let \(N=2\), \(\boldsymbol\Lambda\in W^{1,\infty}\) and \(\alpha\) be as in Definition~\ref{defAlpha}.  For \(\tilde{t}\in H^{-1/2}(\Gamma_C)\cap \mathscr{M}(\Gamma_C)\), we denote by \(\mathbf{u}\) the unique minimizer on \(V\) of:
	\begin{equation}
		\label{eq:functAlpha}
		\mathbf{v} \mapsto \frac{1}{2}\int_\Omega \boldsymbol\varepsilon(\mathbf{v}):\boldsymbol\Lambda\boldsymbol\varepsilon(\mathbf{v})\,{\rm d}x - \int_{\Gamma_C} \tilde{t} \,(v_t-\alpha v_n)=\frac{1}{2}\int_\Omega \boldsymbol\varepsilon(\mathbf{v}):\boldsymbol\Lambda\boldsymbol\varepsilon(\mathbf{v})\,{\rm d}x - \int_{\Gamma_C} \tilde{t} \,v^2,
	\end{equation}
	and by \(u_n=u^1\in H^{1/2}(\Gamma_C)\) the normal part of the trace of \(\mathbf{u}\) on \(\Gamma_C\).  The linear mapping 
	\(\tilde{t}\mapsto u_n=u^1=:L^\alpha\tilde{t}\) is continuous for the strong topologies of \(H^{-1/2}(\Gamma_C)\) and \(H^{1/2}(\Gamma_C)\) (and also for the weak topologies, as it is linear). 
	
	Then, for any sequence \(t^k\in K\) converging weakly in \(H^{-1/2}(\Gamma_C)\) towards a limit \(t\in K\), and any sequence \(\tilde{t}^k\rightharpoonup \tilde{t}\) in \(H^{-1/2}(\Gamma_C)\cap \mathscr{M}(\Gamma_C)\) such that \(|\tilde{t}^k| \leq -t^k\), we have:
	\[
	\lim_{k \rightarrow +\infty} \bigl\langle L^\alpha\tilde{t}^k , t^k\bigr\rangle = \bigl\langle L^\alpha\tilde{t} , t\bigr\rangle.
	\]
\end{theo}

\noindent\textbf{Proof.}  The minimizer of~\eqref{eq:functAlpha} fulfills \(t_1=t_n+\alpha t_t=0\) on \(\Gamma_C\).  Then, proofs of Lemma \ref{thm:represGreen} and of Theorem~\ref{thm:cornerstone2D} are straightforwardly adapted based on this property together with the form of the Neumann-to-Dirichlet operator provided by Theorem~\ref{thm:anisotropicD2N}.  \qed

\bigskip
\begin{coro}
	Let \(N=2\), \(\boldsymbol\Lambda\in W^{1,\infty}(\Omega)\), and \(f\in W^{1,\infty}(\Gamma_C)\) fulfilling condition~\eqref{eq:condAnisotropic}.  Let \(\tau\in K\), \(t^k\rightharpoonup t\) be a weakly convergent sequence in \(K\) such that \(\mathscr{A}^\alpha(t^k,\tau)\rightharpoonup F\) converges weakly in \(H^{1/2}(\Gamma_C)\) towards some limit \(F\).  Then, \(\lim_{k \rightarrow+\infty} \langle \mathscr{A}^\alpha(t^k,\tau),t^k\rangle = \langle F,t\rangle\).  In other words, under the previous hypotheses, property (iv) of Definition~\ref{thm:defLerayLions} holds true and the operator \(A^\alpha\) is a Leray-Lions operator.
\end{coro}

\noindent\textbf{Proof.} Adapt the proof of Corollary~\ref{thm:coroPropiv} by invoking Theorem~\ref{thm:cornerstone2Dalpha} instead of Theorem~\ref{thm:cornerstone2D}. \qed

\bigskip
Bringing all together, and invoking Proposition~\ref{thm:LerayLionsImpliesPseudomonotone} in Appendix~A, we have proved:

\begin{theo}
	Let $N=2$, \(\boldsymbol\Lambda\in W^{1,\infty}(\Omega)\), and \(f\in W^{1,\infty}(\Gamma_C)\) fulfilling condition~\eqref{eq:condAnisotropic}.  Then, the operator \(A^\alpha\) is a Leray-Lions operator and is therefore pseudomonotone.  In addition, it is coercive.
\end{theo}

We can now rely on Brézis's theorem (Theorem~\ref{thm:Brezis} in Appendix~A) to solve variational inequalities based on the operator \(A\).

\begin{theo}
	\label{theoexistGen}
	Let \(N=2\), \(\boldsymbol\Lambda\in W^{1,\infty}(\Omega)\), and \(f\in W^{1,\infty}(\Gamma_C)\) fulfilling condition~\eqref{eq:condAnisotropic}, \(\mathbf{F}\in L^2(\Omega,\mathbb{R}^N)\), \(\mathbf{T}\in L^2(\Gamma_T,\mathbb{R}^N)\), \(g\in H^{1/2}(\Gamma_C)\) and \(\mathbf{w}_t\in H^{1/2}(\Gamma_C;\mathbb{R}^N)\).  Then, there exists \(t\in K\) satisfying the variational inequality:
	\[
	\forall \hat{t}\in K, \qquad \bigl\langle A^\alpha t-g,\hat{t}-t\bigr\rangle \geq 0.
	\]
\end{theo}

It is readily checked that the unique minimizer \(\mathbf{u}\in V\) of \(E_{t,t}\) on \(V\) (see Definition~\ref{defu}) with $t$ given by Theorem \ref{theoexistGen}, solves the formal time-incremental Signorini-Coulomb problem~\eqref{eq:SignoriniCoulombDiscret} (meaning that, whenever the minimizer is smooth enough, the conditions are fulfilled pointwisely).

\subsection{The Neumann-to-Dirichlet operator of the anisotropic elastic 2D half-space}

\label{sec:N2Danisotrop}

The aim of this section is to provide a synthetic and self-contained proof of Theorem~\ref{thm:anisotropicD2N}.  This result is formally known in Solid Mechanics.  A formal version can be found, for example, in \cite[Section~8.5]{Ting}.  For the historical background of this result, the reader should refer to that book. 

\bigskip
\noindent\textbf{Proof of Theorem~\ref{thm:anisotropicD2N}.}   

\noindent\textbf{Step 1.} \textit{Preliminaries and notation.}

The orthonormal coordinate system \((x,y)\) will be used in the 2D space.  We consider an arbitrary displacement \(\mathbf{u}\) with components \(u_x(x,y)\) and \(u_y(x,y)\).  We adopt the following notation for the three independent entries of the matrix of the symmetric gradient \(\boldsymbol\varepsilon\):
\begin{align*}
	\varepsilon_1 & :=\varepsilon_{xx} = u_{x,x},\\
	\varepsilon_2 & :=\varepsilon_{yy} = u_{y,y},\\
	\varepsilon_3 & :=\varepsilon_{xy} = \varepsilon_{yx} = (u_{x,y} + u_{y,x})/2 ,\\
\end{align*}
where indices after a comma mean a derivative, as usual.  A similar convention is adopted for the stress matrix \(\boldsymbol\sigma\):
\[
\sigma_1 := \sigma_{xx},\qquad \sigma_2 := \sigma_{yy},\qquad \sigma_3 := \sigma_{xy}=\sigma_{yx},
\]
so that an arbitrary anisotropic elastic modulus tensor \(\boldsymbol\Lambda\) can now be represented by a positive definite symmetric \(3\times 3\) matrix:
\[
\begin{pmatrix}
	\sigma_1 \\
	\sigma_2 \\
	\sqrt{2}\sigma_3 
\end{pmatrix}
=
\begin{pmatrix}
	\Lambda_{11} & \Lambda_{12} & \Lambda_{13} \\
	\Lambda_{12} & \Lambda_{22} & \Lambda_{23} \\
	\Lambda_{13} & \Lambda_{23} & \Lambda_{33} 
\end{pmatrix}
\begin{pmatrix}
	\varepsilon_1 \\
	\varepsilon_2 \\
	\sqrt{2}\varepsilon_3 
\end{pmatrix}
.
\]
With this notation, the elastic equilibrium equation $\text{div}\,\boldsymbol\sigma+\mathbf{F}=\mathbf{0}$ in the whole space with arbitrary body forces \(\mathbf{F}=(F_x,F_y,0)\) read as:
\begin{equation}
	\label{eq:eqElEq}
	\begin{aligned}
		\frac{\Lambda_{33}}{2}u_{x,yy}+\frac{\sqrt{2}\Lambda_{23}}{2}u_{y,yy}+\sqrt{2}\Lambda_{13}u_{x,xy} + \Bigl(\frac{\Lambda_{33}}{2}+\Lambda_{12}\Bigr)u_{y,xy}+\Lambda_{11}u_{x,xx} + \frac{\sqrt{2}\Lambda_{13}}{2}u_{y,xx} = - F_x,\\
		\frac{\sqrt{2}\Lambda_{23}}{2}u_{x,yy} + \Lambda_{22}u_{y,yy} + \Bigl(\frac{\Lambda_{33}}{2}+\Lambda_{12}\Bigr)u_{x,xy} + \sqrt{2}\Lambda_{23}u_{y,xy} + \frac{\sqrt{2}\Lambda_{13}}{2}u_{x,xx} + \frac{\Lambda_{33}}{2}u_{y,xx} = - F_y.
	\end{aligned}
\end{equation}
Adopting convention~\eqref{eq:convFourier} for defining the Fourier transform, and denoting by \(\hathat{u}_x(\xi,\eta)\) the Fourier transform of \(u_x\) with respect to the pair \((x,y)\), we obtain the algebraic system:
\begin{equation}
	\label{eq:systDoubleFourier}
	\begin{aligned}
		\biggl[\Lambda_{11}\xi^2+\sqrt{2}\Lambda_{13}\xi\eta+\frac{\Lambda_{33}}{2}\eta^2\biggr]\hathat{u}_x + \biggl[\frac{\sqrt{2}\Lambda_{13}}{2}\xi^2+\Bigl(\frac{\Lambda_{33}}{2}+\Lambda_{12}\Bigr)\xi\eta+\frac{\sqrt{2}\Lambda_{23}}{2}\eta^2\biggr]\hathat{u}_y = \hathat{F}_x,\\
		\biggl[\frac{\sqrt{2}\Lambda_{13}}{2}\xi^2+\Bigl(\frac{\Lambda_{33}}{2}+\Lambda_{12}\Bigr)\xi\eta+\frac{\sqrt{2}\Lambda_{23}}{2}\eta^2\biggr]\hathat{u}_x + \biggl[\frac{\Lambda_{33}}{2}\xi^2+\sqrt{2}\Lambda_{23}\xi\eta+\Lambda_{22}\eta^2\biggr]\hathat{u}_y = \hathat{F}_y.
	\end{aligned}
\end{equation}
Thanks to the positive definiteness of \(\mathbf{\Lambda}\), the symmetric matrix \(\mathbf{M}(\xi,\eta)\) appearing in this system satisfies a strong ellipticity property:
\begin{multline}
	\label{eq:defDelta}
	\exists\alpha>0,\; \forall (\xi,\eta),\mathbf{X}\in\mathbb{R}^2,\quad {}^t\mathbf{X}\,\mathbf{M}(\xi,\eta)\,\mathbf{X}\geq \alpha (\xi^2+\eta^2)|\mathbf{X}|^2 \\
	\quad\implies\quad \forall(\xi,\eta)\neq 0,\quad \Delta(\xi,\eta):=\text{det}\,\mathbf{M}(\xi,\eta) > 0.
\end{multline}

\medskip
\noindent\textbf{Step 2.} \textit{General form of the solutions of the elastic equilibrium equations in a half-space.}

From now on, we look for \(u_x\) and \(u_y\) that are tempered distributions in the half-space \(\{y<0\}\).  Denoting by \(\hat{u}_x(\xi,y)\) the Fourier transforms of \(u_x\) with respect to the variable \(x\) only, the system of ordinary differential equations for \(\hat{u}_x\) and \(\hat{u}_y\) reads as:
\begin{equation}
	\label{eq:systEDO}
	\begin{aligned}
		\frac{\Lambda_{33}}{2}\hat{u}_{x,yy}+\frac{\sqrt{2}\Lambda_{23}}{2}\hat{u}_{y,yy}+i\xi\sqrt{2}\Lambda_{13}\hat{u}_{x,y} + i\xi\Bigl(\frac{\Lambda_{33}}{2}+\Lambda_{12}\Bigr)\hat{u}_{y,y} - \xi^2\Lambda_{11}\hat{u}_{x} - \xi^2 \frac{\sqrt{2}\Lambda_{13}}{2}\hat{u}_{y} = 0,\\
		\frac{\sqrt{2}\Lambda_{23}}{2}\hat{u}_{x,yy} + \Lambda_{22}\hat{u}_{y,yy} + i\xi\Bigl(\frac{\Lambda_{33}}{2}+\Lambda_{12}\Bigr)\hat{u}_{x,y} + i\xi\sqrt{2}\Lambda_{23}\hat{u}_{y,y} - \xi^2 \frac{\sqrt{2}\Lambda_{13}}{2}\hat{u}_{x} - \xi^2 \frac{\Lambda_{33}}{2}\hat{u}_{y} = 0.
	\end{aligned}
\end{equation}
We analyze first the solutions in $(\mathbb{R}\setminus\{0\})\times\mathbb{R}^-$.  This system admits nontrivial solutions of the form \(\hat{\mathbf{u}}(\xi,y)=\mathbf{v}(\xi)\,e^{i\lambda \xi y}\), if and only if \(\lambda\) is a solution of the characteristic equation:
\begin{multline}
	\label{eq:charEq}
	P(\lambda) := \frac{\Lambda_{22}\Lambda_{33}-\Lambda_{23}^2}{2}\lambda^4 + \sqrt{2}\bigl(\Lambda_{13}\Lambda_{22}-\Lambda_{23}\Lambda_{12}\bigr)\lambda^3 + \bigl(\Lambda_{11}\Lambda_{22}+\Lambda_{13}\Lambda_{23}-\Lambda_{33}\Lambda_{12}-\Lambda_{12}^2\bigr)\lambda^2 \\
	+\sqrt{2}\bigl(\Lambda_{11}\Lambda_{23}-\Lambda_{13}\Lambda_{12}\bigr)\lambda + \frac{\Lambda_{11}\Lambda_{33}-\Lambda_{13}^2}{2} = 0.
\end{multline}
We have \(\Delta(\xi,\eta) = \xi^4\,P(\eta/\xi)\), for \(\xi\neq 0\), where \(\Delta(\xi,\eta)\) is the determinant defined in~\eqref{eq:defDelta}.  The polynomial \(P\) has therefore no real root.  As a consequence, the solutions of the characteristic equation~\eqref{eq:charEq} are two pairs of conjugate complex numbers, these two pairs being either distinct (case 1) or identical (case 2).
\begin{itemize}
	[leftmargin=1.7cm]
	\item[\textbf{Case 1.}] The characteristic equation~\eqref{eq:charEq} has two distinct pairs of complex conjugate roots, say \(\alpha_1\pm i\beta_1\) and \(\alpha_2\pm i\beta_2\), where we can assume \(\beta_1>0\) and \(\beta_2>0\) without loss of generality.  On the one hand, all the solutions in \((\mathbb{R}\setminus\{0\})\times \mathbb{R}^-\) of the system~\eqref{eq:systEDO} are of the form:
	\begin{align*}
		\hat{u}_x(\xi,y) & = \sum_{k=1}^2A_k(\xi)\,f_k(\xi)\,e^{\beta_k |\xi|y}e^{i\alpha_k\xi y} +\sum_{k=1}^2\tilde A_k(\xi)\,f_k(\xi)\,e^{-\beta_k |\xi|y}e^{i\alpha_k\xi y}, \\
		\hat{u}_y(\xi,y) & = -\sum_{k=1}^2\,A_{k}(\xi)\,g_k(\xi)\,e^{\beta_k |\xi|y}e^{i\alpha_k\xi y}-\sum_{k=1}^2\,\tilde{A}_{k}(\xi)\,g_k(\xi)\,e^{-\beta_k |\xi|y}e^{i\alpha_k\xi y},
	\end{align*}
	with
	\begin{align*}
		& f_k(\xi):=\frac{\sqrt{2}\Lambda_{23}}{2}\Bigl(\beta_k^2-\alpha_k^2+2i\beta_k\alpha_k\text{\rm sgn}(\xi)\Bigr) - \Bigl(\frac{\Lambda_{33}}{2}+\Lambda_{12}\Bigr)\Bigl(\alpha_k-i\beta_k\text{\rm sgn}(\xi)\Bigr) - \frac{\sqrt{2}\Lambda_{13}}{2},\\
		& g_k(\xi):=\frac{\Lambda_{33}}{2}\Bigl(\beta_k^2-\alpha_k^2+2i\beta_k\alpha_k\text{\rm sgn}(\xi)\Bigr)- \sqrt{2}\Lambda_{13}\Bigl(\alpha_k-i\beta_k\text{\rm sgn}(\xi)\Bigr) - \Lambda_{11},
	\end{align*}
	for four arbitrary coefficients \(A_k(\xi)\), \(\tilde{A}_{k}(\xi)\).  The fact that \(\hat{u}_x\) and \(\hat{u}_y\) are tempered distributions on \(\{y<0\}\) requires that \(\tilde{A}_1\) and \(\tilde{A}_2\) should vanish identically.  
	
	On the other hand, using the fact that any distribution on \(\mathbb{R}^2\) which is supported in $\{0\}\times\mathbb{R}$ is of the form $\sum_{i=1}^n \delta^{(i)}(\xi)\,T_i(y)$ (where \(T_i\) are distributions in one variable only), we find that all the solutions of \eqref{eq:systEDO} supported in $\{0\}\times\mathbb{R}^-$ are of the form:
	\begin{align*}
		\hat{u}_x(\xi,y) & = (M_1+M_2y)\delta(\xi) + M_3\delta'(\xi),\\
		\hat{u}_y(\xi,y) & = (N_1+N_2y)\delta(\xi) + N_3\delta'(\xi),
	\end{align*}
	where the \(M_i\), \(N_i\) are complex constants and \(\delta\) is the Dirac measure.  Finally, all the tempered distributions \(\hat{u}_x\) and \(\hat{u}_y\) in \(\{y<0\}\) that solve system~\eqref{eq:systEDO} are of the form:
	\begin{align*}
		\hat{u}_x(\xi,y) & = \sum_{k=1}^2\Bigl(L_{xk1}+i\,\text{\rm sgn}(\xi)\,L_{xk2}\Bigr)A_{k}(\xi)\,e^{\beta_k |\xi|y}e^{i\alpha_k\xi y}  + (M_1+M_2y)\delta(\xi) + M_3\delta'(\xi),\\
		\hat{u}_y(\xi,y) & = \sum_{k=1}^2\Bigl(L_{yk1}+i\,\text{\rm sgn}(\xi)\,L_{yk2}\Bigr)A_{k}(\xi)\,e^{\beta_k |\xi|y}e^{i\alpha_k\xi y} + (N_1+N_2y)\delta(\xi) + N_3\delta'(\xi),
	\end{align*}
	for two arbitrary coefficients \(A_1(\xi)\), \(A_2(\xi)\) defined in $\mathbb{R}\setminus\{0\}$ and six arbitrary constants \(M_i\), \(N_i\).  Above, \(L_{xkl}\) and \(L_{ykl}\) denote eight real constants that are uniquely determined by \(\mathbf{\Lambda}\).  We now make explicit the dependence of $\hat u_x$ and $\hat u_y$ on the Dirichlet datum.  To this aim, setting:
	\begin{equation}
		\label{eq:hatU}
		\begin{aligned}
			\hat{U}_x(\xi,y) & := \sum_{k=1}^2\Bigl(L_{xk1}+i\,\text{\rm sgn}(\xi)\,L_{xk2}\Bigr)A_{k}(\xi)\,e^{\beta_k |\xi|y}e^{i\alpha_k\xi y} ,\\
			\hat{U}_y(\xi,y) & := \sum_{k=1}^2\Bigl(L_{yk1}+i\,\text{\rm sgn}(\xi)\,L_{yk2}\Bigr)A_{k}(\xi)\,e^{\beta_k |\xi|y}e^{i\alpha_k\xi y} ,
		\end{aligned}
	\end{equation}
	we observe that the mapping \((A_1(\xi),A_2(\xi))\mapsto (\hat{U}_x(\xi,0),\hat{U}_x(\xi,0))\) is linear, and the corresponding matrix has all entries of the form \(L_{ij1}+iL_{ij2}\text{\rm sgn}(\xi)\).  Such matrix is invertible and its inverse matrix has all entries of the same structure.  Indeed, assume for simplicity \((A_1,A_2)\in C^\infty_c(\mathbb{R}^+)\).  Then, one defines the inverse Fourier transform  \(\mathbf{U}(x,y)\) of $\hat{\mathbf{U}}(\xi,y):=(\hat{U}_x(\xi,y),\hat{U}_y(\xi,y))$ with respect to $\xi$, and by Plancherel theorem and the energy equality one finds that:
	\begin{equation}
		\label{eq:energyequality}
		\int_{\{y=0\}}\hat{\mathbf{U}}(\xi,0)\cdot \overline{\hat{\boldsymbol{\sigma}}(\mathbf{U})\mathbf{n}}\,{\rm d}\xi=\int_{\{y=0\}}\mathbf{U}(x,0)\cdot \boldsymbol{\sigma}(\mathbf{U})\mathbf{n}\,{\rm d}x=\int_{\{y<0\}}\boldsymbol{\varepsilon}(\mathbf{U}):\boldsymbol{\Lambda}\boldsymbol{\varepsilon}(\mathbf{U})\,{\rm d}x.
	\end{equation}
	Therefore, \(\hat{\mathbf{U}}(\xi,0)\equiv \mathbf{0}\) implies that $\mathbf{U}$ is a rigid motion and, taking into account formulae~\eqref{eq:hatU}, that $\mathbf{U}\equiv \mathbf{0}$.  In conclusion, $A_k(\xi)\equiv 0$ and  the mapping \((A_1(\xi),A_2(\xi))\mapsto (\hat{U}_x(\xi,0),\hat{U}_y(\xi,0))\) is injective.  Finally, all the tempered distributions \(\hat{u}_x\) and \(\hat{u}_y\) in \(\{y<0\}\) that solve system~\eqref{eq:systEDO} are of the form:
	\begin{equation}
		\label{eq:casOne}
		\begin{aligned}
			\hat{u}_x(\xi,y) & = \hspace*{-0.6em}\sum_{\substack{k\in\{1,2\}\\
			l\in\{x,y\}}}\hspace*{-0.4em}\Bigl(L_{xkl1}+iL_{xkl2}\,\text{\rm sgn}(\xi)\Bigr)\,U_l(\xi)\,e^{\beta_k |\xi|y}e^{i\alpha_k\xi y}  + (M_1+M_2y)\delta(\xi) + M_3\delta'(\xi),\\
			\hat{u}_y(\xi,y) & = \hspace*{-0.6em}\sum_{\substack{k\in\{1,2\}\\
			l\in\{x,y\}}}\hspace*{-0.4em}\Bigl(L_{ykl1}+iL_{ykl2}\,\text{\rm sgn}(\xi)\Bigr)\,U_l(\xi)\,e^{\beta_k |\xi|y}e^{i\alpha_k\xi y} + (N_1+N_2y)\delta(\xi) + N_3\delta'(\xi),
		\end{aligned}
	\end{equation}
	for two arbitrary coefficients \(U_x(\xi)\), \(U_y(\xi)\) defined in $\mathbb{R}\setminus\{0\}$ and six arbitrary constants \(M_i\), \(N_i\).  Above, \(L_{ijkl}\) denote sixteen real constants that are uniquely determined by \(\mathbf{\Lambda}\) and that ensure that \(\hat{U}_i(\xi,0)=U_i(\xi)\).
	
	\item[\textbf{Case 2.}] The characteristic equation~\eqref{eq:charEq} has a pair of double complex conjugate roots \(\alpha\pm i\beta\), where \(\beta>0\).  This degenerate situation is the one encountered in the case of isotropic elasticity.  In that case, one can conclude by a similar argument as in case 1, that all the tempered distributions \(\hat{u}_x\) and \(\hat{u}_y\) in \(\{y<0\}\) that solve system~\eqref{eq:systEDO} are of the form:
	\begin{equation}
		\label{eq:caseTwo}
		\begin{aligned}
			\hat{u}_x(\xi,y) & = \Biggl(U_x(\xi) + y|\xi|\biggl[\sum_{k\in\{x,y\}} \Bigl(L_{xk1}+iL_{xk2}\,\text{\rm sgn}(\xi)\Bigr)U_k(\xi)\biggr]\Biggr)e^{\beta |\xi|y}e^{i\alpha\xi y} \\
			&\hspace*{6cm}\mbox{} + (M_1+M_2y)\delta(\xi) + M_3\delta'(\xi),\\
			\hat{u}_y(\xi,y) & = \Biggl(U_y(\xi) + y|\xi|\biggl[\sum_{k\in\{x,y\}} \Bigl(L_{yk1}+iL_{yk2}\,\text{\rm sgn}(\xi)\Bigr)U_k(\xi)\biggr]\Biggr)e^{\beta |\xi|y}e^{i\alpha\xi y} \\
			&\hspace*{6cm}\mbox{}+ (N_1+N_2y)\delta(\xi) + N_3\delta'(\xi),
		\end{aligned}
	\end{equation}
	for two arbitrary coefficients \(U_x(\xi)\), \(U_y(\xi)\) defined in $\mathbb{R}\setminus\{0\}$ and six constants \(M_i\), \(N_i\).  Above, \(L_{xkl}\) and \(L_{ykl}\) denote eight real constants that are uniquely determined by \(\mathbf{\Lambda}\). 
\end{itemize}

\medskip
\noindent\textbf{Step 3.} \textit{General form of the Neumann-to-Dirichlet operator of the half-space.}

One compute the Fourier transform \(\hat{t}_x\) and \(\hat{t}_y\) of the components of the surface force by use of the formulae:
\begin{align*}
	\hat{t}_x(\xi) & = \hat{\sigma}_{xy}(\xi) = i\xi\frac{\sqrt{2}\Lambda_{13}}{2}\hat{u}_x(\xi,0) + \frac{\sqrt{2}\Lambda_{23}}{2}\frac{\partial \hat{u}_y}{\partial y}(\xi,0)+ i\xi\frac{\Lambda_{33}}{2}\hat{u}_y(\xi,0) + \frac{\Lambda_{33}}{2}\frac{\partial \hat{u}_x}{\partial y}(\xi,0),\\
	\hat{t}_y(\xi) & = \hat{\sigma}_{yy}(\xi) = i\xi\Lambda_{12}\hat{u}_x(\xi,0) 
	+ \Lambda_{22}\frac{\partial \hat{u}_y}{\partial y}(\xi,0) + i\xi\frac{\Lambda_{23}}{2}\hat{u}_y(\xi,0) + \frac{\Lambda_{23}}{2}\frac{\partial \hat{u}_x}{\partial y}(\xi,0).
\end{align*}
Applying the above formulae to either~\eqref{eq:casOne} or~\eqref{eq:caseTwo}, we get:
\begin{equation}
	\label{eq:eqCaseOne}
	\begin{aligned}
		\hat{t}_x(\xi) & = |\xi| \Bigl(\bigl[L_{xx}^++iL_{xx}^-\text{\rm sgn}(\xi)\bigr]U_x(\xi)+\bigl[L_{xy}^++iL_{xy}^-\text{\rm sgn}(\xi)\bigr]U_y(\xi)\Bigr)+M\delta(\xi),\\
		\hat{t}_y(\xi) & = |\xi| \Bigl(\bigl[L_{yx}^++iL_{yx}^-\text{\rm sgn}(\xi)\bigr]U_x(\xi)+\bigl[L_{yy}^++iL_{yy}^-\text{\rm sgn}(\xi)\bigr]U_y(\xi)\Bigr)+N\delta(\xi),
	\end{aligned}
\end{equation}
where the \(L_{ij}^\pm\) are eight real constants that are uniquely determined by \(\mathbf{\Lambda}\), and \(M\) and \(N\) are complex constants that depend on \(\mathbf{\Lambda}\), but also on the \(M_i\)'s and the \(N_i\)'s.  Considering \(M=N=0\) and \(U_x,U_y\in C^\infty_c(\mathbb{R})\), the energy equality~\eqref{eq:energyequality} shows that the bilinear mapping:
\[
\bigl(\mathbf{U}^1,\mathbf{U^2}\bigr) \mapsto   \int_{-\infty}^{+\infty}
\bigl(U_x^1\;U_y^1\bigr)
\begin{pmatrix}
	L_{xx}^++iL_{xx}^-\text{\rm sgn}(\xi)  &  L_{xy}^++iL_{xy}^-\text{\rm sgn}(\xi) \\
	L_{yx}^++iL_{yx}^-\text{\rm sgn}(\xi)  &  L_{yy}^++iL_{yy}^-\text{\rm sgn}(\xi)
\end{pmatrix}
\begin{pmatrix}
	\overline{U}_x^2 \\
	\overline{U}_y^2
\end{pmatrix}
\,|\xi|\,{\rm d}\xi
\]
should be positive definite, with Hermitian symmetry.  Therefore, the matrix:
\[
\begin{pmatrix}
	L_{xx}^++iL_{xx}^-\text{\rm sgn}(\xi)  &  L_{xy}^++iL_{xy}^-\text{\rm sgn}(\xi) \\
	L_{yx}^++iL_{yx}^-\text{\rm sgn}(\xi)  &  L_{yy}^++iL_{yy}^-\text{\rm sgn}(\xi)
\end{pmatrix}
\]
should be positive definite, Hermitian, for all \(\xi\neq 0\).  Hence, \(L_{xx}^-=L_{yy}^-=0\), \(L_{yx}^+=L_{xy}^+\), \(L_{yx}^-=-L_{xy}^-\), \(L_{xx}^+>0\), \(L_{yy}^+>0\) and its inverse matrix has the form:
\begin{equation}
	\label{eq:matrix}
	\begin{pmatrix}
		\pi\mathscr{C}_4 & \pi\mathscr{C}_2 - 2i \mathscr{C}_3\,\text{\rm sgn}(\xi) \\
		\pi\mathscr{C}_2 + 2i \mathscr{C}_3\,\text{\rm sgn}(\xi) & \pi\mathscr{C}_1
	\end{pmatrix}
\end{equation}
where the \(\mathscr{C}_i\)'s are uniquely determined by \(\mathbf{\Lambda}\), \(\mathscr{C}_1>0\) and \(\mathscr{C}_4>0\).  Finally, taking \(\hat{t}_x(\xi)=T_x/\sqrt{2\pi}\), \(\hat{t}_y(\xi)=T_y/\sqrt{2\pi}\) for \((T_x,T_y)\in\mathbb{R}^2\) in \eqref{eq:eqCaseOne}, all the tempered distributions \(U_x\), \(U_y\) on \(\mathbb{R}\setminus\{0\}\) satisfying \eqref{eq:eqCaseOne} where necessarily \(M=N=0\), are given by:
\begin{align*}
	U_x(\xi) & = \mathscr{C}_4\,T_x\,\sqrt{\frac{\pi}{2}}\,\frac{1}{|\xi|} +  \mathscr{C}_2\,T_y\,\sqrt{\frac{\pi}{2}}\,\frac{1}{|\xi|} - i \mathscr{C}_3\,T_y\,\sqrt{\frac{2}{\pi}}\,\frac{1}{\xi},\\
	U_y(\xi) & =  \mathscr{C}_1\,T_y\,\sqrt{\frac{\pi}{2}}\,\frac{1}{|\xi|} +  \mathscr{C}_2\,T_x\,\sqrt{\frac{\pi}{2}}\,\frac{1}{|\xi|} + i \mathscr{C}_3\,T_y\,\sqrt{\frac{2}{\pi}}\,\frac{1}{\xi}.
\end{align*}
Using~\eqref{eq:FourierTransform}, it turns out that all the tempered distributions \(u_x\), \(u_y\) in the half-space, that satisfy the elastic equilibrium equations with vanishing body forces and surface forces equal to \((T_x,T_y)\delta\) on the boundary have a trace on the boundary given by:
\begin{align*}
	u_x & = - \mathscr{C}_4\,T_x\,\log|x| - \mathscr{C}_2\,T_y\,\log|x| + \mathscr{C}_3\,T_y\,\text{\rm sgn}(x) + a_x(x),\\
	u_y & = - \mathscr{C}_1\,T_y\,\log|x| - \mathscr{C}_2\,T_x\,\log|x| - \mathscr{C}_3\,T_x\,\text{\rm sgn}(x) + a_y(x),
\end{align*}
where \(\mathbf{a}\) denote an arbitrary affine displacement in the half-space, which is compatible with boundary free of surface force (such an arbitrary displacement is determined by four arbitrary real constants: three components of an overall rigid motion and one component of a constant \(\sigma_{xx}\) stress).  The contribution of \(\Gamma_{\rm Eul}\) in formulae~\eqref{eq:FourierTransform} has been included in \(a_x\) and \(a_y\).  This result is enough to obtain the form of the Neumann-to-Dirichlet operator of the anisotropic homogeneous half-space as given in Theorem~\ref{thm:anisotropicD2N}.

\medskip
\noindent\textbf{Step 4.} \textit{The maps \(\mathbf{\Lambda}\mapsto \mathscr{C}_i\) are of class \(C^\infty\).}

The preceding Steps~1, 2, 3 contain actually an effective method of calculating the \(\mathscr{C}_i\)'s in terms of the entries of \(\mathbf{\Lambda}\) and the complex roots of the characteristic equation~\eqref{eq:charEq}.  The expression can even be made explicit if desired.  All the \(\mathscr{C}_i\)'s are rational functions of the entries of \(\mathbf{\Lambda}\) and the complex roots of the characteristic equation~\eqref{eq:charEq}.  This does not readily give the expected conclusion, as the roots of a polynomial are \(C^\infty\) functions of the coefficients of the polynomial only in the case where the roots are simple and the root functions do not even need to be Lipschitz-continuous at a multiple root.

We are therefore going to provide an alternative expression of the \(\mathscr{C}_i\)'s in terms of \(\mathbf{\Lambda}\), better suited to prove the regularity of the dependence.  The method runs as follows.  We are going to look for a tempered distribution \(\mathbf{u}\) \emph{on the whole space} \(\mathbb{R}^2\), satisfying the elastic equilibrium with body forces of the form \(\mathbf{F}(x)\delta(y)\) and adjust \(\mathbf{F}(x)\) so that the restriction of \(\mathbf{u}\) to the half-space \(\{y<0\}\) is a fundamental solution for the Neumann problem on the half-space.  We will prove that all the fundamental solutions can be obtained by means of an appropriate choice of the distribution \(\mathbf{F}(x)\) and, thus, get an alternative expression of the \(\mathscr{C}_i\)'s in terms of \(\mathbf{\Lambda}\).

The inversion of the algebraic system~\eqref{eq:systDoubleFourier} fulfilled by the double Fourier transform \(\hathat{\mathbf{u}}\) gives:
\begin{equation}
	\label{eq:DoubleFourierU}
	\begin{aligned}
		\hathat{u}_x(\xi,\eta) & = \frac{1}{\sqrt{2\pi}\Delta(\xi,\eta)}\Biggl\{\biggl[\frac{\Lambda_{33}}{2}\xi^2+\sqrt{2}\Lambda_{23}\xi\eta+\Lambda_{22}\eta^2\biggr]\hat{F}_x(\xi) \\
		& \hspace*{4cm}\mbox{}-\biggl[\frac{\sqrt{2}\Lambda_{13}}{2}\xi^2+\Bigl(\frac{\Lambda_{33}}{2}+\Lambda_{12}\Bigr)\xi\eta+\frac{\sqrt{2}\Lambda_{23}}{2}\eta^2\biggr]\hat{F}_y(\xi)\Biggr\},\\
		\hathat{u}_y(\xi,\eta) & = \frac{1}{\sqrt{2\pi}\Delta(\xi,\eta)}\Biggl\{-\biggl[\frac{\sqrt{2}\Lambda_{13}}{2}\xi^2+\Bigl(\frac{\Lambda_{33}}{2}+\Lambda_{12}\Bigr)\xi\eta+\frac{\sqrt{2}\Lambda_{23}}{2}\eta^2\biggr]\hat{F}_x(\xi) \\
		& \hspace*{4cm}\mbox{}+\biggl[\Lambda_{11}\xi^2+\sqrt{2}\Lambda_{13}\xi\eta+\frac{\Lambda_{33}}{2}\eta^2\biggr]\hat{F}_y(\xi)\Biggr\},
	\end{aligned}
\end{equation}
where the determinant \(\Delta(\xi,\eta)\) was seen to be expressed in terms of the characteristic polynomial~\eqref{eq:charEq} by the formula \(\Delta(\xi,\eta) = \xi^4P(\eta/\xi)\), so that \(\Delta(\xi,\eta)>0\), for \(\xi\neq 0\).  Actually, the inversion of the algebraic system~\eqref{eq:systDoubleFourier} provides \(\hathat{\mathbf{u}}\) up to an arbitrary additive term of the form \(\mathbf{C}_0\delta(\xi)\delta(\eta)+\mathbf{C}_1\delta'(\xi)\delta(\eta)+\mathbf{C}_2\delta(\xi)\delta'(\eta)\), with \(\mathbf{C}_i\in \mathbb{R}^2\).  This term represents an arbitrary affine displacement which plays no role in the sequel and will therefore be ignored.  As the right-hand sides are integrable functions of \(\eta\), one can use the formula:
\begin{equation}
	\label{eq:zeroValueSimpleFT}
	\hat{u}_x(\xi,0) = \frac{1}{\sqrt{2\pi}}\int_{-\infty}^{+\infty} \hathat{u}_x(\xi,\eta)\,{\rm d}\eta,
\end{equation}
to get:
\begin{equation}
	\label{eq:surfDepl}
	\begin{aligned}
		\hat{u}_x(\xi,0) & = \frac{1}{2\pi|\xi|}\Biggl\{\biggl[\frac{\Lambda_{33}}{2}I_0+\sqrt{2}\Lambda_{23}I_1+\Lambda_{22}I_2\biggr]\hat{F}_x(\xi) \\
		& \hspace*{4cm}\mbox{}-\biggl[\frac{\sqrt{2}\Lambda_{13}}{2}I_0+\Bigl(\frac{\Lambda_{33}}{2}+\Lambda_{12}\Bigr)I_1+\frac{\sqrt{2}\Lambda_{23}}{2}I_2\biggr]\hat{F}_y(\xi)\Biggr\},\\
		\hat{u}_y(\xi,0) & = \frac{1}{2\pi|\xi|}\Biggl\{-\biggl[\frac{\sqrt{2}\Lambda_{13}}{2}I_0+\Bigl(\frac{\Lambda_{33}}{2}+\Lambda_{12}\Bigr)I_1+\frac{\sqrt{2}\Lambda_{23}}{2}I_2\biggr]\hat{F}_x(\xi) \\
		& \hspace*{4cm}\mbox{}+\biggl[\Lambda_{11}I_0+\sqrt{2}\Lambda_{13}I_1+\frac{\Lambda_{33}}{2}I_2\biggr]\hat{F}_y(\xi)\Biggr\},
	\end{aligned}
\end{equation}
where:
\[
I_0:=\int_{-\infty}^{+\infty} \frac{{\rm d}x}{P(x)}>0,\qquad I_1:=\int_{-\infty}^{+\infty} \frac{x\,{\rm d}x}{P(x)},\qquad I_2:=\int_{-\infty}^{+\infty} \frac{x^2\,{\rm d}x}{P(x)}>0,
\]
depend only on \(\mathbf{\Lambda}\) and are \(C^\infty\) functions of \(\mathbf{\Lambda}\), by dominated convergence.  From~\eqref{eq:DoubleFourierU}, we get:
\begin{align*}
	\hathat{\sigma}_{xy}(\xi,\eta) & = \frac{i}{\Delta(\xi,\eta)}\Biggl\{\biggl[ \frac{\Lambda_{22}\Lambda_{33}\!-\!\Lambda_{23}^2}{2}\eta^3+\frac{\sqrt{2}}{2}\bigl(\Lambda_{13}\Lambda_{22}\!-\!\Lambda_{23}\Lambda_{12}\bigr)\eta^2\xi+\frac{\Lambda_{13}\Lambda_{23}\!-\!\Lambda_{33}\Lambda_{12}}{2}\eta\xi^2\biggr]\frac{\hat{F}_x(\xi)}{\sqrt{2\pi}}\\
	& \hspace*{1.3cm}\mbox{}+\biggl[\frac{\Lambda_{13}\Lambda_{23}\!-\!\Lambda_{33}\Lambda_{12}}{2}\eta^2\xi+\frac{\sqrt{2}}{2}\bigl(\Lambda_{23}\Lambda_{11}\!-\!\Lambda_{13}\Lambda_{12}\bigr)\eta\xi^2+\frac{\Lambda_{11}\Lambda_{33}\!-\!\Lambda_{13}^2}{2}\xi^3\biggr]\frac{\hat{F}_y(\xi)}{\sqrt{2\pi}}\Biggr\},\\
	\hathat{\sigma}_{yy}(\xi,\eta) & = \frac{i}{\Delta(\xi,\eta)}\Biggl\{\!\!-\!\biggl[ \frac{\Lambda_{22}\Lambda_{33}\!-\!\Lambda_{23}^2}{2}\eta^2\xi+\frac{\sqrt{2}}{2}\bigl(\Lambda_{13}\Lambda_{22}\!-\!\Lambda_{23}\Lambda_{12}\bigr)\eta\xi^2+\frac{\Lambda_{13}\Lambda_{23}\!-\!\Lambda_{33}\Lambda_{12}}{2}\xi^3\biggr]\frac{\hat{F}_x(\xi)}{\sqrt{2\pi}}\\
	& \hspace*{2.0cm}\mbox{}+\biggl[\frac{\Lambda_{22}\Lambda_{33}\!-\!\Lambda_{23}^2}{2}\eta^3+\sqrt{2}\bigl(\Lambda_{13}\Lambda_{22}\!-\!\Lambda_{23}\Lambda_{12}\bigr)\eta^2\xi\\
	& \hspace*{1.5cm}\mbox{}+\Bigl(\Lambda_{11}\Lambda_{22}-\Lambda_{12}^2+\frac{\Lambda_{13}\Lambda_{23}\!-\!\Lambda_{33}\Lambda_{12}}{2}\Bigr)\eta\xi^2+\frac{\sqrt{2}}{2}\bigl(\Lambda_{11}\Lambda_{23}\!-\!\Lambda_{12}\Lambda_{13}\bigr)\xi^3\biggr]\frac{\hat{F}_y(\xi)}{\sqrt{2\pi}}\Biggr\}.	
\end{align*}
As the function \(\eta\mapsto\eta^3/\Delta(\xi,\eta)\) is not integrable, we cannot simply rely on formula~\eqref{eq:zeroValueSimpleFT} to compute \(\hat{\sigma}_{xy}(\xi,0)\).  Instead, we write:
\[
\frac{\Lambda_{22}\Lambda_{33}\!-\!\Lambda_{23}^2}{2}\frac{\tilde{\eta}^3}{P(\tilde{\eta})} = \frac{\tilde{\eta}}{1+\tilde{\eta}^2}+ \biggl[\frac{\Lambda_{22}\Lambda_{33}\!-\!\Lambda_{23}^2}{2}\frac{\tilde{\eta}^3}{P(\tilde{\eta})}-\frac{\tilde{\eta}}{1+\tilde{\eta}^2}\biggr],
\]
where the second term in the right-hand side is an integrable rational function of \(\tilde{\eta}\).  Making use of the knowledge of the inverse Fourier transform of \(\tilde{\eta}/(1+\tilde{\eta}^2)\):
\[
\frac{1}{\sqrt{2\pi}}\int_{-\infty}^{+\infty}\frac{\tilde{\eta}\,e^{iy\tilde{\eta}}}{1+\tilde{\eta}^2}\,{\rm d}\tilde{\eta} = i\sqrt{\frac{\pi}{2}}\,e^{-|y|}\,\text{sgn}(y),
\]
we obtain:
\[
\lim_{y \rightarrow 0\pm} \frac{\Lambda_{22}\Lambda_{33}\!-\!\Lambda_{23}^2}{2}\int_{-\infty}^{+\infty} \frac{\eta^3\,e^{iy\eta}}{\Delta(\xi,\eta)}\,{\rm d}\eta = \pm i\pi + I_3\,\text{sgn}\,\xi, 
\]
where:
\[
I_3 := \int_{-\infty}^{+\infty} \biggl[\frac{\Lambda_{22}\Lambda_{33}\!-\!\Lambda_{23}^2}{2}\frac{x^3}{P(x)}-\frac{x}{1+x^2}\biggr]\,{\rm d}x,
\]
is a \(C^\infty\) function of \(\mathbf{\Lambda}\) only.  As a result:
\begin{equation}
	\label{eq:systInvertible}
	\begin{aligned}
		\hat{\sigma}_{xy}(\xi,0-) & = \Biggl\{\pi+i\,\text{sgn}\,\xi\biggl[I_3+\frac{\sqrt{2}}{2}\bigl(\Lambda_{13}\Lambda_{22}\!-\!\Lambda_{23}\Lambda_{12}\bigr)I_2+\frac{\Lambda_{13}\Lambda_{23}\!-\!\Lambda_{33}\Lambda_{12}}{2}I_1\biggr]\Biggr\}\frac{\hat{F}_x(\xi)}{2\pi}
		\\
		& \phantom{\mbox{}=\mbox{}}\mbox{} + i\,\text{sgn}\,\xi\Biggl\{\frac{\Lambda_{13}\Lambda_{23}\!-\!\Lambda_{33}\Lambda_{12}}{2}I_2+\frac{\sqrt{2}}{2}\bigl(\Lambda_{23}\Lambda_{11}\!-\!\Lambda_{13}\Lambda_{12}\bigr)I_1+\frac{\Lambda_{11}\Lambda_{33}\!-\!\Lambda_{13}^2}{2}I_0\Biggr\}\frac{\hat{F}_y(\xi)}{2\pi},\\
		\hat{\sigma}_{yy}(\xi,0-) & = -i\,\text{sgn}\,\xi\Biggl\{\frac{\Lambda_{22}\Lambda_{33}\!-\!\Lambda_{23}^2}{2}I_2+\frac{\sqrt{2}}{2}\bigl(\Lambda_{13}\Lambda_{22}\!-\!\Lambda_{23}\Lambda_{12}\bigr)I_1+\frac{\Lambda_{13}\Lambda_{23}\!-\!\Lambda_{33}\Lambda_{12}}{2}I_0\Biggr\}\frac{\hat{F}_x(\xi)}{2\pi}\\
		& \phantom{\mbox{}=\mbox{}}\mbox{} + \Biggl\{\pi+i\,\text{sgn}\,\xi\biggl[I_3+\sqrt{2}\bigl(\Lambda_{13}\Lambda_{22}\!-\!\Lambda_{23}\Lambda_{12}\bigr)I_2\\
		& \hspace*{1.5cm}\mbox{}+\Bigl(\Lambda_{11}\Lambda_{22}-\Lambda_{12}^2+\frac{\Lambda_{13}\Lambda_{23}\!-\!\Lambda_{33}\Lambda_{12}}{2}\Bigr)I_1+\frac{\sqrt{2}}{2}\bigl(\Lambda_{11}\Lambda_{23}\!-\!\Lambda_{12}\Lambda_{13}\bigr)I_0\biggr]\Biggr\}\frac{\hat{F}_y(\xi)}{2\pi}.
	\end{aligned}
\end{equation}
The next step is now to find \((\hat{F}_x,\hat{F}_y)\) such that:
\[
\hat{\sigma}_{xy}(\xi,0-) = \frac{T_x}{\sqrt{2\pi}},\qquad \hat{\sigma}_{yy}(\xi,0-) = \frac{T_y}{\sqrt{2\pi}},
\]
where \((T_x,T_y)\) is arbitrarily fixed in \(\mathbb{R}^2\).   We therefore have to prove that system~\eqref{eq:systInvertible} is invertible and this is where the preliminary Steps~1, 2, 3 will be essential.  Take \(\hat{\sigma}_{xy}(\xi,0-) = \hat{\sigma}_{yy}(\xi,0-) = 0\) in system~\eqref{eq:systInvertible}.  Then, on the one hand, the restriction of \(\mathbf{u}(x,y)\) to the half-space \(\{y<0\}\) solves an elasticity problem in that half-space with free boundary.  The analysis in Step~1, 2, 3 shows that the restriction of the corresponding \(\mathbf{u}(x,y)\) to the half-space \(\{y<0\}\) must be affine so that \(\hat{\mathbf{u}}(\xi,0) = \mathbf{C}_0 \delta(\xi)+ \mathbf{C}_1\delta'(\xi)\), where \(\mathbf{C}_0,\mathbf{C}_1\in\mathbb{R}^2\).  On the other hand, by system~\eqref{eq:systInvertible}, we obtain:
\[
\hat{\sigma}_{xy}(\xi,0+) = - \hat{F}_x(\xi),\qquad \hat{\sigma}_{yy}(\xi,0+) =  - \hat{F}_y(\xi),
\]
so that the restriction of \(\mathbf{u}(x,y)\) to the half-space \(\{y>0\}\) solves a Neumann elastic problem on that half-space with boundary data \((\hat{F}_x(\xi),\hat{F}_y(\xi))\).  As \(\hat{\mathbf{u}}(\xi,0) = \mathbf{C}_0 \delta(\xi)+ \mathbf{C}_1\delta'(\xi)\), the analysis in Steps~1, 2, 3 shows that \((\hat{F}_x,\hat{F}_y)=(0,0)\).  Hence, the $2\times 2$ matrix appearing in system~\eqref{eq:systInvertible} is invertible.  We can therefore compute
\((\hat{F}_x,\hat{F}_y)\) in terms of \((T_x,T_y)\):
\[
\begin{pmatrix}
	\hat{F}_x(\xi) \\
	\hat{F}_y(\xi) 
\end{pmatrix}
= 
\begin{pmatrix}
	C_1 + i \,C_2 \,\text{sgn}\,\xi & C_3 + i \,C_4 \,\text{sgn}\,\xi \\
	C_5 + i \,C_6 \,\text{sgn}\,\xi & C_7 + i \,C_8 \,\text{sgn}\,\xi
\end{pmatrix}
\begin{pmatrix}
	T_x \\
	T_y 
\end{pmatrix}
,
\]
where the \(C_i\)'s are rational functions of the entries of \(\mathbf{\Lambda}\) and of the \(I_i\)'s only.  Injecting this expression in~\eqref{eq:surfDepl}, we obtain:
\[
\begin{pmatrix}
	\hat{u}_x(\xi,0) \\
	\hat{u}_y(\xi,0) 
\end{pmatrix}
= \frac{1}{\sqrt{2\pi}|\xi|}
\begin{pmatrix}
	C_1' + i \,C_2' \,\text{sgn}\,\xi & C_3' + i \,C_4' \,\text{sgn}\,\xi \\
	C_5' + i \,C_6' \,\text{sgn}\,\xi & C_7' + i \,C_8' \,\text{sgn}\,\xi
\end{pmatrix}
\begin{pmatrix}
	T_x \\
	T_y 
\end{pmatrix}
,
\]
where the \(C_i'\)'s are rational functions of the entries of \(\mathbf{\Lambda}\) and of the \(I_i\)'s only.  As a consequence, the \(C_i'\)'s are \(C^\infty\) functions of \(\mathbf{\Lambda}\).  As the above matrix equals~\eqref{eq:matrix}, we have therefore proved that all the \(\mathscr{C}_i\)'s are \(C^\infty\) functions of \(\mathbf{\Lambda}\). \qed

\section*{Appendix~A:~Pseudomonotone and Leray-Lions operators}

\setcounter{section}{1}
\setcounter{defi}{0}
\renewcommand\thesection{\Alph{section}}

Let \(B\) be a reflexive Banach space, \(B^*\) its dual space, \(K\subset B\) a nonempty closed convex subset, \(f\in B^*\) and \(A:K \rightarrow B^*\).  In~\cite{BrezisIeq}, Brézis sought minimal conditions on \(A\) to ensure the solvability of the problem of finding \(u\in K\) satisfying the variational inequality:
\begin{equation}
	\label{eq:varIneq}
	\forall v\in K,\qquad \bigl\langle Au-f, v-u\bigr\rangle \geq 0.
\end{equation}
More precisely, Brézis sought the minimal conditions on \(A\) ensuring that a sequence of Galerkin approximations for the solution of~\eqref{eq:varIneq} can be built based on Brouwer's theorem and that a limit to a solution of~\eqref{eq:varIneq} can be taken.

\begin{defi}
	\label{thm:genDefi}
	A mapping \(A:K \rightarrow B^*\) is said \emph{hemicontinuous} if each real-valued function: \(\lambda \mapsto \bigl\langle A\bigl(\lambda u +(1-\lambda)v\bigr)\,,\, v-u\bigr\rangle\) is continuous on \([0,1]\).  The mapping \(A:K \rightarrow B^*\) is said \emph{bounded} if it maps bounded subsets into bounded subsets.  The mapping \(A:K \rightarrow B^*\) is said \emph{monotone} if \(\langle Au - Av,u-v\rangle\geq 0\), for all \(u,v\in K\).
\end{defi}

\begin{defi}
	\label{thm:defPseudoMonot}
	A mapping \(A:K \rightarrow B^*\) is said \emph{pseudomonotone} (in the sense of Brézis) if it is bounded and if, for all weakly converging sequence \(u_n \rightharpoonup u\) in \(K\), such that \(\limsup \langle Au_n, u_n-u\rangle \leq 0\), then we have:
	\[
	\forall v\in K,\qquad \liminf \langle Au_n, u_n-v\rangle \geq \langle Au, u-v\rangle.
	\]
	A mapping \(A:K \rightarrow B^*\) is said \emph{coercive} if:
	\[
	\lim_{\substack{\|u\|_{B}\rightarrow +\infty,\\
	u\in K}} \frac{\langle Au,u\rangle}{\|u\|_{B}} = +\infty.
	\]
\end{defi}

The motivation for the previous definition lies in the following result due to Brézis~\cite[Corollary 30]{BrezisIeq}.

\begin{theo}
	[Brézis] 
	\label{thm:Brezis}
	Let \(A:K \rightarrow B^*\) be a pseudomonotone, coercive operator.  Then, for all \(f\in B^*\), the variational inequality~\eqref{eq:varIneq} has at least one solution.
\end{theo}

Pseudomonotone operators are a generalization of monotone operators in the following sense:

\begin{prop}
	Let \(A:K \rightarrow B^*\) be a bounded, hemicontinuous, monotone operator.  Then, \(A\) is pseudomonotone.
\end{prop}

Part of the success of Brézis's result is due to the fact that the class of pseudomonotone operators encompasses the subclass of the so-called Leray-Lions operators previously introduced in~\cite{LerayLions} in view of the study of certain nonlinear elliptic partial differential equations.

\begin{defi}
	\label{thm:defLerayLions}
	A mapping \(A:K \rightarrow B^*\) is said to be a \emph{Leray-Lions} operator, if it is bounded and can be written as \(Au=\mathscr{A}(u,u)\), for all \(u\in K\), where \(\mathscr{A}:K\times K \rightarrow B^*\) has the following properties.
	\begin{itemize}
		\item[(i)] For all \(u\in K\), the mapping:
		\[
		\left\{
		\begin{array}{rcl}
			K & \rightarrow & B^*\\
			v & \mapsto & \mathscr{A}(u,v)
		\end{array}
		\right.
		\]
		is bounded, hemicontinuous and satisfies the monotonicity property:
		\[
		\forall v\in K,\qquad \bigl\langle \mathscr{A}(u,u) -\mathscr{A}(u,v)\,,\, u-v\bigr\rangle \geq 0. 
		\]
		\item[(ii)] For all \(v\in K\), the mapping \(u\mapsto \mathscr{A}(u,v)\) is bounded, hemicontinuous from \(K\subset B\) to \(B^*\).
		\item[(iii)] For all \(v\in K\), if \(u_n \rightharpoonup u\) is a weakly converging sequence in \(K\), such that \(\lim \langle\mathscr{A}(u_n,u_n)-\mathscr{A}(u_n,u),u_n-u\rangle=0\), then we have \(\mathscr{A}(u_n,v) \rightharpoonup \mathscr{A}(u,v)\) weakly in \(B^*\).
		\item[(iv)] For all \(v\in K\), if \(u_n \rightharpoonup u\) is a weakly converging sequence in \(K\), such that \(\mathscr{A}(u_n,v) \rightharpoonup F\) weakly in \(B^*\), then we have \(\lim\langle\mathscr{A}(u_n,v),u_n\rangle = \langle F,u\rangle\).
	\end{itemize}
\end{defi}

The proof of the following result is easy (see, for example,~\cite[Lemma 4.13]{Troianiello})).

\begin{prop}
	\label{thm:LerayLionsImpliesPseudomonotone}
	Every Leray-Lions operator \(A:K \rightarrow B^*\) is pseudomonotone.
\end{prop}

\section*{Acknowledgements}

We warmly thank François Murat for a very enlightening discussion and, in particular, for pointing out his article \cite{Murat}.




\begin{thebibliography}{9}

	\bibitem{Andersson} \textsc{L.E. Andersson} (2000), Existence results for quasistatic contact problems with Coulomb friction. \textit{Applied Mathematics \& Optimization}, \textbf{42}, pp~169--202.

	\bibitem{AnderssonKlarbring} \textsc{L.E. Andersson} and \textsc{A. Klarbring} (2001), A review of the theory of static and quasi-static frictional contact problems in elasticity. \textit{Philosophical Transactions of the Royal Society of London. Series A: Mathematical, Physical and Engineering Sciences}, \textbf{359}, pp~2519--2539.

	\bibitem{BallardIurlano} \textsc{P. Ballard} and \textsc{F. Iurlano} (2022), Homogenization of Friction in a 2D Linearly Elastic Contact Problem. \textit{Journal of Elasticity}, \textbf{150}, pp~261–-325. 
	
	\bibitem{bibiJiri} \textsc{P. Ballard}, \textsc{J. Jarušek} (2011), Indentation of an elastic half-space by a rigid flat punch as a model problem for analysing contact problems with Coulomb friction. \textit{Journal of Elasticity} \textbf{103}, pp~15--52.
	
	\bibitem{BrezisIeq}	H.~\textsc{Brézis} (1968), Équations et inéquations non linéaires dans les espaces vectoriels en dualité. \textit{Annales de l'Institut Fourier} \textbf{18},~(1), pp~115--175.
	
	\bibitem{DalMasoPlast} G.~\textsc{Dal Maso}, A.~\textsc{DeSimone} and M.G.~\textsc{Mora} (2006), Quasistatic evolution problems for linearly elastic–perfectly plastic materials. \textit{Archive for Rational Mechanics and Analysis} \textbf{180}, pp~237--291.
	
	\bibitem{DuvautLions} G.~\textsc{Duvaut} and J.L.~\textsc{Lions} (1972), \textit{Les inéquations en mécanique et en physique}.  Dunod, Paris.
	
	\bibitem{EckJarusek} \textsc{C. Eck} and \textsc{J. Jarušek} (1998), Existence results for the static contact problem with Coulomb friction. \textit{Mathematical Models and Methods in Applied Sciences} \textbf{8}, pp~445--468.
	
	\bibitem{BookJiri} \textsc{C. Eck}, \textsc{J. Jaru\v{s}ek} and \textsc{M. Krbec} (2005),  \textit{Unilateral Contact Problems in Mechanics.  Variational Methods and Existence Theorems.} Monographs \& Textbooks in Pure \& Appl.  Math.  No. 270 (ISBN 1-57444-629-0).  Chapman \& Hall/CRC (Taylor \& Francis Group), Boca Raton - London - New York - Singapore.
	
	\bibitem{FrancfortMarigo} G.A.~\textsc{Francfort} and J.J.~\textsc{Marigo} (1998), Revisiting brittle fracture as an energy minimization problem. \textit{Journal of the Mechanics and Physics of Solids} \textbf{46}, No~8, pp~1319--1342.
	
	\bibitem{JarusekPhD1} J.~\textsc{Jarušek} (1983), Contact problems with bounded friction.  Coercive  case. \textit{Czechoslovak Mathematical Journal} \textbf{33}, pp~237--261.
	
	\bibitem{JarusekPhD2} J.~\textsc{Jarušek} (1984), Contact problems with bounded friction.  Semicoercive case. \textit{Czechoslovak Mathematical Journal} \textbf{34}, pp~619--629.

	\bibitem{Klarbring} A.~\textsc{Klarbring} (1990), Examples of non-uniqueness and non-existence of solutions to quasistatic contact problems with friction. \textit{Ingenieur-Archiv} \textbf{60}, No~8, pp~529--541.

	\bibitem{Michelle} G.~\textsc{Lebeau} and M.~\textsc{Schatzman} (1984), A wave problem in a half-space with a unilateral constraint at the boundary. \textit{Journal of Differential Equations} \textbf{54}, No~3, pp~309--361.
	
	\bibitem{LerayLions} J.~\textsc{Leray} \& J.L.~\textsc{Lions} (1965), Quelques résultats de Vi\v{s}ik sur les problèmes elliptiques non linéaires par les méthodes de Minty-Browder. \textit{Bulletin de la Société Mathématique de France} \textbf{93}, pp~97--107.
	
	\bibitem{LionsMagenes} J.L.~\textsc{Lions} and E.~\textsc{Magenes} (1968), \textit{Problèmes aux Limites non Homogènes et Applications.  Volume~1}.  Dunod, Paris. 
	
	\bibitem{FirstMielke} A.~\textsc{Mielke} (2005), Evolution of rate-independent systems. In: \textit{Handbook of Differential Equations: Evolutionary Equations, Vol.~2 (C.M.~Dafermos and E.~Feireisl Editors)}, pp~461--559.  North-Holland, Amsterdam.

	\bibitem{SecondMielke} A.~\textsc{Mielke} and F.~\textsc{Theil} (2004), On rate-independent hysteresis models. \textit{Nonlinear Differential Equations and Applications} \textbf{11}, No~2, pp~151--189.

	\bibitem{Mielke} A.~\textsc{Mielke} \& T.~\textsc{Roub\'\i\v{c}ek} (2015), \textit{Rate-Independent Systems.  Theory and Application}.  Springer New York Heidelberg Dordrecht London.
	
	\bibitem{Murat} F.~\textsc{Murat} (1981) L'injection du cône positif de \(H^{-1}\) dans \(W^{-1,q}\) est compacte pour tout \(q<2\).  \textit{Journal de Mathématiques Pures et Appliquées} \textbf{60}, pp~309--322.
	
	\bibitem{Ting}  T.C.T. \textsc{Ting} (1996), \textit{Anisotropic Elasticity.  Theory and Applications.} Oxford university Press, New York and Oxford.
	
	\bibitem{Troianiello} G.M.~\textsc{Troianiello} (1987), \textit{Elliptic Differential Equations and Obstacle Problems}.  The University Series in Mathematics, Plenum Press, New York and London.
	
\end{thebibliography}
\end{document}